%% file: hs.tex
\documentclass [11pt]{amsart}
\usepackage {amsmath, amssymb, amscd, graphicx, color, url}

\newtheorem {theorem}{Theorem}[section]
\newtheorem {definition}[theorem]{Definition}
\newtheorem {lemma}[theorem]{Lemma}
\newtheorem {proposition}[theorem]{Proposition}
\newtheorem {corollary}[theorem]{Corollary}

\newtheorem {question}[theorem]{Question}
\newtheorem {remark}[theorem]{Remark}

\newtheorem {fact}[theorem]{Fact}

\def\ss {{\mathcal{S}}}
\def\zz {{\mathbb{Z}}}
\def\rr {{\mathbb{R}}}
\def\cc {{\mathbb{C}}}
\def\qq {{\mathbb{Q}}}
\def\ppp {{\mathbb{P}}}
\def\tt {{\mathbb{T}}}
\def\ze {{\mathcal{Z}}}

\def\gz {{\mathcal{G}}}
\def\oo {{\mathcal{O}}}

\def\pp {{\mathcal{P}}}
\def\ii {{\mathcal{I}}}
\def\ttt {{\mathcal {T}}}
\def\lfrak {{\mathfrak{L}}}
\def\tfrak {{\mathfrak{T}}}
\def\vv {{\operatorname{\mathbf{v}}}}
\def\ww {{\operatorname{\mathbf{w}}}}
\def\ka {{L}}
\def\ll {{\mathcal{L}}}
\def\nn {{\mathcal{N}}}
\def\rrr {{\mathcal{R}}}
\def\ss {{\mathcal{S}}}
\def\ccc {{\mathcal{C}}}
\def\dd {{\mathcal{D}}}
\def\yy {{\mathcal{Y}}}
\def\kk {{\mathcal{K}}}
\def\ve {{\mathcal{V}}}
\def\xh {{\mathrm{x}}}
\def\yh {{\mathrm{y}}}
\def\vh {{\mathrm{v}}}
\def\uh {{\mathrm{u}}}
\def\sig {{\hat{S}}}
\def\sigc {{\check{S}}}
\def\spec {{\text{Spec\ }}}

\def\hi {{\operatorname{Hilb}}}
\def\sym {{\operatorname{Sym}}}
\def\sing {{\operatorname{Sing}}}
\def\ggg {{\mathfrak{g}}}
\def\sl {{\mathfrak{sl}}}
\def\spc {{\mathfrak{s}}}
\def\gl {{\mathfrak{gl}}}
\def\conf {{\operatorname{Conf}}}
\def\resc {{\operatorname{resc}}}

\def\fin {{\hfill \square}}

\def\hom {{\operatorname{Hom}}}

\def\aalpha {{\bf \alpha}}
\def\bbeta {{\bf \beta}}
\def\tha {{\tt_{\hat{\aalpha}}}}
\def\thb {{\tt_{\hat{\bbeta}}}}
\begin{document}

\title [Nilpotent slices, Hilbert schemes, and the Jones polynomial] 
{Nilpotent slices, Hilbert schemes, and the Jones polynomial} \author 
[Ciprian Manolescu]{Ciprian Manolescu} 
\thanks {The author is supported by a Clay Research Fellowship.}

\begin {abstract} Seidel and Smith have constructed an invariant of links
as the Floer cohomology for two Lagrangians inside a complex affine
variety $Y.$ This variety is the intersection of a semisimple orbit with a
transverse slice at a nilpotent in the Lie algebra $\sl_{2m}.$ We exhibit
bijections between a set of generators for the Seidel-Smith cochain
complex, the generators in Bigelow's picture of the Jones polynomial, and 
the generators of the Heegaard Floer cochain complex for the double branched 
cover. This is done by presenting $Y$ as an open subset of the Hilbert scheme 
of a Milnor fiber.  \end {abstract}

\address {Department of Mathematics, Princeton University\\ Princeton, NJ 
08540}
\email {cmanoles@math.princeton.edu}

\maketitle

\section {Introduction}

Khovanov cohomology \cite{Kh} is an invariant of links in the form of a 
bigraded 
abelian group $Kh^{*, *}(\ka)$ whose graded Euler characteristic is the 
unnormalized Jones polynomial of the link $\ka:$
$$ \sum_{i,j \in \zz} (-1)^{i+j+1} t^{j/2} \dim(Kh^{i,j}(\ka) \otimes \qq) 
= (t^{1/2}+t^{-1/2}) V_{\ka}(t).$$

By definition, the groups $Kh^{*,*}$ can be combinatorially computed starting from
a specific diagram of the link. Nevertheless, they have been found to be quite
powerful and to have much in common with more subtle invariants of links coming
from gauge theory and symplectic geometry. For example, Rasmussen \cite{R} has used
Khovanov's theory to give a new proof of Milnor's conjecture on the slice genus of
torus knots, a result proved previously by Kronheimer and Mrowka using gauge
theory. Also, in \cite{OS3}, Ozsv\'ath and Szab\'o have constructed a 
spectral
sequence relating $Kh^{*,*}(\ka)$ to the Heegaard Floer homology of the double 
cover of $S^3$ branched over $\ka.$ Heegaard Floer homology is a version of the 
Floer homology defined in symplectic geometry by counting pseudoholomorphic 
curves.

Seidel and Smith \cite{SS} have proposed a remarkable interpretation of
Khovanov cohomology itself in terms of symplectic geometry. They represent a
link $\ka$ as the closure of an $m$-stranded braid $b \in Br_m$ or,
equivalently, as the plat closure of $b \times 1^m \in Br_{2m}.$ Let 
$w$ be the writhe of the braid diagram. $b \times 1^m$ can
be represented as a loop $l$ in the configuration space $\conf^{2m}(\cc)$ 
of $2m$ distinct points in the plane, with $l$ 
starting at some basepoint $\tau.$ They construct a symplectic fibration
over $\conf^{2m}(\cc),$ whose fiber at $\tau$ is a symplectic manifold $Y 
=\yy_{m, \tau},$ and then introduce a Lagrangian $\ll \subset Y,$ 
well-defined up to
isotopy. Applying to $L$ the monodromy map along $l$ yields another
Lagrangian $\ll' \subset Y.$ They define the bigraded groups (the
``symplectic Khovanov cohomology'') $$ Kh_{symp}^*(\ka) = HF^{*+m+w}(\ll,
\ll') $$ as the Lagrangian Floer cohomology applied to $\ll$ and $\ll'$
and with a shift in degree.

Seidel and Smith prove that $Kh_{symp}$ is an invariant of the link $\ka,$ and 
conjecture that it equals the
original Khovanov cohomology, after a collapsing of its bigrading:
$$ Kh_{symp}^k(\ka) \cong \bigoplus_{i-j=k} Kh^{i,j}(\ka) \hskip14pt (?) 
$$

The Seidel-Smith cohomology does not come with a bigrading from which one can read
off the Jones polynomial. The goal of the present paper is to define such a
bigrading at the level of the cochain complex, as well as to shed some light on the
rather mysterious Seidel-Smith construction by giving a more concrete 
characterization of the objects involved.

We start by describing the symplectic manifold $Y$ as an open subset of
a Hilbert scheme. $Y = \yy_{m, \tau}$ is in fact an affine variety,
defined as the set of matrices in $\ss_m$ with fixed characteristic
polynomial $P_{\tau}(t).$ Here $\ss_m$ is an affine subspace of the Lie 
algebra
$\sl_{2m}$ transverse to the orbit of a nilpotent element $N_m$ with two
Jordan blocks of size $m,$ and the coefficients of $P$ are described by
$\tau \in \conf^{2m}(\cc).$ More generally, instead of $N_m$ we can 
consider a nilpotent with two Jordan blocks of sizes $n$ and $2m-n,$ 
respectively ($n \leq m),$ and similarly construct an affine variety 
$\yy_{n, \tau}.$ Note that $m$ is implicit in this notation, $\tau$ being 
an unordered set of $2m$ points. We prove:
\begin {theorem}
\label {hilb}
There is an injective holomorphic map from $\yy_{n, \tau}$ to the  
Hilbert scheme $\hi^n(S_{\tau}),$ where $S_{\tau}$ is the affine surface 
described by the equation $u^2 + v^2 + P_{\tau}(z)=0$ in $\cc^3.$ 
\end {theorem}
Since $\yy_{n, \tau}$ and the Hilbert scheme have the same complex 
dimension $2n,$ it follows that the former can be identified to an open 
subset of the second. 

The case $n=m$ is the only one relevant for the construction of the
Seidel-Smith homology. Nevertheless, Theorem~\ref{hilb} describes a
phenomenon that can be interesting by itself. Indeed, both $\yy_{n, \tau}$
and $\hi^n(S_{\tau})$ are quiver varieties in the sense of Nakajima
\cite{N1}. One may ask whether there are other pairs of quiver varieties
of the same dimension such that one is an open subset of the other.

Next, we give an explicit description of two Lagrangians that can be used to define the
symplectic Khovanov cohomology. This is possible because the Hilbert scheme
$\hi^m(S_{\tau})$ is a certain iterated blow-up of the symmetric product $\sym^m(S_{\tau})$
along subsets of the diagonal $\Delta.$ Thus we can identify $\sym^m(S_{\tau}) - \Delta$ to
an open subset in $\hi^m(S_{\tau}),$ and on $\sym^m(S_{\tau})$ one has nice holomorphic
coordinates. Indeed, a point in $\sym^m(S_{\tau}) - \Delta$ is characterized as an
unordered collection of $m$ distinct points $(u_k, v_k, z_k)  \in S_{\tau}, k=1,\dots, m.$

Let us choose $m$ disjoint arcs $\alpha_k: [0,1] \to \cc,$ joining 
together the $2m$ points of $\tau$ in pairs. The braid $b \times 1^m \in Br_{2m},$ 
whose plat closure is the link $L,$ induces a diffeomorphism of the plane 
that maps the arcs $\alpha_k$ into $2m$ arcs $\beta_k: [0,1] \to \cc,$ 
again joining the points of $\tau.$ For each $\alpha_k,$ one can construct 
the Lagrangian 2-sphere:
$$ \Sigma_{\alpha_k} = \{ (u, v, z) \in S_{\tau} \hskip3pt : z = 
\alpha_k(t) \text{ 
for some }t \in [0,1]; \hskip5pt u, v \in \sqrt{-P_{\tau}(z)} \rr \}$$
in $S_{\tau}$ with the standard K\"ahler metric. Of course, the same 
construction can be done for the beta curves, resulting in spheres 
$\Sigma_{\beta_k}.$

\begin {theorem}
\label {lagr}
We can deform the K\"ahler metric on $\yy_{m, \tau}$ and the Lagrangians $\ll, \ll'$ in the 
Seidel-Smith construction in such a way that the Floer cohomology groups are preserved 
under this deformation, and the resulting Lagrangians are 
$$ \kk = \Sigma_{\alpha_1} \times \Sigma_{\alpha_2} \times \dots 
\times \Sigma_{\alpha_m} \ ; \ 
 \kk' = \Sigma_{\beta_1} \times \Sigma_{\beta_2} \times \dots
\times \Sigma_{\beta_m} $$
in $\bigl (\sym^m(S_{\tau}) - \Delta\bigr) \cap \yy_{m, \tau} \subset 
\hi^m(S_{\tau}).$
\end {theorem}

This enables us to describe explicitly the 
intersection $\kk \cap \kk'$ 
and, in turn, this gives a set of generators for the Seidel-Smith 
cohomology. More precisely, we form a set $\ze$ in the following way: we 
assume that the $\alpha$ and $\beta$ are simple curves that intersect 
transversely in their interior. For every intersection point $x \in \alpha_i 
\cap \beta_j,$ we introduce an element $e_x \in \ze$ in the case when $x 
\in 
\tau,$ and two elements $e_x, e'_x$ when $x \not \in \tau.$ We define maps 
$A, B : \ze \to \{1,2, \dots, m\}$ by taking an element to the indices 
$i,j$ of 
the corresponding $\alpha$ and $\beta$ curves. 

At this point, we note that this picture is very similar to the one given 
by Bigelow \cite{B} in his definition of the Jones polynomial. By 
presenting a link as the plat closure of a braid, Bigelow obtains its 
Jones polynomial as a signed count of the intersection points of two 
half-dimensional submanifolds in a covering of a subset of $\sym^m(\cc)$, 
with certain gradings. 

We exploit this similarity and show:
\begin {theorem}
\label {bijection}
There is a natural correspondence between a set of generators for 
$CF^*(\kk, \kk')$ and the set $\gz$ of $m$-tuples $(z_1, \dots, z_m)$ of 
unordered elements of $\ze$ with
$A(z_i) \neq A(z_j)$ and $B(z_i) \neq B(z_j)$ for $i\neq j.$ This set can 
also be identified with the set of intersection points in Bigelow's 
picture of the Jones polynomial. 
\end {theorem}

Bigelow defined two gradings $Q, T: \gz \to \zz$ whose difference doubled and shifted by a
constant is the ``Jones grading'' $J = 2(T - Q) + m +w.$ The grading $J$ of an element in
$\gz$ tells the coefficient in the Jones polynomial to which that generator contributes in
the signed count.  We define a third grading, $P: \gz \to \zz$ (the ``projective
grading''), starting from the Maslov grading $\tilde P$ of the generator in $CF^*(\kk,
\kk'),$ and then normalizing to $P = \tilde P - (m+w).$ We also show how to read $P$ from 
the concrete picture of the alpha and beta curves intersecting in the plane.

Our discussion can be related to the work of Ozsv\'ath and Szab\'o \cite{OS3}, following
the ideas of Seidel and Smith \cite{S}. The involution $(u, v, z) \to (u, -v, z)$ on
$S_{\tau}$ induces a corresponding involution on $\hi^m(S_{\tau}).$ This involution
preserves the Lagrangians $\kk, \kk'$ from Theorem~\ref{lagr}, and looking at its fixed
point set we find two totally real tori $\tha, \thb$ sitting inside the symmetric power
of a Riemann surface. The Floer homology $HF_*(\tha, \thb)$ is exactly the Heegaard
Floer homology of $\dd(L) \# (S^1 \times S^2),$ where $\dd(L)$ is the double branched
cover of $S^3$ over $L.$ We can also consider the restriction of the involution to the
variety $\yy_{m, \tau}.$ Its fixed point set $W$ is an open subset of the symmetric
product, and the complement of $W$ is a codimension one subvariety $\nabla$ called the
``anti-diagonal.'' This point of view is useful because even though the homological grading on the
Heegaard Floer chain complex in the symmetric product is only defined modulo $2,$ it can
be improved to a $\zz$ grading by considering the restriction to $W.$

\begin {theorem}
\label {newth}
After some noncanonical choices (to be described in Section~\ref{sec:dbc}), the set $\gz$ 
can 
be identified with a set of generators for the Heegaard Floer homology
of $\dd(L) \# (S^1 \times S^2).$ There is a well-defined integer grading on $\gz$ induced 
by the Maslov grading of the intersection points of $\tha$ and $\thb$ taken 
inside the open set $W.$ This grading equals $\tilde P + T - Q.$
\end {theorem}

Denote the Maslov grading from Theorem~\ref{newth} by $\tilde R = \tilde P + T - Q.$ If we 
renormalize to $R = \tilde R -
(m+w)/2,$ we obtain the simple formula $R = P + (J/2).$ This renormalization is similar
to that of $P = \tilde P - (m+w).$ In the case of the identity braid $b = 1^m \in Br_m,$
they make the $P$ and $R$ gradings of the resulting generators to be symmetric around
zero.  The fact that the correction term for $R$ is half of that for $P$ is related to
the fact that restriction to the fixed point set cuts in half the dimension of the
objects involved.

Let us end with some open questions. First, Theorem~\ref{bijection} does not
say anything about the differentials in the Seidel-Smith cochain complex.
We conjecture that these differentials preserve the $J$ grading, and that
therefore the bigrading $(J, P)$ descends to the level of cohomology. This
would imply the existence of a bigraded Floer theory analagous to
Khovanov's. Some hope in the direction of this conjecture could come from
looking at the localization of Floer cohomology under the involution
considered above. As suggested by Seidel and Smith, localization could also give a 
geometric interpretation of the spectral sequence in \cite{OS3}. While we cannot prove the 
conjecture at this point, a quick corollary of Theorem~\ref{newth} is that the 
difference of the Maslov gradings in $\yy_{m, \tau}$ and in the
fixed point set $W,$ after some normalization, gives exactly the Jones 
polynomial i.e. the (conjectural) graded Euler characteristic for $Kh_{symp}:$ 

\begin {corollary}
\label {coro}
Given a bridge presentation of a link $L,$ consider the corresponding 
set of generators $\gz$ with the two gradings $P$ and $R$ as above. 
Then the Jones polynomial of $L$ can be expressed in terms of symplectic geometric 
quantities by the formula:

$$ V_L(t) = -(t^{1/2} + t^{-1/2})^{-1} \cdot \sum_{\gamma \in \gz} 
(-1)^{P({\gamma})}  t^{R(\gamma) - P(\gamma)}.$$

\end {corollary}

Second, the Lagrangians $\kk$ and $\kk'$ are subsets of the Hilbert scheme
$\hi^m(S_{\tau}).$ The possible relevance of $\hi^m(S_{\tau})$ to
low-dimensional topology was pointed out by Khovanov in \cite[section
6.5]{Kh2}. One can equip the Hilbert scheme with a suitable K\"ahler
metric and look at the Floer cohomology of $\kk, \kk'$ inside
$\hi^m(S_{\tau}).$ An interesting question is whether this cohomology
turns out to be a link invariant.

Finally, following \cite{S}, we can also look at the Floer homology of the tori $\tha, \thb$ inside
$W.$ We conjecture that this homology, together with its $R$ grading, is another link invariant.

\medskip \noindent \textbf{Acknowledgements.} I would like to thank Peter Kronheimer for his guidance and 
encouragement, Jacob Rasmussen for helping me become familiar with both Khovanov homology and Bigelow's 
work, and Mikhail Khovanov, Gang Tian and Davesh Maulik for helpful conversations. I am grateful to Cliff 
Taubes and to the participants in his informal Harvard seminar for creating a forum for discussing the 
Seidel-Smith construction in the summer of 2004. I am also indebted to the referee for many useful 
comments and for suggesting an improvement in Proposition 2.7.

\section {Nilpotent slices and Hilbert schemes}

In this section we recall the definitions of the manifold $\yy_{m, \tau}$
appearing in the Seidel-Smith construction, of the Hilbert scheme 
$\hi^n(S_{\tau}),$ and then we prove Theorem 1. 

Let us introduce some notation. For $m > 0,$ $\sym^{2m}(\cc)$ is the 
symmetric product 
of $\cc$ and can be identified with $\cc^{2m}$ via symmetric polynomials. 
$\sym^{2m}_0(\cc) \cong \cc^{2m-1}$ is the space of unordered sets of $2m$ 
complex numbers with sum zero.

The configuration spaces
$$ \conf^{2m}(\cc) \subset \sym^{2m}(\cc) \text{ ; }  \conf^{2m}_0(\cc) 
\subset \sym^{2m}_0(\cc) $$
consist of the $2m$-tuples made of distinct complex numbers. The map
\begin {eqnarray}
\label {proj}
 p: \sym^{2m} (\cc) &\to& \sym^{2m}_0 (\cc), \\
 (a_1, \dots, a_{2m}) &\to& (a_1 - \frac{1}{2m}\sum a_k, \dots, a_{2m} -  
\frac{1}{2m}\sum a_k)  \nonumber
\end {eqnarray}
is a trivial $\cc$-bundle, and its restriction to $\conf^{2m}(\cc)$ also 
exhibits this space as a trivial $\cc$-bundle over $\conf^{2m}_0(\cc).$

Throughout this section, we restrict our attention to $\tau \in 
\sym^{2m}_0(\cc),$ with the understanding that everything extends 
trivially for $\tau \in \sym^{2m}(\cc)$ by considering the objects 
associated to $p(\tau).$

\subsection {Transverse slices at nilpotent orbits}
\label {sec:yy}
We start by presenting the definition of the manifold $\yy_{m, \tau}$ 
from \cite{SS}. As noted in the introduction, 
we work in a slightly more general setting and define a manifold $\yy_{n, 
\tau}$ for any $0 \leq n \leq m$ and $\tau \in \conf^{2m}_0(\cc).$ In 
fact, we can do the same for every $\tau \in \sym^{2m}_0(\cc),$ but in 
that case the result is a possibly singular variety $\yy_{n, \tau}.$ 

Consider the complex algebraic group $G=SL_{2m}$ and its Lie algebra 
$\ggg 
= \sl_{2m}.$ The adjoint quotient map 
$$\chi: \ggg \to \ggg/G = \cc^{2m-1}$$
is defined by taking $A \in \ggg$ to the coefficients of the 
characteristic 
polynomial $det(tI-A).$ The fiber of $\chi$ over $\tau \in 
\conf^{2m}_0(\cc) \subset \sym^{2m}_0(\cc) = \cc^{2m-1}$ is a smooth 
manifold (the adjoint orbit of a semisimple element in $\ggg$). In 
fact, the map $\chi$ can be shown to be a differentiable fiber bundle when 
restricted to the preimage of $\conf^{2m}_0(\cc).$

Choose a nilpotent element $N_n \in \ggg$ with two Jordan 
blocks of sizes $n$ and $2m-n$, respectively:
$$ N_n =\left( \begin {tabular}{c|c}
$ \begin {array}{ccccc} 0 & 1 &  &  & \\
 & 0 & 1 & \ & \\
 & \ & \ & \ldots & \ \\
 & \ & \  & \  &  1 \\
 & \ & \ & \ & 0
\end {array}$  &

 \\
\hline
                                                                                
  &
                                                                                
$ \begin {array}{cccccc} 0 & 1 &  & &  & \\
 & 0 & 1 & \ & & \\
 & \ & \ & \ldots & \ & \\
 & \ & \ & \ldots & \ & \\
 & \ & \ &  & \  &  1 \\
 & \ & \ & & \ & 0
\end {array} $
                                                                                
\end {tabular} \right).
$$

The orbit of $N_n$ under the adjoint action of $G$ is a manifold $\oo_n$ 
whose tangent space at $N_n$ can be described as:
$$ T_{N_n}\oo_n = N_n + ad(\ggg)N_n = \{ N_n + [N_n, B] : B \in \ggg\}.$$

\begin {definition}
A {\bf transverse slice} at $N_n \in \ggg$ is a local complex submanifold 
$\ss \subset \ggg, N_n \in \ss,$ such that the tangent spaces of 
$\ss$ and $\oo_n$ at $N_n$ are complementary.
\end {definition}

In our discussion we choose a particular slice $\ss_n,$ the affine 
subspace consisting of matrices of the form:
$$ A =\left( \begin {tabular}{l|l}
$ \begin {array}{lcccc} a_1 & 1 &  &  & \\
a_2 & 0 & 1 & \ & \\
\ldots & \ & \ & \ldots & \ \\
a_{n-1} & \ & \  & \  &  1 \\
a_n & \ & \ & \ & 0
\end {array}$  & 

$\begin {array}{lccccc} b_1 &  & & &  & \\
b_2 &  &  &  & &\\
\ldots &  &  &  &  & \ \\
b_{n-1} & \ & \ &  & \  &   \\
b_n & \ & \ & \ & &
\end {array}  $

 \\
\hline 
                                                                                
$ \begin {array}{lcccc} 0 &  &  &  & \\
\ldots &  &  &  & \\
0 &  &  &  & \ \\
c_1 & \ &  & \  &   \\
\ldots & & & & \\
c_n & \ & \ & \ &
\end {array}$ &

$ \begin {array}{lccccc} d_1 & 1 &  & &  & \\
d_2 & 0 & 1 & \ & & \\
\ldots & \ & \ & \ldots & \ & \\
\ldots & \ & \ & \ldots & \ & \\
 & \ & \ &  & \  &  1 \\
d_{2m-n} & \ & \ & & \ & 0
\end {array} $

\end {tabular} \right),
$$
where $a_k, b_k, c_k, d_k \in \cc, a_1 + d_1 =0.$ 

Consider the polynomials 
\begin {eqnarray*}
A(t) &=& t^n -a_1 t^{n-1} + a_2 t^{n-2} - \dots + (-1)^n a_n; \\
B(t) &=& b_1t^{n-1} - b_2t^{n-2} + \dots +(-1)^{n-1}b_n ;\\
C(t) &=& c_1t^{n-1} - c_2t^{n-2} + \dots +(-1)^{n-1}c_n;\\
D(t) &=& t^{2m-n} -d_1 t^{2m-n-1} +  \dots + (-1)^{2m-n} 
d_{2m-n}.
\end {eqnarray*}

It is easy to check that:
$$ det(tI-A) = A(t)D(t) - B(t)C(t).$$

Also, a straightforward computation shows that $\ss_n$ is complementary to
$T\oo_n$ at $N_n,$ and thus is indeed a transverse slice.  In the case
$n=m,$ $\ss_m$ is the slice considered by Seidel and Smith in
\cite{SS}, with a reordering of the coordinates.

Let us look at the restriction of the adjoint map to the slice:
$$ \chi|_{\ss_n} : \ss_n \to \sym^{2m}_0(\cc).$$
This is again a differentiable fiber bundle when restricted over 
$\conf^{2m}_0(\cc),$ and its fibers are affine varieties of complex 
dimension $2n.$

For every $\tau \in \sym^{2m}_0(\cc),$ we define:
$$ \yy_{n, \tau} = \chi|_{\ss_n}^{-1}(\tau).$$

Explicitly, $\yy_{n, \tau}$ is an affine variety in $\ss_n = 
\cc^{2m+2n-1}.$ In terms of the coordinates $a_k, b_k, c_k, d_k,$
it is described by a set of $2m-1$ algebraic equations that can all 
be grouped into one:
\begin {equation}
\label {eq}
A(t)D(t) - B(t)C(t) = P_{\tau}(t).
\end {equation}
Here $P_{\tau}$ is the polynomial with roots given by $\tau \in 
\conf_0^{2m}(\cc),$ counted with multiplicities, and the equations 
correspond to identifying the $2m-1$ lowest coefficients of $t$ in 
(\ref{eq}).

\subsection {An important example.} 
\label {sec:ex}
Let $n=1.$ Then $\ss_n = \cc^{2m+1}$ with coordinates $z=a_1 = -d_1, 
b_1, c_1, d_2, \dots, d_{2m-1}.$ Our manifold is described as:
$$ \yy_{1, \tau} = \{(z, b_1, c_1, d_2, \dots, d_{2m-1}) : (t-z)D(t) - 
b_1c_1 = P_{\tau}(t)\}.$$

Note that once we know that $P_{\tau}(z) + b_1c_1 = 0,$ the polynomial 
$D(t)$ is determined uniquely. Hence, by making the coordinate change $b_1 
= u+iv, c_1 = u-iv,$ we can write:
$$ \yy_{1, \tau} = \{ (u, v, z) \in \cc^3 : u^2 + v^2 + P_{\tau} (z) = 0\} 
$$

This is a complex surface, the Milnor fiber associated to the $A_{2m}$ 
singularity, and it will play an important role in the discussion to 
follow. We denote it by $S_{\tau}.$ It is smooth if and only if the 
polynomial $P_{\tau}$ has no multiple roots, i.e. for $\tau \in 
\conf_0^{2m}(\cc).$

\subsection {The Hilbert scheme.} 
\label {sec:hi}

We recall a few standard results about Hilbert schemes of points on 
surfaces. The interested reader can consult the book by Nakajima \cite{N2}
for a thorough exposition of the subject.

Let $X$ be a quasi-projective scheme over an algebraically closed field 
$k.$ Fix a polynomial $P \in \zz[t].$ A fundamental result of Grothendieck 
\cite{G} is the existence of a quasi-projective scheme $\hi^P(X)$ that  
parametrizes flat families of closed subschemes of $X$ with Hilbert 
polynomial $P.$ In other words, $\hi^P(X)$ comes with a universal family 
$\mathcal{Z}$ such that for every flat family $Z$ of closed subschemes of 
$X$ with Hilbert polynomial $P,$ parametrized by a scheme $U,$ there is a 
unique morphism $\phi: U \to \hi^P(X)$ such that $Z$ is the pullback 
$\phi^*(\mathcal{Z}).$  

For our purposes we restrict to the case when $k=\cc$ and $P$ is the 
constant polynomial $n.$ In this case $\hi^n(X)$ parametrizes closed 
$0$-dimensional subschemes of $X$ of length $n.$ The typical example of 
a subscheme of this form is a subvariety consisting of $n$ distinct 
points of $n.$ We get nonreduced examples by letting some of the $n$ 
points collide.

\begin {definition}
$\hi^n(X)$ is called the Hilbert scheme of $n$ points on $X.$
\end {definition}

Here are a few basic facts about Hilbert schemes of points, collected from 
\cite{N2}. We assume for simplicity that $X$ is reduced.

\begin{fact}
There is a natural morphism from the Hilbert scheme to the 
symmetric product of $X,$ defined by
$$ \pi: \hi^n(X) \to \sym^n(X),$$ 
\begin {equation}
\label {hc}
 \pi (Z) = \sum_{x \in X} \text{length}(Z_x) [x]. 
\end {equation}
This is called the Hilbert-Chow morphism.
\end {fact} 

\begin {fact}
\label {fact1}
When $\dim_{\cc} X =1,$ $\pi$ is an isomorphism. Therefore, $\hi^n(X) = 
\sym^n(X).$
\end {fact}

\begin {fact} [Fogarty's theorem \cite{Fo}]
\label {facto}
When $\dim_{\cc} X = 2$ and $X$ is smooth, $\hi^n(X)$ is smooth of 
complex dimension $2n,$ and the Hilbert-Chow morphism is a resolution of 
singularities. If $X$ is irreducible, then $\hi^n(X)$ is also 
irreducible.
\end {fact}

In higher dimensions there are examples when the Hilbert scheme is not 
smooth even for $X$ smooth.

We are interested in the case when $\dim_{\cc} X =2,$ so let us explain 
more carefully what happens then. The symmetric product $\sym^n(X)$ is singular 
along its diagonal:
$$ \Delta = \{(x_1, \dots, x_n) \in \sym^n(X) : x_i = x_j \text{ for 
some } i\neq j\}.$$
The diagonal has complex codimension $2$ in $\sym^n(X),$ and the 
Hilbert-Chow 
morphism $\pi: \hi^n(X)$ $\to \sym^n(X)$ is one-to-one over $\sym^n(X) - 
\Delta.$ 

\begin {fact}
When $\dim X = 2,$ the preimage $\pi^{-1}(\Delta)$ has complex codimension 
$1$ in $\hi^n(X).$
\end {fact}

This is easiest to visualize when $n=2$ and $X$ is smooth. The Hilbert
scheme $\hi^2(X)$ parametrizes $0$-dimensional subschemes of length $2.$
The reduced ones are of the form $(x_1, x_2) \in \sym^2(X) - \Delta.$ The
nonreduced ones are in the preimage $\pi^{-1}(\Delta).$ A point $z$ which 
maps 
to $(x, x)$ under $\pi$ is a subscheme defined by $\oo_X / \ii,$ where 
$\ii \subset \oo_X$ is an ideal of the form
$$ \ii = \{f \in \oo_X : f(x) =0, df_x(v) =0\},$$
for $v \neq 0 \in  T_x X.$ In other words, points in $\pi^{-1}(\Delta)$ 
are described by a point $x \in X$ and a direction $v$ at $X.$ The fiber 
of $\pi$ over $(x, x) \in \Delta$ is a copy of $\ppp^1$ and, in fact, the 
Hilbert scheme is a blow-up of the symmetric product: 
$$ \hi^2(X) = \text{Blow}_{\Delta}(\sym^2(X)) = \text{Blow}_{\Delta}(X 
\times X) / \zz_2.$$

In general, $\hi^n(X)$ is more difficult to describe. We refer to
\cite{N2} for an explicit stratification in the case $X= \cc^2.$ This is
relevant for the local behaviour of $\hi^n(X)$ for any surface $X.$

\subsection {An open holomorphic embedding.}
\label {sec:open}

We are now ready to prove Theorem~1.1, which says that there is an 
injective holomorphic map $j$ from $\yy_{n, \tau}$ into the Hilbert scheme 
of $n$ points on the affine surface $S_{\tau}$ considered in 
section~\ref{sec:ex}. We construct $j$ as an algebraic morphism. 

By the defining property of Hilbert schemes, a morphism into 
$\hi^n(S_{\tau})$ is the same as a flat family of subschemes of $S_{\tau}$ 
with Hilbert polynomial $n,$ parametrized by the domain $\yy_{n, \tau}.$ 

Recall from section~\ref{sec:yy} that $\yy_{n, \tau} = \spec R,$ where 
$R$ is a commutative ring, the quotient of the polynomial ring in 
the $2m+2n - 1$ coordinates $a_k, b_k, c_k, d_k$ by the ideal generated by 
the algebraic relations in (\ref{eq}). We can think of $A(t) = t^n 
- a_1t^{n-1} + \dots + (-1)^n a_n$ as an
element in $R[t],$ and the same is true for $B(t), C(t), D(t).$
Consider the polynomials
\begin {equation}
\label {bcuv}
 U(t) = \frac{1}{2}\bigl( B(t) + C(t)\bigr) \in R[t]; \hskip6pt V(t) = 
\frac{1}{2i}\bigl( B(t) - C(t)\bigr) \in R[t].
\end {equation}

Then (\ref{eq}) can be rewritten as:
\begin {equation}
\label{neweq}
U(t)^2 + V(t)^2 + P_{\tau}(t) = A(t)D(t). 
\end {equation}
Set
$$ \rrr = R[u, v, z] / (u^2 + v^2 + P_{\tau}(z)). $$
 
Consider the map 
$$ \psi: \rrr \to R[t], \hskip10pt \psi(Q(u,v, z)) = Q(U(t), V(t), t)$$ 
and let us compose it with the natural projection $p: R[t] \to 
R[t]/(A(t)).$

It is easy to see that the composition 
$$p \circ \psi: \rrr \to R[t]/(A(t))$$
is surjective. Let $\ii \subset \rrr$ be its kernel. Then $\rrr/\ii$ is 
isomorphic to $R[t]/(A(t)) \cong R^n$ as an $R$-module.

We define the closed subscheme 
$$ Z = \spec \rrr/\ii \subset \spec \rrr = \yy_{n, \tau} \times S_{\tau}. 
$$

Since $\rrr/\ii$ is a free $n$-dimensional module over $R,$ it follows 
that
the composition $$ Z \subset \yy_{n, \tau} \times S_{\tau} \to \yy_{n,
\tau} = \spec R$$ 
exhibits $Z$ as a flat family of 0-dimensional subschemes of
$S_{\tau}$ of length $n.$ This defines the desired morphism
$$ j: \yy_{n, \tau} \to \hi^n (S_{\tau}).$$ 

From now on we think of $\yy_{n, \tau}$ as an affine variety, with its
reduced scheme structure. Then the points in $\yy_{n, \tau}$ are 4-tuples
of polynomials $(A, D, U, V)$ in $\cc[t]$ satisfying (\ref{neweq}), and
the points in $\hi^n(S_{\tau})$ can be identified with ideals in 
$\oo= \cc[u, v, z]/(u^2 + v^2 + P_{\tau}(z))$ such that $\dim_{\cc} (\oo/I) =n.$
Explicitly, the morphism $j$ is given by: 
\begin {equation}
\label {jaduv}
j(A, D, U, V) = \{ Q(u, v, z): A(t) \text{ divides } Q(U(t), V(t), t)  \}.
\end {equation}

Note that $R_0 = \cc[z]$ is a subring of $\cc[u, v, z],$ so $R_1 = R_0/(R_0 \cap (u^2 + v^2 + 
P_{\tau}(z))) \cong \cc[z]$ is a subring of $\oo = \cc[u, v, z]/(u^2 + v^2 + P_{\tau}(z)),$ the ring of 
functions on the affine variety $S_{\tau}.$ Given an ideal $I \subset \oo$ describing a subscheme $X = 
\spec \oo/I$ in $\hi^n(S_{\tau}),$ the 
intersection $I \cap R_1$ corresponds to a subscheme of $\cc,$ the image of $X$ under the map:
\begin {equation}
i: S_{\tau} \longrightarrow \cc, \ i(u, v, z) = z.
\end {equation}
It is clear that $R_1/(I \cap R_1)$ injects into $\oo/I,$ hence $i(X)$ must have length at most $n.$

With this background in place, we prove the following:

\begin{proposition}
The morphism $j$ is an open embedding. The image of $j$ in $\hi^n(S_{\tau})$ consists of the subschemes 
$X$ such that $i(X)$ is a subscheme of $\cc$ of length exactly $n.$   
\end {proposition}

\noindent{\it Proof.} Since $\yy_{n, \tau}$ and $\hi^n(S_{\tau})$ have 
the same dimension and $\hi^n(S_{\tau})$ is 
irreducible (Fact~\ref{facto}), in order to prove that $j$ is an open 
embedding it suffices to show that it is 
injective. Pick $(A_i, D_i, U_i, V_i), i=1,2,$ that map to the same ideal 
$I \subset \oo $ under $j.$ Then
$$ I \cap R_1 = \{ Q \in \cc[z]: A_i(z) \text{ divides } Q(z) \},$$
where $i$ can be either $1$ or $2.$ Since the $A_i$ are monic polynomials 
of degree $n < 2m=$ degree $(P_{\tau}),$ it follows that $A_1 = A_2.$

Next, note that the polynomials $u - U_1(z)$ and $u-U_2(z)$ are in $I,$ 
hence $U_1(z) - U_2(z) \in I \cap R_1.$ But $U_1 - U_2$ has degree at most $n-1 < n 
=$ degree $A,$ so $A$ dividing $U_1 - U_2$ modulo $P_{\tau}$ implies $U_1 
= U_2.$ Similarly $V_1 = V_2.$ Also, the relation (\ref{neweq}) determines $D$ uniquely from $A, U, V$ 
and $P_{\tau}.$ Therefore, we must also have $D_1 = D_2,$ and this shows 
that $j$ is injective. 

If $I$ is in the image of $j,$ then $I \cap R_1$ is of the form $\{ Q \in \cc[z]: A(z) \text{ divides } 
Q(z) \},$ where $A$ is a monic polynomial of degree $n.$ It follows that $\dim R_1 / (I \cap R_1) = n.$

Conversely, let $I$ be an ideal corresponding to a subscheme in $\hi^n(S_{\tau})$ such that $\dim R_1 / 
(I \cap R_1) = n.$ We claim that $I$ lies in the image of the embedding $j.$ Since $R_1 \cong \cc[z],$ 
the ideal $I \cap R_1$ is necessarily generated by a unique monic polynomial $A(z)$ of degree $n.$
Because of the equality of dimensions, the inclusion $R_1/(I \cap R_1) \hookrightarrow \oo/I$ is a 
bijection. Consider the projection
$$ \oo = \cc[u, v, z]/(u^2 + v^2 + P_{\tau}(z)) \  \to \oo/I \cong R_1/(I \cap R_1) \cong 
\cc[z]/(A(z)).$$
 
The images of the elements $u, v \in \oo$  under this projection can be represented by polynomials $U(z), 
V(z)$ of degree at most $n-1,$ while the image of $z$ is obviously the polynomial $z.$ Hence, we get a 
relation 
$$ U(z)^2 + V(z)^2 + P_{\tau}(z) = A(z)D(z),$$
where $D(z)$ is a uniquely determined monic polynomial of degree $2m-n.$  Because of the form of 
$P_{\tau},$ the second leading coefficients of $A$ and $D$ must add up to zero. We get that $I = j(A, D, 
U, V),$ so indeed $I$ lies in the image of $j. \hfill \fin$
\medskip

This completes the proof of Theorem~1.1. From now on we can think of 
$\yy_{n, \tau}$ as an open subset of the Hilbert scheme. Recall from 
section~\ref{sec:hi} that the Hilbert-Chow morphism 
$$ \pi: \hi^n(S_{\tau}) \to \sym^n(S_{\tau}) $$
is 1-to-1 away from the diagonal $\Delta.$ We can easily write down the 
composition:
$$ \pi \circ j : \yy_{n, \tau} \to \sym^n(S_{\tau}). $$

This leads to probably the simplest way to think about the effect of $j$:

\begin {remark}
Given a point $(A, D, U, V)$ in $ \yy_{n, \tau},$ its image under $\pi 
\circ j$ is the unordered collection of $n$ points $(u_k, v_k, z_k) \in  
\sym^n(S_{\tau}),$ $k=1, \dots, n,$ where $z_k$ are the roots of $A(t),$ 
$u_k = U(z_k)$ and $v_k = V(z_k).$ 
\end {remark}

A quick corollary of this discussion is:
\begin {corollary}
\label {cor}
The intersection of $\yy_{n, \tau}$ with the open subset 
$$U_{n, \tau} = \pi^{-1}(\sym^n(S_{\tau}) - \Delta) \subset 
\hi^n(S_{\tau})$$
is the complement in $U_{n, \tau}$ of the codimension one subset
$$ \Xi = \pi^{-1} \bigl(\{(u_k, v_k, z_k) \in \sym^n(S_{\tau}) - \Delta 
: z_i = z_j \text{ for some } i \neq j \} \bigr). $$
\end {corollary} 

\begin {remark}
When $n=1,$ we have $\yy_{1, \tau} = \hi^1(S_{\tau}) = S_{\tau},$ and $j$
is an isomorphism.
\end {remark}

\subsection {Families.} Until now we have described the open embedding $j$ 
for a fixed $\tau.$ As we vary $\tau$ in $\sym^{2m}_0(\cc),$ the varieties 
$\yy_{n, \tau}$ form the transverse slice $\ss_n.$ When $n=1,$ the 
family $\yy_{n, \tau} = S_{\tau}$ induces a family of Hilbert 
schemes 
$$\mathcal{H}_n \to \sym^{2m}_0(\cc),$$
whose fiber over $\tau$ is $\hi^n(S_{\tau}).$ The same arguments used in  
section~\ref{sec:open} carry over to show the following statement for 
families:

\begin {proposition}
\label {fam}
There is an open algebraic embedding $\hat{j}: \ss_n \to 
\mathcal{H}_n$ whose restriction to each $Y_{n, \tau}$ is the embedding 
$j: \yy_{n, \tau} \to \hi^n(S_{\tau})$ considered above. 
\end {proposition}

\subsection {Quiver varieties.}

This section is not necessary for understanding the rest of the paper, 
and the reader may skip it if she wishes. Our only goal here is to put our 
discussion in a larger context, by noting that both $Y_{n, \tau}$ 
and $\hi^n(S_{\tau})$ are particular examples of a more general concept. 
This concept is that of a quiver variety, introduced by Nakajima in \cite{N1}.
Our discussion follows \cite{N1} closely, but is not intended as a 
substitute for that paper. We simply want to provide a rough idea of what 
is involved for the reader who is not familiar with this rich subject.

A {\it quiver} $A$ is a finite oriented graph with no oriented cycles. Let 
$K$ be the set of vertices and $H$ the set of pairs 
consisting of an edge 
and an orientation on it, which may or may not be the one we chose 
initially. For $h\in H,$ we set $\epsilon(h)=1$ if the orientation agrees 
with the initial one, and $\epsilon(h) = -1$ otherwise. We also denote 
by $\text{out}(h), \text{in}(h) \in K$ the outgoing and incoming vertices 
of $h,$ respectively, and by $\bar h$ the same edge with reversed 
orientation.

For each vertex $k \in K,$ we choose a pair of Hermitian vector spaces 
$V_k, W_k$ of dimensions $v_k, w_k,$ respectively. We write the set of 
dimensions as vectors $\vv, \ww$ of length the cardinality of $K.$ We 
form the complex vector space
$$ \mathbb{M}(\vv, \ww) = \Bigl( \bigoplus_{h \in H} 
\hom(V_{\text{out}(h)},
V_{\text{in}(h)}) \Bigr) \oplus \Bigl( 
\bigoplus_{k \in K} \hom(W_k, V_k) \oplus \hom(V_k, W_k )\Bigr).$$

An element of $M$ consists of components $B_h, i_k, j_k.$ 
The group $G_{\vv} = \prod_{k \in K} U(V_k) $ with Lie algebra 
$\mathfrak{g}_{\vv}$ acts on $\mathbb{M}(\vv, \ww)$ by:
$$ (B_h, i_k, j_k) \to (g_{\text{in}(h)} B_h g_{\text{out}(h)}^{-1}, 
g_ki_k, j_k g_k^{-1} ). $$

Consider the map
$$ \mu = \mu_{\rr} \oplus \mu_{\cc} :  \mathbb{M}(\vv, \ww) \to 
\mathfrak{g}_{\vv} \oplus (\mathfrak{g}_{\vv} \otimes \cc),$$
$$ \mu_{\rr} (B, i, j) = \frac{i}{2} \Bigl( \sum_{h\in H: k = 
\text{in}(h)} 
B_h B_h^* - B_{\bar h}^* B_{\bar h} - i_k i_k^* - j_k^* j_k \Bigr)_k \in   
\bigoplus_k \mathfrak{u}(V_k) = \mathfrak{g}_{\vv}, $$
$$\mu_{\cc} (B, i, j) = \frac{i}{2} \Bigl( \sum_{h\in H: k = \text{in}(h)}
\epsilon(h) B_h B_{\bar h}  + i_k j_k \Bigr)_k \in
\bigoplus_k \mathfrak{gl}(V_k) = \mathfrak{g}_{\vv} \otimes \cc. $$

Let $Z_{\vv} \subset \mathfrak{g}_{\vv}$ be the center, and choose an 
element $\zeta = (\zeta_{\rr}, \zeta_{\cc}) \in Z_{\vv} \oplus (Z_{\vv} 
\otimes \cc).$

\begin {definition}
The quiver variety associated to the quiver $A,$ the vectors $\vv, \ww,$ 
and the parameter $\zeta$ is
$$ \mathfrak{M}_{\zeta}(\vv, \ww) = \{(B, i, j) \in \mathbb{M}(\vv, \ww) : 
\mu(B, i, j) = - \zeta \} / G_{\vv}.$$
\end {definition}

In fact, one can show that $ \mathbb{M}(\vv, \ww)$ is naturally a 
quaternionic vector space, and that the action of $G_{\vv}$ respects its 
hyper-K\"ahler structure. The map $\mu$ is the hyper-K\"ahler moment map 
associated to that action, and the quiver variety is the so-called 
hyper-K\"ahler quotient. As such, $\mathfrak{M}_{\zeta}(\vv, \ww)$ 
inherits a 
hyper-K\"ahler structure. It is a noncompact space and, in general, it can 
have singularities. It turns out that different choices of orientations on 
the same graph give rise to isomorphic quiver varieties.

When $\zeta_{\rr} = 0,$  $\mathfrak{M}_{\zeta}(\vv, \ww)$ is homeomorphic 
to the affine algebraic 
quotient of $\mu_{\cc}^{-1}(-\zeta_{\cc})$ by the action of the 
complexified group $G_{\vv}^{\cc} = \coprod_{k \in K} GL(V_k).$ 

There are numerous interesting examples of quiver varieties: cotangent 
bundles of generalized flag manifolds, nilpotent adjoint orbits in the Lie 
algebra $\gl_m,$ and their intersections with various transverse slices. 
For the last two their quiver variety description follows from a result of 
Kronheimer \cite{Kr2}, which should generalize for other orbits as 
well. In particular, we expect the variety $\yy_{n, \tau}$ from 
section~\ref{sec:yy} to be isomorphic to 
$\mathfrak{M}_{\zeta}(\vv, \ww)$ for the extended Dynkin graph of type 
$\tilde A_{2m-1},$ i.e. a polygon with $2m$ vertices numbered $0, 1, 
\dots, 
2m-1,$ with the following data:

$ \bullet \ v_k = k \text{ for }k \leq n; \ \ v_k = n \text{ for }n 
\leq k \leq 2m-n;\ \ v_k = 2m-k \text{ for } k \geq 2m-n;$

$ \bullet \ w_k = 1\text{ for }k = n, 2m-n\text{ and }n < m; w_n = 2 
\text{ for }n=m; w_k=0 \text{ otherwise}; $
 
$ \bullet \ \zeta_{\rr}=0; \text{ some }\zeta_{\cc} \text{ depending on 
}\tau.$

(As noted by Seidel and Smith \cite{SS}, this result is only conjectural 
at this point.)

\medskip
Another important class of examples comes from instanton theory and, in 
fact, this was the starting point for the whole subject. The moment map 
equations $\mu(B, i, j) = -\zeta$ above are modelled from the ADHM 
equations 
which describe the moduli space of instantons on $\rr^4.$ More 
generally, one can consider moduli spaces of instantons on the so-called 
ALE (asymptotically locally Euclidean) spaces \cite{Kr}. An ALE space 
is a minimal resolution of a singularity of the type $\cc^2/\Gamma, \Gamma 
\subset \text{SU}(2)$ finite, endowed with a certain hyper-K\"ahler 
metric. An example is the Milnor fiber $S_{\tau}$ from 
section~\ref{sec:ex}, which corresponds to the cyclic group $\Gamma= 
\zz/2m\zz.$ The ALE 
spaces can be described as quiver varieties \cite{Kr}, and the same is 
true for various moduli spaces of instantons on them \cite{KN}. 

Furthermore, by a similar argument to that given by Kronheimer and
Nakajima in \cite{KN}, the Hilbert scheme of an ALE space is also a quiver
variety. These Hilbert schemes have been studied by Qin and Wang in
\cite{QW}, \cite{W}. In particular, $\hi^n(S_{\tau})$ is isomorphic to
$\mathfrak{M}_{\zeta}(\vv, \ww)$ for the extended Dynkin graph of type
$\tilde A_{2m-1}$ with the following data:
                                                                                
$ \bullet \ v_k = n \text{ for all }k; $
                                                                                
$ \bullet \ w_k = 1\text{ for }k = 0, \  w_k=0 \text{ otherwise}; $
                                                                                
$ \bullet  \text{ some }\zeta_{\rr}, \zeta_{\cc} \text{ depending on
}\tau.$
 
\medskip 

Thus, for every $n \leq m,$ Theorem~\ref{hilb} describes an open holomorphic embedding of a quiver
variety into another one, where both varieties are constructed from the same quiver but with
different data $(\vv, \ww).$ These are the only nontrivial examples of this kind known to the
author. We leave the following as an open problem:

\begin {question}
What are the pairs of quiver varieties $(M, M')$ such that $M$ admits an 
open holomorphic embedding into $M'$ ? 
\end {question}

\section {The Seidel-Smith construction}

Seidel and Smith \cite{SS} 
defined a link invariant as the Floer cohomology of two specific 
Lagrangian submanifolds in $\yy_{n, \tau},$ for $n=m.$ In this section we 
summarize their construction. We skip most of the proofs, and we refer 
the reader to \cite{SS} for more details.

\subsection {K\"ahler metrics and parallel transport.} 
\label{sec:kal} Recall that the 
varieties $Y_{m, \tau}$ are fibers of the map:
$$ \chi|_{\ss_m}: \ss_m \to \sym^{2m}_0 (\cc).$$

Seidel and Smith endow the affine slice $\ss_m$ with a K\"ahler form
$\Omega$, and give the fibers the K\"ahler form induced by restriction.  
The form $\Omega$ is constructed as $\Omega = -dd^c\psi,$ for some smooth 
function $\psi: \ss_m \to \rr.$ The function $\psi$ is not
unique, but there is a weakly contractible parameter space of
possible choices, and in the end different choices give rise to the same
link invariant.

The main requirements for $\psi$ are related to its behaviour at infinity 
in $\ss_m.$ Basically, $\psi$ is asked to satisfy the following four 
conditions:
\begin {eqnarray}
\label {cond1}
&\bullet & -dd^c \psi > 0, \text{ so that } \Omega \text{ is K\"ahler.} \\
\label {cond2}
&\bullet & \psi  \text{ is proper and bounded below. } \\
\label{cond3}
& \bullet & \text{ Outside a compact set of }\ss_m, \|\nabla \psi\| \leq 
C \psi \text{ for some }C > 0; \\
\label {cond4}
& \bullet & \text{ The fiberwise critical set of } \psi \text{ maps 
properly to 
}\sym^{2m}_0(\cc) \text{ under } \chi.
\end {eqnarray}

The construction of $\psi$ with these properties starts with the 
observation that there is a natural action $\lambda$ of the multiplicative 
group $\rr_+$ on $\ss_m.$ In terms of the coordinates $a_k, b_k, c_k, 
d_k$ from section~\ref{sec:yy}, $k=1, \dots, m,$ the element $r \in \rr_+$ 
acts by: 
\begin {equation}
\label {action}
\lambda_r : (a_k, b_k, c_k, d_k) \to (r^k a_k, r^k b_k, r^k c_k, r^k d_k).
\end {equation}

This action takes the fiber $Y_{m, \tau}$ to $Y_{m, r \tau},$ where $\tau= 
(\tau_1, \dots, \tau_{2m}) \in \sym^{2m}_0(\cc)$ and $r$ acts on $\tau$ by 
multiplication on each component.

Fix some real number $\alpha > m.$ For each $k=1, \dots, m$ consider the 
function $\xi_k(z) = |z|^{\alpha/k}$ on $\cc,$ and add to it a compactly 
supported function $\eta_k$ on $\cc$ such that $\psi_k = \eta_k + \xi_k$ 
is $C^{\infty}$ and satisfies $-dd^c \psi_k > 0$ on $\cc.$ Apply $\psi_k$ 
to the coordinates 
$a_k, b_k, c_k, d_k$ for each $k$ and sum up these functions to obtain the
function $\psi$ on $\ss_m.$ Conditions $(\ref{cond1}), (\ref{cond2})$ and 
$(\ref{cond3})$ can be checked immediately, while for condition 
$(\ref{cond4})$ one needs to use the asymptotical homogeneity of $\psi$ 
with respect to the $\rr_+$ action (see Lemma~\ref{hom} below).

The reasons for the choice of $\psi$ are that conditions 
$(\ref{cond1})-(\ref{cond4})$ allow one to define {\it rescaled parallel 
transport maps,} in the following sense.

The map $\chi$ is a fibration over $\conf^{2m}_0(\cc)$ when restricted to 
$$W = \chi^{-1}(\conf^{2m}_0(\cc) \cap \ss_m).$$
Take a path $\gamma: [0,1] \to \conf^{2m}_0(\cc)$ on the base. The 
parallel transport vector field $H_{\gamma}$ on the pullback $\gamma^*W 
\to [0,1]$ consists of the sections of $TW|_{W_{\gamma(s)}}$ which project 
to $\gamma'(s),$ and are orthogonal to the tangent space along the fibers in 
the given K\"ahler metric. We would like to define a symplectic 
isomorphisms between the fibers by integrating $H_{\gamma}.$ However, 
this is not possible, because the fibers are not compact, and integral 
lines of $H_{\gamma}$ may go to infinity in finite time. The only thing we 
can say is that for every compact $P \subset W_{\gamma(s)}$ and $\epsilon 
> 0$ small 
(depending on $P$), there is a symplectic embedding 
$$ h_{\gamma}: P \to W_{\gamma(s+\epsilon)},$$
called {\it naive parallel transport}. 

However, Seidel and Smith show that one can say more when the 
K\"ahler metric on $\ss_m$ has well-chosen behaviour at infinity, such as 
our $-dd^c \psi,$ with $\psi$ satisfying (\ref{cond1})-(\ref{cond4}). 
In this case for every compact $P \subset W_{\gamma(0)},$ we can define a 
symplectic embedding:
$$ h_{\gamma}^{\resc}: P \to W_{\gamma(1)}.$$
This is the rescaled parallel transport, and is defined (roughly) by
subtracting from the usual parallel transport vector field a multiple of a
certain Liouville vector field on the fibers, integrating the resulting
vector field for all times, and then rescaling back on the image.
For a bigger compact set $P'$ containing $P,$ there is another rescaled 
parallel transport map to $W_{\gamma(1)},$ but its restriction to $P$ is 
isotopic to the previous one, in the class of symplectic embeddings. 
Therefore, given a closed Lagrangian 
submanifold $\ll \subset W_{\gamma(0)},$ there is a closed Lagrangian 
$h_{\gamma}^{\resc}(\ll) \subset W_{\gamma(1)},$ well-defined up to 
Lagrangian isotopy.

It is not hard to see from the definition in \cite{SS} that rescaled 
parallel transport behave nicely under composition of paths: 
if $\gamma: [0,2] \to \conf^{2m}_0(\cc)$ is smooth, then 
$h^{\resc}_{\gamma|_{[1,2]}} \circ h^{\resc}_{\gamma|_{[0,1]}}$ is 
isotopic to the $h^{\resc}_{\gamma}$ for the full path. 

Also, rescaled parallel transport has the same result as the naive
parallel transport for small paths, where the latter is defined. It
follows that in some cases, we can describe the image
$h_{\gamma}^{\resc}(\ll)$ of a closed Lagrangian $\ll \subset
W_{\gamma(0)}$ in terms of naive parallel transport as follows: we break
the path $\gamma$ into small pieces corresponding to a partition $0=t_0 <
t_1 < \dots < t_n = 1.$ We assume that $h_{\gamma|_{[t_0, t_1]}}(\ll)$ is
well-defined, and we use a Lagrangian isotopy to deform it into another
Lagrangian $\ll_1 \subset W_{\gamma(t_1)},$ such that $h_{\gamma|_{[t_1,
t_2]}}(\ll_1)$is well-defined. We iterate the process and end up with a
Lagrangian $\ll_n \subset W_{\gamma(1)},$ the same as
$h_{\gamma}^{\resc}(\ll)$ up to isotopy. Of course, there is no guarantee
that the partition and the Lagrangians $\ll_k$ considered above exist, but
we will see that in some cases they do. This is the strategy that we use
in section~\ref{sec:thm2} below.

\subsection {$A_1$ fibered singularities.} The following result describes the 
local structure of the singular part of $\yy_{m, \tau}$ when two of the 
points in $\tau$ become zero. It is seen to have a fibered singularity 
of $(A_1)$ type. 

Let $\bar \tau = (\mu_3, \dots, \mu_{2m}) \in \conf^{2m-2}_0(\cc)$ and 
consider the disk $D \subset \sym^{2m}_0(\cc)$ corresponding to 
eigenvalues 
$\tau(\zeta)= (-\sqrt{\zeta}, \sqrt{\zeta}, \mu_3, \dots, \mu_{2m})$ with 
$\zeta$ small. We have:

\begin {lemma}
\label {sing}
(i) The singular set $\sing(\yy_{m, \tau(0)})$ of $\yy_{m, \tau(0)}$ is 
canonically isomorphic to $\yy_{m-1, \bar \tau}.$ 

(ii) There is a neighborhood of $\sing(\yy_{m, \tau(0)})$ inside 
$\chi^{-1}(D) \cap \ss_m,$ and an isomorphism of that with a neighborhood 
of $(\yy_{m-1, \bar \tau}) \times \{0\}^3$ inside $(\yy_{m-1, \bar \tau}) 
\times \cc^3.$ This isomorphism is compatible with the one in (i) and 
fits into a commutative diagram:
$$\begin {CD}
\chi^{-1}(D) \cap \ss_m @>{\operatorname{local }\cong}>> (\yy_{m-1, \bar 
\tau}) \times \cc^3 \\
@V{\chi}VV @V{u^2 + v^2 +z^2}VV \\
D @>\zeta>> \cc
\end {CD}$$
where $u, v, z$ are the coordinates on $\cc^3.$
\end {lemma}

In terms of the coordinates $a_k, b_k, c_k, d_k$ on $\yy_{m, \tau(0)} 
\subset \ss_m,$ the isomorphism in (i) corresponds to simply setting $a_m 
= b_m = c_m = d_m =0,$ and identifying the other coordinates with the 
corresponding ones on $\yy_{m-1, \bar \tau} \subset \ss_m.$ Note that if 
we choose the same $\alpha$ and the same functions $\psi_k$ in the 
construction of the K\"ahler form for both $m$ and $m-1,$ and give the 
singular set the K\"ahler form induced by restriction, then the 
isomorphism in (i) is in fact a symplectomorphism.

\subsection {Relative vanishing cycles.}
\label{sec:vanc}

The next lemma deals with the construction
of relative vanishing cycles near the singularity considered above. Let
$X$ be a complex manifold. (In our applications $X$ will be taken to be
$\yy_{m-1, \bar \tau}.$) Give $Y = X \times \cc^3$ any K\"ahler metric,
and consider the map 
\begin {equation}
\label{cu}
 q: Y=X \times \cc^3 \to \cc, \ \ q(x,u,v,z) = u^2+v^2+z^2.
\end {equation}

We equip the fibers $Y_w=q^{-1}(w)$ with the induced metrics. The critical
point set of $q$ is $\{u=v=z=0\} \cong X,$ and the real part $Re(q)$ is a
Morse-Bott function. Let $Q \subset q^{-1}(\rr^{\geq 0})$ be the stable
manifold of $Re(q),$ and $l:Q \to X$ the map which assigns to a point in
$Q$ its limit under the negative gradient flow of $Re(q).$

\begin {lemma}
\label {rvc}
Let $K \subset X \cong X \times \{0\}$ be a compact Lagrangian 
submanifold. Then for sufficiently small $w > 0, L_w = l^{-1}(K) \cap Y_w$ 
is a Lagrangian submanifold of $Y_w$ diffeomorphic to $K \times S^2,$ 
and called the relative vanishing cycle associated to $K.$ $L_w$ can also 
be described in terms of parallel transport along the linear path $\gamma: 
[0,1] \to \cc, \ \gamma(s) = (1-s)t,$ as the set of points in $Y_w$ which 
are taken to a point in $K$ by the naive parallel transport along 
$\gamma,$ in the limit $s \to 1.$
\end {lemma} 

Note that by multiplying $q$ with some constant in $S^1,$ we can define 
stable manifolds which lie over other half-lines in $\cc,$ and 
corresponding relative vanishing cycles in $Y_w$ for all sufficiently 
small $w \in \cc^*.$

\subsection {The Lagrangians in terms of parallel 
transport.}\label{sec:lag}  As
mentioned in the introduction, Seidel and Smith defined their link
invariant with the help of two Lagrangians $\ll, \ll'\subset \yy_{m,
\tau}.$ Both $\ll$ and $\ll'$ are described in \cite{SS} in terms of the 
following construction, which makes use of the rescaled parallel transport 
maps and of vanishing cycles.
                                                                                
Let $\tau = (\mu_1, \dots, \mu_{2m}) \in \conf^{2m}(\cc).$ In 
section~\ref{sec:yy} we made the remark that all our constructions can be 
done for $\conf^{2m}(\cc)$ and $\sym^{2m}(\cc)$ via their projections $p$ 
to $\conf^{2m}(\cc)$ and $\sym^{2m}(\cc),$ respectively. In particular, 
$\yy_{m, \tau}$ is well-defined. We will associate 
a Lagrangian $\ll(\delta) \subset \yy_{m, \tau}$ to a
a collection of $m$ disjoint oriented arcs $\delta = (\delta_1, \dots, 
\delta_m)$
in $\cc,$ joining together the points of $\tau$ in pairs. Without loss of 
generality, we assume that $\delta_k$ runs from $\mu_{2k-1}$ to 
$\mu_{2k}.$ $\ll(\delta)$ will be well-defined up to Lagrangian isotopy.

Consider the path $\tilde \gamma: [0,1] \to \sym^{2m}(\cc)$ starting at
$\tau$ which keeps $\mu_3, \dots, \mu_{2m}$ fixed, and moves $\mu_1$ and
$\mu_2$ towards each other following $\delta_1.$ They collide at the
midpoint $\delta_1(1/2).$ We assume that the arc $\delta_1$ is a straight
horizontal line near its midpoint, and the two points move towards each 
other with
equal speed for $s$ close to $1.$ We compose that with the
natural projection to get a path $\gamma : [0,1] \to \sym^{2m}_0(\cc),
\gamma(s) = p \circ \tilde \gamma(s).$

The construction of $\ll(\delta)$ is done inductively on $m.$  We 
could start with the trivial case of a point as subset of a point for 
$m=0,$ but let us do the case $m=1$ for completeness. Then the adjoint 
quotient map is 
$$ \chi: \ss_1 = \sl_2(\cc) = \cc^3 \to \sym^2_0(\cc) \cong \cc$$ given 
by $\chi(u,v,z) = u^2+ v^2 +z^2.$ Here we identify the base with $\cc$ by 
$\tau = (\mu, -\mu) \to \mu^2,$ and let $P_{\tau}(z) = z^2 - \mu^2.$ 
The map $\chi$ has a unique critical point in the 
fiber over $\gamma(0)=0.$ 
In the nearby fibers $\gamma(s)$ for $s$ close to $1$ we have an 
associated 
vanishing cycle, a Lagrangian two-sphere. We then use the reverse rescaled 
parallel transport maps along $\gamma$ to move this to a Lagrangian 
$\ll(\delta) \subset \yy_{1, \tau}.$ 

The inductive step is similar. Let $\tau' = (\mu_1', \mu_2', \mu_3',
\dots, \mu_{2m}') = \gamma(1) \in \sym^{2m}_0(\cc),$ with $\mu_1' =
\mu_2'.$ For small $s$ we have $\gamma(1-s) = \tau'_s =(\mu_1' - s,
\mu_1' +s, \mu_3', \dots, \mu_{2m}').$ Let also $\bar \tau = (\mu_3'',
\dots, \mu_{2m}'') \in \sym^{2m-2}_0(\cc)$ be the image of $(\mu_3, \dots,
\mu_{2m})$ under the projection $p.$ It follows that $\tau''_s =(-s,s,
\mu_3'', \dots, \mu_{2m}'')$ is in $\sym^{2m}_0(\cc)$ as well. Assume we
already have a compact Lagrangian in $\yy_{m-1, \bar \tau}$ as in section
\ref{sec:vanc}. Using Lemmas~\ref{sing} and \ref{rvc} we get a relative
vanishing cycle in $\yy_{m, \tau''_s}$ for small $s.$ We move it to 
$\yy_{m, \gamma(1-s)}$ using rescaled parallel transport along the linear 
path going from $\tau''_s$ to $\tau'_s.$ Rescaled parallel
transport is then used again in reverse along $\gamma$ to give a 
Lagrangian submanifold in $\yy_{m,\tau},$ which is the desired $\ll(\delta).$

Note: Seidel and Smith interpolate between $\tau''_s$ and $\tau'_s$ for 
$s=0$ in the space composed of the singular sets of $\yy_m$'s, then use 
the vanishing cycle construction at $\tau'_s$ rather than $\tau''_s.$ 
However, the two constructions are easily seen to be equivalent up to 
isotopy.

\subsection {Links as braid closures.}
\label {sec:braid}
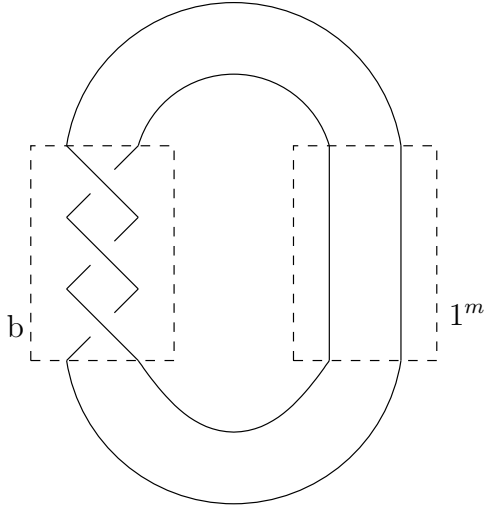
\begin {figure}
\begin {center}
\input {circ.pstex_t}
\end {center}
\caption {The left-handed trefoil.}
\label {fig:circ}
\end {figure}

Given a link $L,$ we can present it as the closure of an $m$-stranded
braid $b \in Br_m.$ (See Figure~\ref{fig:circ} for a presentation of the
left-handed trefoil.)
                                                                                
Equivalently, $L$ is the plat closure of the braid $b \times 1^m \in Br_m
\times Br_m \hookrightarrow Br_{2m}.$ We represent $b \times 1^m$ by a
loop $l: [0,1] \to \conf^{2m}(\cc).$ Our convention is that braids act on 
the $2m$ punctured plane on the right, with geometric braids reading from top 
to bottom, in the sense that the first generator of the braid from the top 
acts first, and so forth. 
                   
\begin {figure}
\begin {center}
\input {match.pstex_t}
\end {center}
\caption {Standard crossingless matching.}
\label {fig:match}
\end {figure}
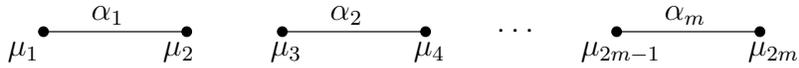
                                                             
Consider the standard crossingless matching of $2m$ points in the plane in
Figure~\ref{fig:match}. The endpoints of the $m$ segments $\alpha_1,
\dots, \alpha_m$ are $\mu_1, \dots, \mu_{2m} \in \rr \subset \cc,$ with
$\tau = (\mu_1, \dots, \mu_{2m}) \in \conf^{2m}_0(\cc).$ The odd points
$\mu_1, \mu_3, \dots, \mu_{2m-1}$ correspond to the strands on the left
side of Figure~\ref{fig:circ}, and the even points $\mu_2, \dots,
\mu_{2m}$ corresponds to the vertical strands on the right.

Note that Seidel and Smith consider a different crossingless matching to
be standard in \cite{SS}, but their picture is equivalent to ours after
isotopy and conjugation with a fixed braid, and conjugation does
not change the link type.
                                                                                
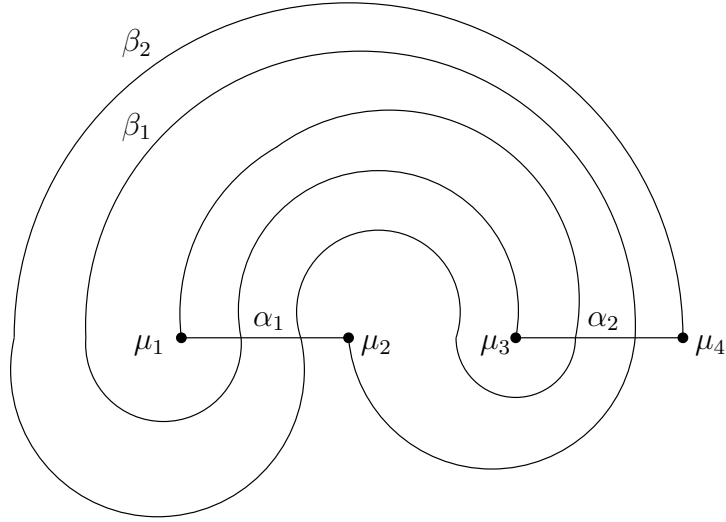
\begin {figure}
\begin {center}
\input {tref.pstex_t}
\caption {Flattened braid diagram of the trefoil.}
\label {fig:flat}
\end {center}
\end {figure}

Given our loop $l: [0,1] \to \conf^{2m}(\cc),$ we can find a smooth family
of crossingless matchings in the plane with endpoints $l(s), s\in [0,1].$
Note that $l(s)$ are $2m$-tuples of points in $\cc,$ and $\mu_2, \dots,
\mu_{2m}$ always appear among these $2m$ points.  At time $s=1$ we get a
crossingless matching composed of $m$ segments $\beta_1, \dots, \beta_m$
joining the points of $\tau$ in pairs. We assume that the $\alpha$ and
$\beta$ are simple curves that intersect transversely in their interior.
Then the original link $L$ can be recovered from the diagram of the alpha
and beta curves intersecting in the plane, which we call a {\it flattened
braid diagram} for $L.$ Indeed, if we choose the alpha curves to serve as
underpasses at each crossing in a flattened braid diagram, we obtain a
usual plane diagram for the link. This is shown in Figure~\ref{fig:flat}
for the trefoil.

The Lagrangians $\ll$ and $\ll'$ considered by Seidel and Smith are then
$\ll(\delta)$ when the crossingless matching $\delta$ is given by the
alpha and beta curves, respectively. 

\begin {figure}
\begin {center}
\input {unknot.pstex_t}
\caption {Two flattened diagrams of the $1$-braided unknot.}
\label {fig:unknot}
\end {center}
\end {figure}
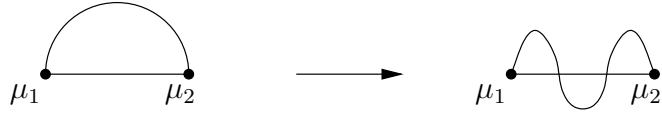

\begin {remark}
The braid $b$ determines the flattened diagram up to isotopies that keep
the $\mu_k$'s fixed. Note that such isotopies can introduce new
intersections between the alpha and beta curves, as shown in
Figure~\ref{fig:unknot} for the $1$-braided unknot.
\end {remark}                                                                                

\subsection {Floer cohomology.} \label{sec:floer} Lagrangian Floer 
cohomology \cite{F} is a
very rich subject. The version used in \cite{SS} takes place in a K\"ahler
manifold $(Y, \Omega)$ such that 
\begin {equation}
\label {manifold}
\Omega \text{ is exact } \ ; \  Y \text{ is Stein } ; \ c_1(Y) = 0 \ ; \ 
H^1(Y) = 0.
\end {equation}

We also need two closed connected Lagrangian submanifolds $\ll,
\ll' \subset Y$ such that
\begin {equation}
\label {lagrangians}
H_1(\ll) = H_1(\ll')=0 \ ; \ w_2 (\ll) = w_2(\ll') = 0
\end {equation}

One can check that $Y = \yy_{m, \tau}$ and $\ll = \ll(\alpha), 
\ll'=\ll(\beta)$ satisfy these conditions. In general, for any $Y$ and 
$\ll, \ll'$ satisfying (\ref{manifold}) and (\ref{lagrangians}) there is a 
well-defined abelian group with a relative $\zz$ grading: $$
HF^*(\ll, \ll') = H(CF^*(\ll, \ll'), d).$$
                                                                                
This is called Floer cohomology and is obtained from a cochain group
$CF^*(\ll, \ll')$ using a differential $d.$ The $\zz$ grading is relative
in the sense that it is well-defined only up to an overall constant shift.
                                                                                
A short overview of the main properties of $HF^*$ relevant to their
construction is given by Seidel and Smith in \cite[section 4(D)]{SS}.  
The most important property is that $HF^*$ is invariant under smooth
deformations of the objects involved, provided that one remains within a
class where $HF^*$ is well-defined. For example, $HF^*$ is invariant under
Lagrangian isotopies of $\ll$ and $\ll'.$ In this paper we are mostly
interested in understanding the cochain complex $CF^*.$ Let us just
quickly note that the differential $d$ is defined by counting
``pseudo-holomorphic disks'' for a family of almost complex structures
taming $\Omega.$ These are solutions to some PDE's similar to the
Cauchy-Riemann equations, with certain boundary conditions.

A set of generators over $\zz$ for the cochain complex $CF^*$ is 
obtained by isotoping one Lagrangian so that the intersection $\ll \cap 
\ll'$ is transverse, and then taking a generator for each point in $\ll 
\cap \ll'.$ The relative grading is obtained from a Maslov index 
calculation. In some cases, including the one of interest to us, it can be 
improved to an absolute $\zz$ grading. This was done by Seidel in  
\cite{Se}, following the ideas of Kontsevich \cite{Ko}.

Let $\lfrak \to Y$ be the natural fiber bundle over $Y$ whose fibers $\lfrak_x$ are the manifolds 
of
Lagrangian subspaces of $T_xY.$ Then $\pi_1(\ll_x) \cong \zz$ has a canonical generator called the
Maslov class. 

Since $c_1(Y) = 0,$ we can pick a complex volume form
$\Theta$ on $Y,$ i.e. a nowhere vanishing section of the canonical 
bundle. This determines a square phase map 
\begin {equation}
\label {sqp1}
 \theta : \lfrak \to
\cc^*/\rr_+, \ \theta(V) = \Theta(e_1 \wedge \dots \wedge e_N)^2 
\end {equation}
for any orthonormal basis $e_1, \dots, e_N$ of $V \subset T_xY.$ Here $2N= \dim_{\rr} Y.$  We can
identify $\cc^*/\rr_+$ with $S^1$ by the contraction $z \to z/|z|.$
Let $\tilde \lfrak \to \lfrak$ be the infinite cyclic covering obtained from the 
universal covering $\rr \to S^1$ by pulling it back under the map $\theta.$ 

Note that the Lagrangian submanifold $\ll$ gives a canonical section $s_{\ll} : \ll
\to \lfrak, s_{\ll}(x) = T_x \ll.$ This produces a map
\begin {equation}
\label {sqp}
\theta_{\ll} = \theta \circ s_{\ll} : \ll \to \cc^*/\rr_+ \cong S^1.
\end {equation}

The condition $H^1(\ll) = 0$ allows us to lift $s_{\ll}$ to a section $\tilde s_{\ll}: \ll 
\to \tilde \lfrak.$ This is equivalent to lifting $\theta_{\ll}$ to a map $\tilde \theta_{\ll}: \ll 
\to \rr.$ 

\begin {definition}
\label {deaf}
A {\bf grading} on $\ll$ is a choice of a lift $\tilde \theta_{\ll} :\ll \to
\rr.$
\end {definition}

If we choose a grading for $\ll'$ as well, then every
point $x$ in the intersection $\ll \cap \ll'$ (which was assumed to be
transverse) has a well-defined absolute Maslov index $\mu(x) \in \zz$ \cite{Se}.
To define it, take a path $\tilde \lambda : [0,1] \to \tilde \lfrak_x$ with endpoints $\tilde 
\lambda (0) = \tilde s_{\ll} (x)$ and $\tilde \lambda (1) = \tilde s_{\ll'}(x).$
The projection of $\tilde \lambda$ to $\lfrak_x$ is a path $\lambda_0: [0,1] \to \lfrak_x$ joining 
$T_x \ll$ and $T_x \ll'.$ Consider 
also the constant path $\lambda_1 : [0,1] \to \lfrak_x, \lambda_1(t) = T_x\ll'.$ Set
\begin {equation}
\label {formk}
 \mu(x) = \mu^{\text{paths}}(\lambda_0, \lambda_1) + \frac{N}{2},
\end {equation}
where $\mu^{\text{paths}}(\lambda_0, \lambda_1)$ is the Maslov index for paths in $\lfrak_x$ defined 
in \cite{RS}.
  
The Maslov grading $\mu$ induces an absolute $\zz$ grading on the cochain complex and
on cohomology. The result does not depend on the choice of $\theta,$
because the condition $H^1(Y) =0$ ensures that different choices are
homotopic in the class of smooth trivializations of the canonical bundle.

If $\ll \to \ll[1]$ denotes the process which subtracts the constant $1$
from the grading, we have 
\begin {equation}
\label {grad}
 CF^*(\ll, \ll'[1]) = CF^*(\ll[-1], \ll') =
CF^{*+1} (\ll, \ll')
\end {equation}
and the same holds for cohomology.

In our case $Y = \yy_{m, \tau}$ and $\ll=\ll(\alpha), \ll'=\ll(\beta)$
there is a way of eliminating even this last $\zz$ indeterminacy in $HF^*$
by choosing the gradings consistently on $\ll$ and $\ll'.$ Specifically,
we start with arbitary complex volume forms on the slice $\ss_m\cong
\cc^{4m-1}$ and on the base $\sym^{2m}_0(\cc) \cong \cc^{2m-1}.$ This
induces a smooth family of complex volume forms on the smooth fibers
$\yy_{m, \tau}, \tau \in \conf^{2m}_0(\cc).$ We know that $\ll(\beta)$ is
obtained from $\ll(\alpha)$ by following the family of crossingless
matchings in the plane determined by the braid $b \times 1^m.$ Given an
arbitrary grading of $\ll = \ll(\alpha),$ we can continue it uniquely to a
smooth family of gradings on the respective Lagrangians, and end with a
grading on $\ll' = \ll(\beta).$ Adding a constant to the grading on $\ll$ 
affects the one on $\ll'$ in the same way, so by (\ref{grad}) the grading 
on $CF^*$ remains the same. Thus it makes sense to write $HF^k(\ll, \ll')$ 
for any $k \in \zz.$

As noted in the introduction, the main result proved in \cite{SS} is: 

\begin {theorem}[Seidel-Smith] 
\label {theoss}
Let us denote by $w$ the writhe of the
braid $b \in Br_m$ i.e. the number of positive minus the number of
negative crossings (in Figure~\ref{fig:circ}, not in the flattened braid
diagram!). Then the Floer cohomology groups $$ Kh_{symp}^* 
= HF^{*+m+w}(\ll, \ll') $$ are link invariants. \end {theorem}

Conditions (\ref{manifold}) and (\ref{lagrangians}) can be weakened in various ways. In many cases 
we 
can still define a version of Lagrangian Floer cohomology at the expense of giving up some nice 
properties: for example, the Floer groups can be only $\zz/2$-graded, or only defined over a Novikov 
ring, etc. Also, typically Floer cohomology is only invariant under a restricted class of 
deformations, e.g. Hamiltonian isotopies of $\ll, \ll'$ instead of all Lagrangian isotopies.
For a discussion of Lagrangian Floer cohomology in a very general setting, we refer to \cite{FOOO}.

Furthermore, sometimes we can define Floer cohomology for half-dimensional submanifolds $\ttt, \ttt' 
\subset Y$ which are not Lagrangian, but only totally real. 

\begin {definition}
A real subspace $V \subset \cc^N$ is called {\bf totally real} (with respect to the standard 
complex structure) if $\dim_{\rr} V = N$ and $V \cap iV = 0.$ A half-dimensional submanifold $\ttt$ 
of an almost complex manifold $(Y, J)$ is called {\bf totally real} if $T_x\ttt \cap J(T_x \ttt) = 
0$ for all $x \in \ttt.$
\end {definition}

Pseudo-holomorphic disks with boundary on totally real submanifolds were studied in \cite{O2}.  
Unlike in the Lagrangian case where certain energy bounds come for free, here one needs to make sure
that these bounds still hold. In some cases they do and then a Floer cohomology can be defined.  
Most notably, this is the setting for the Heegaard Floer theory of Ozsv\'ath and Szab\'o \cite{OS1}.

On the other hand, the problem of defining the Maslov index is no more difficult in the totally real
case than in the Lagrangian case. The relative index is treated in \cite{O1} and \cite{KL}, and then 
this can easily be improved to an absolute index in the spirit of \cite{Se}. Let $\cc^N$ be endowed
with its standard symplectic form and complex structure. Denote by $\lfrak_n$ and $\tfrak_n$ the
spaces of Lagrangian and totally real subspaces of $\cc^N,$ respectively. Then the inclusion 
\begin {equation}
\label {heq}
\lfrak_n \cong U(N)/ O(N) \ \hookrightarrow \ GL(N, \cc)/GL(N, \rr) \cong \tfrak_n
\end {equation}
is a homotopy equivalence (see \cite{KL} or \cite{O1}).

We discuss here the case which will be relevant for us in section~\ref{sec:gradf}. $Y$ is a K\"ahler
manifold with $c_1(Y)=0,$ and $\Theta$ a complex volume form as before. Let $\ttt, \ttt' \subset Y$ 
be totally real and intersecting transversely. We do not want to assume neither that $H^1(\ttt) = 
H^1(\ttt') = 0$ nor that $H^1(Y) =0.$ Just like
in the Lagrangian case, there is a natural bundle $\tfrak \to Y$ whose fibers $\tfrak_x$ are the
manifolds of totally real subspaces of $T_xY.$ There is also a section $s_{\ttt}: \ttt \to \tfrak$
and the square phase maps $\theta: \tfrak \to \cc^* /\rr_+ \cong S^1$ and $\theta_{\ttt}:\ttt \to
S^1$ are still well-defined. 

We construct the infinite cyclic covering $\tilde \tfrak \to \tfrak$ as before. Since $H^1(\ttt)$
may be nonzero, it is possible that $s_{\ttt}$ does not lift to a section in $\tilde \tfrak.$
However, let us assume it does. A {\it grading} on $\ttt$ is a choice of such a lift. Assuming that
$\ttt'$ also has a grading, we can define an absolute Maslov grading $\mu(x) \in \zz$ for every $x
\in \ttt \cap \ttt'.$ Indeed, we can construct paths $\lambda_0$ and $\lambda_1$ in $\tfrak_x$ as
before, get corresponding paths in $\lfrak_x$ using the homotopy equivalence (\ref{heq}), and
then use formula (\ref{formk}) for those. Note that this time the result can depend on $\theta,$
because $H^1(Y)$ may be nonzero.

\section {A different K\"ahler metric}
                                                                                
In this section we go through the Seidel-Smith construction again, but
using a different choice of the K\"ahler metric on $\yy_{m, \tau}.$ This
will make it easier to write down the Lagrangians in terms of
the embedding of $\yy_{m, \tau}$ into the Hilbert scheme. The goal is to
prove Theorem~\ref{lagr}.

\subsection {A choice of K\"ahler form.}
\label {sec:newk}                                                                                
                                                                                
Corollary~\ref{cor} says that via the Hilbert-Chow
morphism $\pi$ from (\ref{hc}), we can identify an open subset of $\yy_{m, 
\tau}$ with
\begin {equation}
\label {vmt}
 \ve_{m, \tau} = \{(u_j, v_j, z_j) \in 
\sym^m(S_{\tau}): z_i \neq z_j
\text{ for all }i\neq j \}.
\end {equation}
                                                                                
The subsets $\ve_{m, \tau}$ form a family $\ve_m$ over
$\sym^{2m}_0(\cc).$ Using Proposition~\ref{fam}, we can identify $\ve_m$
with an open subset of $\ss_m.$ For a point $(\mu_1, \dots, \mu_{2m}) \in 
\sym^{2m}_0(\cc),$ we form the symmetric polynomials 
$$\nu_j = \sum_{i_1 < \dots < i_j} \mu_{i_1} \mu_{i_2} \dots \mu_{i_j}$$
for $2 \leq j \leq 2m.$ The $\nu_j$ are coordinates on $\sym^{2m}_0(\cc)$ 
viewed as an affine space $\cc^{2m-1}.$ 

We construct a new K\"ahler form $\tilde \Omega = -dd^c \tilde \psi$ on
$\ss_m,$ guided by the following requirement: We fix a large 
relatively compact, open subset $U$ of $\ve_m.$ We want the restriction of 
$\tilde \Omega$ to $U$ to have the form
\begin {equation}
\label {form}
\tilde \Omega|_U =  \frac{i}{2} \sum_{j=1}^m \bigl(du_j \wedge d\bar u_j + 
dv_j \wedge 
d\bar v_j + dz_j \wedge d\bar z_j \bigr) + \frac{i}{2}\sum_{j=2}^{2m} 
d\nu_j \wedge d\bar \nu_j. 
\end {equation}
This corresponds to $\tilde \psi$ of the form
$$ \sum_{j=1}^m (|u_j|^2 + |v_j|^2 + |z_j|^2) + \sum_{j=2}^{2m} 
|\nu_j|^2.$$

We also want $\tilde \psi$ to satisfy conditions 
(\ref{cond1})-(\ref{cond4}), so that $\tilde \Omega$ has well-defined 
rescaled parallel transport maps and the construction of the Lagrangians 
can proceed as before. 

Recall the construction of $\psi$ from section~\ref{sec:kal}. In it we use
smooth functions $\psi_k: \cc \to \rr$ $(k=1, \dots, m)$ with 
$-dd^c \psi_k > 0$ everywhere
and $\psi_k(z) = \xi_k(z)= |z|^{\alpha/k}$ for $z$ near infinity. 
We define similar functions $\xi_k, \psi_k$ for $k = m+1, \dots, 2m$ as 
well. If $\alpha$ is
chosen so that $\alpha > 4m,$ then given any relatively compact
open subset of $\cc,$ we can construct $\psi_k$ with these properties 
so that $\psi_k(z) = |z|^2$ on the specified open subset. We take this 
choice of $\psi_k$ and define a function on the Milnor fiber 
$S_{\tau} : (u^2 + v^2 + P_{\tau}(z) =0)$ by 
$$ (u, v, z) \to \psi_m(u) + \psi_m(v) + \psi_1(z).$$
Summing up these functions for all coordinates gives a function on the 
symmetric product $\sym^{m}(S_{\tau}).$ We pull it back to the Hilbert 
scheme $\hi^m(S_{\tau})$ via the Hilbert-Chow morphism $\pi$, then add the 
terms
$$ \psi_2 (\nu_2) + \psi_3(\nu_3) + \dots + \psi_{2m} (\nu_{2m})$$
to obtain a plurisubharmonic function on the family of Hilbert 
schemes $\hi^m(S_{\tau})$ over $\sym^0_{2m}(\cc).$ By 
Proposition~\ref{fam}, $\ss_m$ is an open subset of this family, hence by 
restriction we obtain a plurisubharmonic function $\rho: \ss_m \to \rr.$
   
We cannot use $\rho$ to construct a K\"ahler metric, because $-dd^c\rho$
is degenerate on the preimages $\pi^{-1}(\Delta)\cap \yy_{m, \tau} $ of
the diagonals $\Delta \subset \sym^m(S_{\tau}),$ i.e. outside $\ve_m =
\cup_{\tau} \ve_{m, \tau}.$ The function $\psi$ which was used to
construct $\Omega$ in section~\ref{sec:kal} comes to the rescue. Let
$\beta: \ss_m \to \rr$ be a bump function that is identically $1$ on a
neighborhood of $\ss_m - \ve_m $ and identically $0$ outside
a slightly larger neighborhood of $\ss_m - \ve_m.$ Then for
$\epsilon > 0$ sufficiently small, the function \begin {equation} \label
{tpsi} \tilde \psi = \rho + \epsilon \beta \cdot \psi \end {equation} is
strictly plurisubharmonic and defines a K\"ahler metric $\tilde \Omega =
-dd^c\tilde \psi$ on $\ss_m,$ and has the form (\ref{form}) on $U.$

Note that if we had used the functions $\xi_k$ instead of $\psi_k$ in the 
construction of $\tilde \psi,$ we would have obtained a 
corresponding $C^1$ function $ \tilde \xi: \ss_m \to \rr$ that is 
identical to $\tilde \psi$ near infinity.

\subsection {Rescaled parallel transport.} 
\label {sec:gg}

Let us check that $\tilde \psi$ satisfies conditions $(\ref{cond1}) - 
(\ref{cond4})$ and therefore gives rise to well-defined rescaled parallel 
transport maps. This runs parallel to the corresponding proof for $\psi$ 
given in section 5(A) of \cite{SS}.

Conditions $(\ref{cond1})$ and $(\ref{cond2})$ are immediate from the 
definition. To check $(\ref{cond3}),$ note that $\|\nabla \psi_k \| \leq 
c+ d\psi_k$ for some $c, d > 0.$ Adding up these inequalities, together 
with the fact that $\psi \gg 0$ outside a compact subset, shows that 
$\|\nabla \psi \| < C \psi$ for some $C > 0.$

For condition $(\ref{cond4}),$ we need to use the $\rr_+$ action 
(\ref{action}) on $\ss_m.$ On the open set $\ve_m \subset \ss_m$ the 
action of $r \in \rr_+$ can be written in terms of $u_j, v_j, z_j, \nu_j$ 
as
$$  \lambda_r : (u_j, v_j, z_j, \nu_j) \to (r^m u_j, r^m v_j, rz_j, r^j 
\nu_j). $$

Because the set $\ve_m$ is preserved by this action, we can choose the bump function
$\beta$ to satisfy $\beta \circ \lambda_r = \beta.$ Then $\tilde \xi$ is homogeneous of 
weight $\alpha,$ i.e. $\tilde \xi \circ \lambda_r = r^{\alpha}\tilde \xi.$

\begin {lemma}
\label {hom}
The function $\tilde \psi$ is asymptotically homogeneous with respect to 
the action $\lambda,$ in the sense that
$$ \frac{\tilde \psi \circ \lambda_r}{r^{\alpha}} \to  \tilde \xi \ 
\text{as } r \to \infty $$
where the convergence is uniform in $C^1$ norm.
\end {lemma}

\noindent\textit {Proof.} This is the analogue of Lemma 40 in \cite{SS}. 
The differences $\psi_k - \xi_k$ are compactly supported, and after 
rescaling the support of $\psi_k (r^k z) - \xi_k(r^k z)$ is 
getting smaller as $r \to \infty.$ This gives a uniform bound on their 
$C^0$ norms. For their derivatives we use the fact that $(\psi_k 
(r^k z) - \xi_k(r^k z))/r^{\alpha} \to 0$ uniformly for $\alpha > 4m > k.$ 
$\hfill \fin$

\medskip 

Now, to see that condition $(\ref{cond4})$ is satisfied, we can apply the
argument used by Seidel and Smith in Lemma 41 from \cite{SS}.  Here is a
sketch: there is a simultaneous resolution of the map $\chi:  \ss_m \to
\sym^{2m}_0(\cc)$ in the form of a differentiable fiber bundle $\hat
{\chi}: \hat{\ss}_m \to \cc^{2m-1}.$ The $\rr_+$ action lifts to one on
$\hat{\ss}_m.$ It preserves the fiber over $0,$ and the asymptotic
homogeneity of the lift of $\tilde \psi$ (Lemma~\ref{hom}) shows that the
critical set of the restriction of this lift to $\hat{\chi}^{-1}(0)$ is
compact. Then one shows that the whole fiberwise critical set of the lift
of $\tilde \psi$ to $\hat{\ss}_m$ maps properly to the base by using a
rescaling argument.

\subsection {Restriction and interpolation}
\label {sec:res}

The transverse slice $\ss_{m-1}$ sits embedded in $\ss_m$ by forgetting 
the coordinates $a_m, b_m, c_m, d_m.$ This is compatible with the 
isomorphism in Lemma~\ref{sing} (i). In terms of the coordinates $u_j, 
v_j, z_j, \nu_i$ on $\ve_m \subset \ss_m,$ the corresponding coordinates 
on $\ve_{m-1} \subset \ss_{m-1}$ are
$$ u'_j = \frac{u_j}{z_j} \ ; \   v'_j = \frac{v_j}{z_j} \ ; \ z'_j =  z_j 
\ ; \ \nu'_i =  \nu_i $$
for $j=1, \dots, m-1$ and $i=2, \dots, 2m-2.$

One disadvantage that our K\"ahler form $\tilde \Omega$ has over Seidel 
and Smith's $\Omega$ is that its restriction to $\ss_{m-1}$ is not of the 
same form. For example, on a big open subset $U' \subset U \cap \ss_{m-1}$ 
avoiding $z_j' = 0,$ instead of $(\ref{form})$ we 
have
\begin {equation}
\label {form2}
 \tilde \Omega|_{U'} = \frac{i}{2} \sum_{j=1}^{m-1} 
\bigl(d(u'_jz'_j) \wedge d(\bar u'_j\bar 
z'_j) + d(v'_j z'_j) \wedge
d(\bar v'_j \bar z'_j)+ dz'_j \wedge d\bar z'_j \bigr) + 
\frac{i}{2}\sum_{j=2}^{2m-2}
d\nu'_j \wedge d\bar \nu'_j. 
\end {equation}

Therefore, if in the recursive construction of the Lagrangians 
from section~\ref{sec:lag} we used the forms $\tilde\Omega_k$ on  
$\ss_k$ for all $k \leq m,$ where $\tilde \Omega_k$ are constructed 
just like $\tilde \Omega_m = \tilde \Omega$ in section~\ref{sec:newk},
we could run into the problem that 
the Lagrangians constructed in the singular set of $\yy_{k, \tau}$ may not 
be Lagrangians for the K\"ahler metric on $\ss_k.$ We deal with this 
problem by using the restrictions of $\tilde \Omega_m = \tilde \Omega$ to 
each $\ss_k$ 
instead. Their properties are not very different from those of $\tilde 
\Omega_k;$ in particular, they have a similar behaviour at infinity.

Let us interpolate linearly between the two functions $\psi, \tilde \psi: 
\ss_m \to \rr,$ by defining $\psi_t = (1-t)\psi + t \tilde \psi$ for $t 
\in [0,1].$ The proof of the following lemma is entirely similar to the 
argument in section~\ref{sec:gg} above, so we omit it.

\begin {lemma}
The functions $\psi_t$ and their restrictions to any $\ss_k \subset 
\ss_m, k < m$ satisfy the conditions $(\ref{cond1}) - (\ref{cond4}).$ Thus 
the corresponding K\"ahler forms have well-defined rescaled parallel 
transport maps over the respective configuration spaces. 
\end {lemma}

Since Floer cohomology is invariant under deformation, it follows that we 
can define Lagrangians $\tilde \ll, \tilde \ll'$ in $\ss_m$ with the form 
$\tilde \Omega$ and 
\begin {equation}
\label {hfl}
 HF^*(\tilde \ll, \tilde \ll') = HF^*(\ll, \ll').
\end {equation}

Therefore, the link invariant $Kh^*_{symp}$ can also be defined using 
$\tilde \Omega$ rather than $\Omega.$

\subsection {The new Lagrangians in an explicit form.}
\label {sec:thm2}

In section~\ref{sec:lag} we associated a Lagrangian $\ll(\delta)\subset
\yy_{m, \tau}$ to every crossingless matching $\delta$ consisting of arcs
$\delta_1, \dots, \delta_m$ joining the points of $\tau$ in pairs. Let
$\tilde \ll(\delta)$ be the corresponding Lagrangian constructed using the
K\"ahler form $\tilde \Omega$ in place of $\Omega.$

To each arc $\delta_k$ we associate the Lagrangian 2-sphere in the Milnor
fiber $S_{\tau}$: $$ \Sigma_{\delta_k} = \{ (u, v, z) \in S_{\tau}
\hskip3pt : z = \delta_k(t) \text{ for some }t \in [0,1]; \hskip5pt u, v
\in \sqrt{-P_{\tau}(z)} \rr \}.$$ These Lagrangians have appeared before
in the work of Khovanov and Seidel \cite{KS}.

Let us change coordinates on $S_{\tau}$ from $(u, v, z)$ to $(x, y, z)$ 
where
$$ x = \frac{u+iv}{\sqrt{2}} \ ; \ y  = \frac{u+iv}{\sqrt{2}}.$$ 
Then the equation for $S_{\tau}$ changes to $2xy + P_{\tau}(z) =0$ and
 $$ \Sigma_{\delta_k} = \{ (x, y, z) \in S_{\tau}
\hskip3pt : |x|=|y|, \ z = \delta_k(t) \text{ for some }t \in [0,1] 
\}.$$
We change coordinates on $\ve_{m, \tau}$ accordingly, from $(u_k, v_k, 
z_k)$ 
to $(x_k, y_k, z_k)$ for $k=1, \dots, m.$

\begin {proposition}
\label {super}
$\tilde \ll(\delta)$ is Lagrangian isotopic to $\kk(\delta)
= \Sigma_{\delta_1} \times 
\Sigma_{\delta_2} \times \dots \times \Sigma_{\delta_m} \subset \ve_{m, 
\tau} \subset \yy_{m, \tau},$ where $\ve_{m, \tau}$ is as in (\ref{vmt}).
\end {proposition}

\noindent \textit {Proof.} 
First of all let us note that we can only be sure that $\kk(\delta)$ is a 
Lagrangian for $\tilde \Omega$ if we know that it lives inside the chosen 
subset $U \subset \ve_m$ where $(\ref{form})$ holds. We will arrange so 
that all of our constructions take place inside $U.$

We do induction on $m$ starting with the trivial case $m=0.$ Let us 
explain the inductive step. We return to the notations from 
section~\ref{sec:lag}. Assume the recursive procedure had already 
given us the Lagrangian $\kk(\bar \delta) \subset \yy_{m-1, \bar \tau},$ 
where $\bar \delta$ is the crossingless matching obtained from $\delta$ by 
erasing $\delta_1$ and then translating everything from $(\mu_3, \dots, 
\mu_{2m})$ to $\bar \tau = (\mu''_3, \dots, \mu''_{2m})$ by (\ref{proj}). 

We apply the vanishing cycle procedure and then rescaled parallel 
transport to obtain a Lagrangian in $\yy_{m, \tau}.$ To control the form 
of the Lagrangians we use a moment map for a torus action. The argument is 
similar to the one used by Seidel and Smith in Lemma (32) of \cite{SS}.

The torus $T^m = (S^1)^m$ acts on $\ve_m$ preserving the fibers $\ve_{m, 
\tau}$ in the following way. The action of an element $(e^{i\theta_1}, 
\dots, e^{i \theta_m}) \in T^m$ is 
$$ (x_k, y_k, z_k) \to (e^{i\theta_k} x_k, e^{-i\theta_k} y_k, z_k) \ k=1, 
\dots, m.$$

If we restrict to the open set $U,$ this action is Hamiltonian with moment 
map $f: U \to \rr^m,$ 
$$ f(\{ (x_k, y_k, z_k); k=1 \dots m \}) = (|x_1|^2 - |y_1|^2, \dots, 
|x_m|^2 - |y_m|^2).$$ 

Furthermore, an easy computation shows that the naive parallel transport 
vector fields for $\chi: (U, \tilde \Omega) \to \sym^{2m}_0(\cc)$ are 
invariant with respect to the torus action and $df$ vanishes on them. 
Since $\kk(\bar \delta) \subset \yy_{m-1, \bar \tau}$ is invariant under the 
action and lies in $f^{-1}(0),$ using Lemma~\ref{rvc} we get that the 
vanishing cycle in $\yy_{m, \tau''_s}$ is also invariant and part of 
$f^{-1}(0).$ We also know that it is diffeomorphic to $(S^2)^m.$ It 
follows that it must be of the form $\kk(\delta''_s),$ where $\delta''_s$ is a 
matching of the points of $\tau''_s.$ Since it must lie close to 
$\kk(\bar \delta)$ for $s$ small, its isotopy class is determined uniquely. 
We can take $\delta''_s$ to consist of the linear path from $-s$ to $s$ 
together with the paths in $\bar \delta.$

The next step is moving $\kk(\delta''_s)$ from $\yy_{m, \tau''_s}$ to 
$\yy_{m, \tau'_s}$ by rescaled parallel transport along a linear path 
$\zeta: [0,1] \to \sym^{2m}_0(\cc).$ We choose a partition $0=t_0 < t_1 < 
\dots < t_N = 1$ of the interval $[0,1]$ and consider the corresponding 
partition of $\zeta.$ If each piece is sufficiently small, we can apply 
naive parallel transport to the Lagrangian so that we again know that it 
is invariant under the $T^m$ action and lies in $f^{-1}(0).$ Its isotopy 
class is determined uniquely as before. In particular, we can isotope it 
into $K$ of the matching consisting of the linear path between the first 
two coordinates and the respective translation of $\bar \delta$ joining 
the others. Then we continue the process using naive parallel transport. 
If the partition is sufficiently fine, everything is kept inside $U.$ By 
the discussion at the end of section~\ref{sec:kal}, the result is the same 
as that of the rescaled parallel transport.

We arrived at some Lagrangian in $\yy_{m, \gamma(1-s)}.$ We choose a fine 
partition of the path $\gamma: [0, 1-s] \to \sym^{2m}_0(\cc)$ and use 
naive parallel transport and isotopies to move the Lagrangian into 
$\yy_{m,\gamma(0)}$ such that at each step we have $\kk$ of some matching 
and everything is kept inside $U.$ The isotopy classes of the matchings 
are uniquely determined, and at the end we get $\kk(\delta) \subset \yy_{m, 
\tau}$ as desired.

There is one caveat about our inductive argument. Strictly speaking, the
inductive hypothesis gives us the Lagrangian $\kk(\bar \delta)$ in
$\yy_{m-1, \bar \tau}$ with the K\"ahler form $\tilde \Omega_{m-1}.$
However, as explained in section~\ref{sec:res}, we would like to use the
restriction of the form $\tilde \Omega = \tilde \Omega_m$ instead. Let us 
interpolate linearly between the two forms by letting $\omega(t)$ be the 
restriction
of $(1-t)\tilde \Omega_{m-1} + t \tilde \Omega$ to $\yy_{m-1, \bar \tau},$
for $t \in [0,1].$ We get a corresponding family of Lagrangians $\ell(t)$
in $(\yy_{m-1, \bar \tau}; \omega(t))$ for all $t,$ with $\ell(0) = \kk 
(\bar \delta).$ Note that $\omega(1)$
is of the form $(\ref{form2})$ on a big open set $U'$ that can be assumed
to contain $\kk(\bar \delta).$ Furthermore, $\kk(\bar \delta)$ and $\kk$ of
other matchings in $U'$ are Lagrangians for $\omega(t)$ for all $t.$ Also,
the torus action is still Hamiltonian with the same moment map for all
$\omega(t).$ 

We break the deformation into small pieces. On the first piece $[0,
\epsilon] \subset [0,1],$ Moser's lemma gives a family of symplectic
embeddings $\phi_{t_1, t_2}$ from a neighborhood of $\kk(\bar \delta)$ with
the form $\omega_{t_1}$ into $U'$ with the form $\omega_{t_2},$ for any $0
\leq t_1 \leq t_2 \leq \epsilon.$ Using the moment map as before we get
that $\phi_{0, t}(\ell(0))$ must be of the form $\kk$ of some crossingless
matching close to $\bar \delta$ for all $t\in [0, \epsilon].$ In
particular $\phi_{0, \epsilon} (\ell(0))$ is Lagrangian isotopic to
$K(\bar \delta)$ with respect to $\omega_{\epsilon}.$ On the other hand,
$\phi_{0, \epsilon}(\ell(0))$ is isotopic to $\ell(\epsilon)$ via the
isotopy $\phi_{t, \epsilon} (\ell(t)),$ $t \in [0, \epsilon].$ It follows
that $\ell(\epsilon)$ is isotopic to $\kk(\bar \delta)$ and by iterating
this procedure we get that the same is true for all $\ell(t), t\in [0,1].$
Thus we can safely start the inductive step with $\omega(1)$ rather than
$\omega(0).$ $\hfill \fin$

\medskip 

Theorem~\ref{lagr} is now a direct consequence of (\ref{hfl}) and
Proposition~\ref{super}. 

\begin {remark} We assumed that the variety $\yy_{m, \tau} \subset \ss_m$
came with the restriction of the K\"ahler form $\tilde \Omega.$ In fact,
since $\kk = \kk(\alpha), \kk'= \kk(\beta)$ are in the chosen open set $U,$ by doing
linear interpolation we see that the Floer cohomology groups are the same
for any exact K\"ahler form that has the form (\ref{form}) on $U \cap
\yy_{m, \tau}.$ \end {remark}

\section {Bigelow's definition of the Jones polynomial}
\label {sec:bigj}

The Jones polynomial \cite{J} is an invariant of oriented links in 
$S^3.$ It takes a link $L$ to a Laurent polynomial $V_L(t) \in 
\zz[t^{\pm1/2}]$ and is usually defined by the normalization 
$V_{\text{unknot}}=1$ and the skein relation:
$$ t^{-1}V_{L_+}(t) - tV_{L_-}(t) = (t^{1/2} - t^{-1/2}) V_{L_0}(t).$$

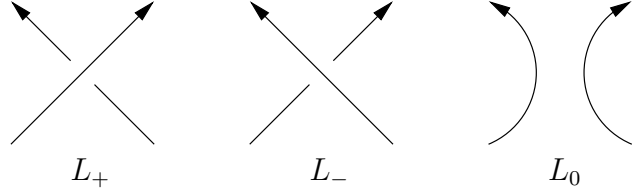
\begin {figure}
\begin {center}
\input {skein.pstex_t}
\caption {The links in the skein relation.}
\label {fig:skein}
\end {center}
\end {figure}

Here $L_+, L_-$ and $L_0$ are links that are identical except in a ball, where 
they look as in Figure~\ref{fig:skein}.

The normalization and the skein relation completely determine the polynomial. 
In this paper we will work mainly with a different normalization. 
Specifically, we make the change of variable $q = -t^{1/2}$ and set
\begin {equation}
\label {ujones}
 J_L(q) = (q + q^{-1}) \cdot V_L 
\end {equation}

This is usually called the {\it unnormalized Jones polynomial.} For
example, the normalized Jones polynomial of the left-handed trefoil is
$V_L = t^{-1} + t^{-3} - t^{-4}$ and the unnormalized one $J_L = q^{-1} +
q^{-3} + q^{-5} - q^{-9}.$

\subsection {Bigelow's picture.} Bigelow \cite{B} gave a different definition 
of the Jones polynomial. His construction is basically a reinterpretation of 
the work of Lawrence \cite{L}, with the resulting picture made more concrete 
and geometric. 

We alert the reader to the fact that Bigelow's conventions are 
different from ours, because he considers braids to act on the punctured plane 
on the left rather than on the right. Our conventions are more in line with 
Seidel and Smith's paper \cite{SS}. Correspondingly, the gradings 
$Q$ and $T$ defined below are in fact minus the gradings in Bigelow's paper.

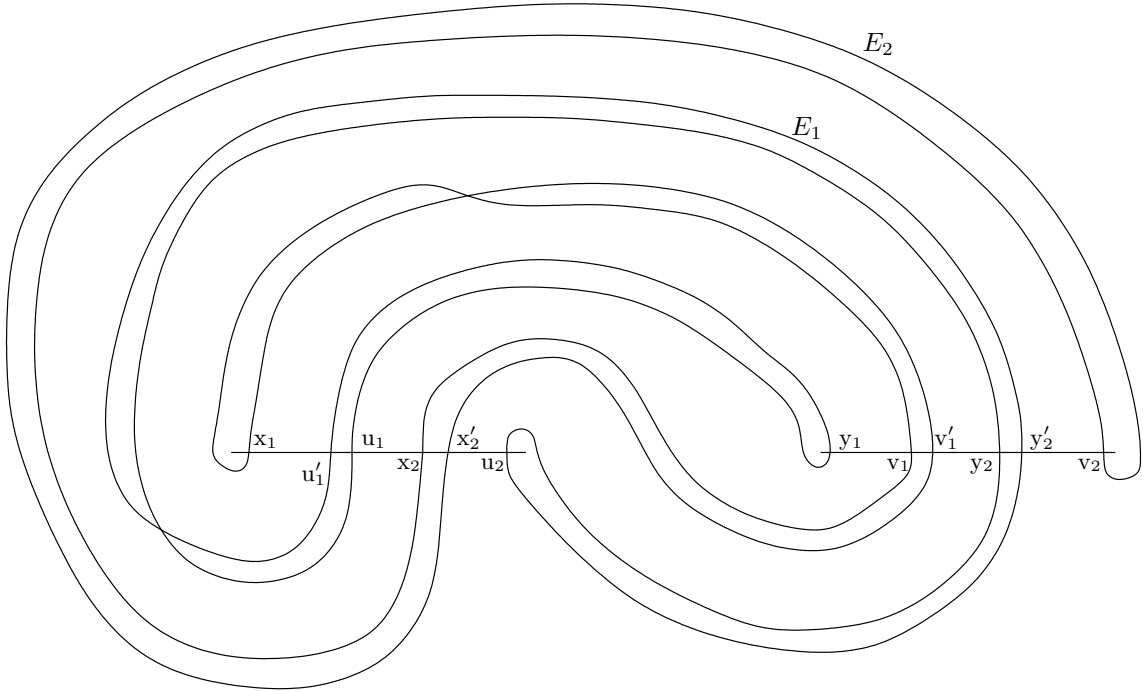
\begin {figure}
\begin {center}
\input {fig8.pstex_t}
\caption {Replacing the beta curves with figure-eights in a flattened braid 
diagram.}
\label {fig:fig8}
\end {center}
\end {figure}

His construction starts with an arbitrary plat representation of the link,
but for our purposes it suffices to consider braid diagrams as in
section~\ref{sec:braid}. We represent $L$ as the closure of $b \in Br_m$
and consider the resulting flattened braid diagram as in
Figure~\ref{fig:flat} for the trefoil. We replace every beta curve with a
figure-eight immersed loop running around it. This is shown in
Figure~\ref{fig:fig8} for the trefoil in Figure~\ref{fig:flat}. We denote
by $E_i$ be the figure-eight going around the curve $\beta_i$ for $i=1,
\dots, m.$

Consider a big disk $D \subset \cc$ containing all the $\alpha, \beta$ and 
$E$ curves. Let $D^*$ be the punctured disk 
$$ D^* = D - \{\mu_1, \dots, \mu_{2m} \}. $$

We can think of the braid group $Br_{2m} = \pi_1 (\conf^{2m}(D))$ as 
the mapping class group of $D^*.$ Thus $b \times 1^m$ induces a boundary 
preserving homeomorphism of $D^*.$ Because of the $1^m$ factor, this 
homeomorphism can in fact be extended over the even puncture points 
$\mu_2, \mu_4, \dots, \mu_{2m}.$

We denote by $\alpha^o$ the alpha curves with the endpoints removed. Then
$$ M_1 = \alpha^o_1 \times \dots \times \alpha^o_m \ \text { and } \ M_2 = 
E_1 \times \dots \times E_m $$
are two manifolds in $\ccc= \conf^m(D^*),$ the first one embedded, the 
second immersed. The Jones polynomial will come up as a graded 
intersection number of the lifts of $M_1$ and $M_2$ in a certain covering 
of $\ccc,$ which we now proceed to define. 

The composition of the map induced by the inclusion $D^* \hookrightarrow 
D$ with the natural abelianization map of the braid group is a 
homomorphism:
$$ \Phi_1 : \pi_1(\conf^m(D^*)) \to  \pi_1(\conf^m(D)) \cong Br_m \to 
\zz.$$

There is also an inclusion $\conf^m(D^*) \hookrightarrow \conf^{3m}(D)$ 
obtained by adding the $2m$ puncture points $\mu_1, \dots, \mu_{2m}.$ We 
get a corresponding homomorphism
$$ \Phi_2: \pi_1(\conf^m(D^*)) \to  \pi_1(\conf^{3m}(D)) \cong Br_{3m} \to
\zz.$$

It is easy to check that the image of $\Phi_1 - \Phi_2$ lies in $2\zz.$ 
The homomorphism
$$ \Phi = \bigl( \frac{1}{2}(\Phi_2- \Phi_1), \Phi_1 \bigr): \pi_1(\ccc) 
\to \zz \oplus \zz $$
is surjective and determines a $\zz \oplus \zz$ covering $\tilde \ccc$ of 
$\ccc.$ Let us denote by $\phi[a,b]: \tilde \ccc \to \tilde \ccc$ the 
covering transformation corresponding to the element $(a, b) \in \zz 
\oplus \zz.$

A more intuitive interpretation of the morphism $\Phi$ is the following. A
loop $\ell$ in $\ccc$ can be represented as a set of $m$ arcs in $D^*.$
The image of $\ell$ in the first $\zz$ factor counts the total winding
number of the $m$ arcs around the puncture points. The image in the second
$\zz$ factor counts twice the winding number of the arcs around each
other, which is basically the linking number of $\ell$ with the diagonal
in $\sym^m(D^*) \hookrightarrow \sym^m(D).$

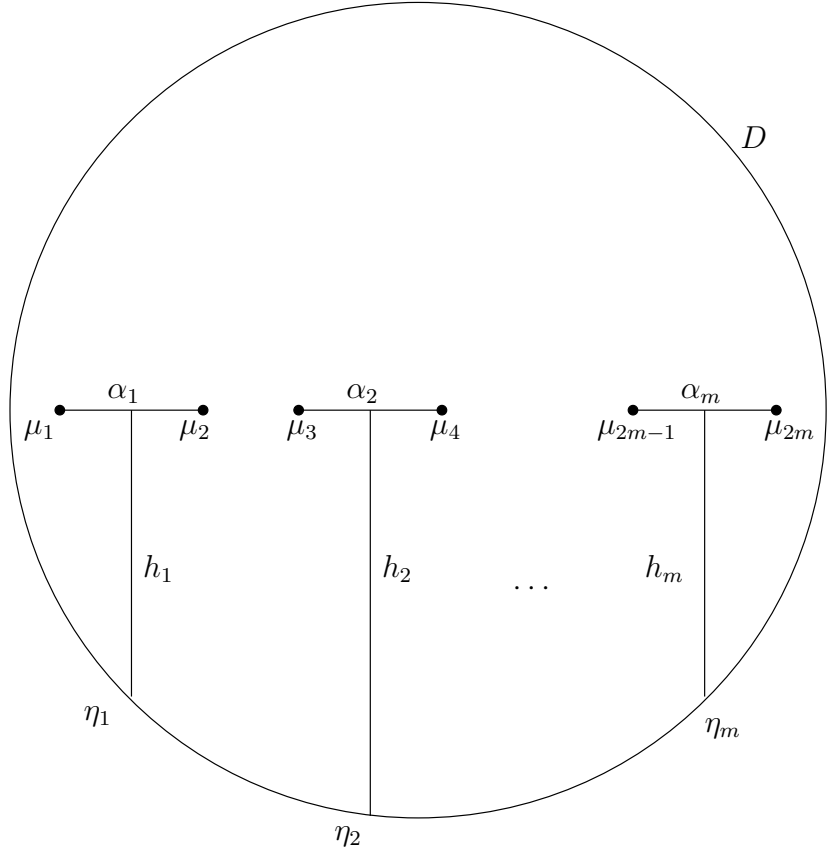
\begin {figure}
\begin {center}
\input {handles.pstex_t}
\caption {Handles for the alpha curves.}
\label {fig:handles}
\end {center}
\end {figure}

In order to lift the manifolds $M_1$ and $M_2$ from $\ccc$ to $\tilde 
\ccc,$ we need to specify some basepoints. We do this by attaching handles 
to each alpha curve. A {\it handle} $h_i$ is a segment in the lower half 
plane 
that starts from a point $\eta_i$ on the boundary of $D$ and ends in the 
midpoint of $\alpha_i,$ as shown in Figure~\ref{fig:handles}.

We take $\eta = (\eta_1, \dots, \eta_m)$ as a basepoint in $\ccc$ and 
also fix a lift $\tilde \eta \in \tilde \ccc$ of $\eta.$ The collection of 
all handles $h= (h_1, \dots, h_m)$ describes a path in $\ccc$ from $\eta$ 
to a point $m \in M_1.$ We lift $h$ to a path in $\tilde \ccc$ from 
$\tilde \eta$ to a point $\tilde m$ in the preimage of $m.$ Let $\tilde 
M_1 \subset \tilde \ccc$ be the lift of $M_1$ which contains 
$\tilde m.$ 

Similarly we can defined a lift $\tilde M_2$ of the $m$-torus $M_2.$
Consider the images of the handles under the homeomorphism of the disk
determined by $b \times 1^m$ which takes the alpha curves to the beta
curves. This gives a handle set for the beta curves that can be turned
into one for the figure eights and used to specify the lift of $M_2.$

There are well-defined algebraic intersection numbers $(\phi[a, b]\tilde 
M_1, \tilde M_2) \in \zz$ between the translates of $\tilde M_1$ and 
$\tilde M_2.$ We denote by $w$ the writhe of the braid $b$ as in 
Theorem~\ref{theoss}. The following result is proved in \cite{B}:

\begin {theorem}
The unnormalized Jones polynomial of the link $L$ can be expressed as
\begin {equation}
\label {jones}
J_L(q) = (-1)^m q^{m+w} \sum_{a, b \in \zz} (-1)^b q^{2(b-a)} (\phi[a, 
b]\tilde M_1, \tilde M_2).
\end {equation}
\end {theorem}

\subsection {The Bigelow generators.}
\label {sec:bigg}

Looking at the formula (\ref{jones}), it is clear that each intersection 
point between $M_1$ and $M_2$ contributes only once to a certain 
coefficient $j_n$ in $J_L(q)=\sum_{n \in \zz} j_n q^n.$ 

\begin {definition} 
\label {defy}
The elements of $\gz = M_1 \cap M_2$ are called {\bf
Bigelow generators.} The integer $n$ corresponding to a Bigelow generator
is called its {\bf Jones grading.} \end {definition}

The formula (\ref{jones}) can be rewritten as:
\begin {equation}
\label {ujon}
 J_L(q) = \sum_{\gamma \in \gz} \sigma(\gamma) q^{J(\gamma)},
\end {equation}
where $\sigma$ is a sign function $\sigma : \gz \to \{\pm 1\}.$ 

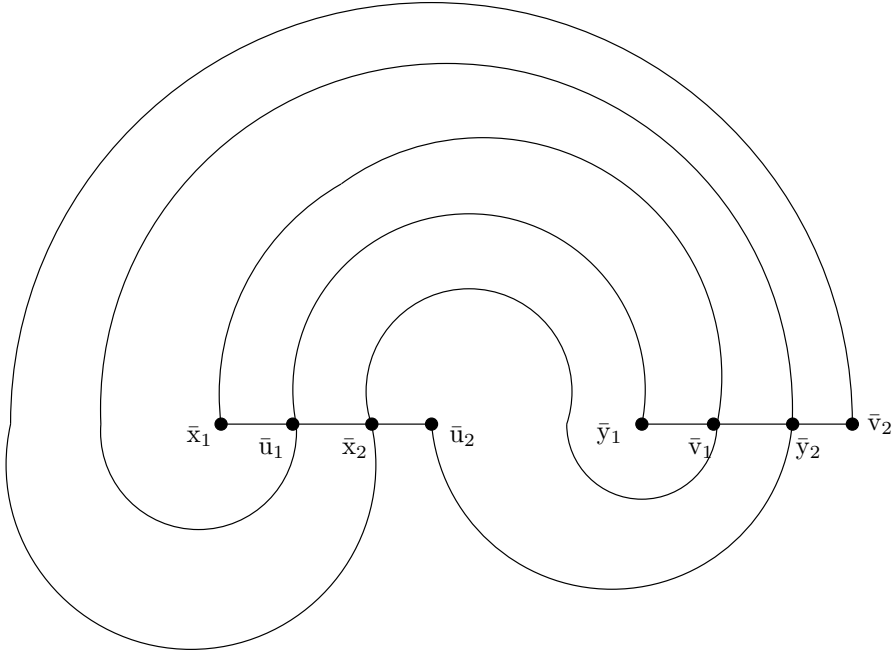
\begin {figure}
\begin {center}
\input {flat2.pstex_t}
\caption {Intersections of the alpha and beta curves.}
\label {fig:flat2}
\end {center}
\end {figure}

Let us explain more carefully the structure of the set $\gz$ and fix some
notation. We started with a flattened braid diagram for the link $L.$ For
the sake of concreteness, we will always write down what happens for the
case of the left-handed trefoil in Figure~\ref{fig:flat}. Let $\bar \ze$
be the set of all intersection points between an alpha curve and a beta
curve. Figure~\ref{fig:flat2} shows them for the trefoil:
$$ \bar \ze = \{\bar \xh_1, \bar \xh_2, \bar \yh_1, \bar \yh_2, \bar 
\uh_1, \bar 
\uh_2, \bar \vh_1, \bar \vh_2 \}.$$

Note that the set of puncture points $\tau = \{\mu_1, \dots, \mu_{2m} \}$
is in fact a subset of $\bar \ze.$ As explained in the introduction, we
construct a set $\ze$ from $\bar \ze$ by doubling the points of $\bar \ze
- \tau.$ That is, for every $x \in \tau$ we introduce an element $e_x \in 
\ze$ and for every $x \in \ze - \tau$ we introduce two elements $e_x, e'_x 
\in \ze.$ There is a natural map 
\begin {equation}
\label {ff}
f: \ze \to \bar \ze
\end {equation}
which takes $e_x$ and $e'_x$ to $x.$

The set $\ze$ can be thought of as the set of intersection points between 
the alpha curves and the figure eights $E_i.$ Indeed, near each $x \in 
\alpha_i \cap \beta_j$ there is one point $e_x \in \alpha_i \cap E_j$ when 
$x \in \tau$ and two points $e_x, e_x'$ otherwise. We distinguish between 
$e_x$ and $e_x'$ by requiring that the loop in the plane starting at 
$e_x,$ following the part of the figure eight with no self-intersections 
around a puncture point up to $e_x'$ and then going back to $e_x$ along 
the alpha curve has winding number $+1$ with the puncture.

In our trefoil example (Figure~\ref{fig:fig8}) we have
$$ \ze = \{\xh_1, \xh_2, \xh_2', \yh_1, \yh_2, \yh_2', \uh_1, \uh_1', 
\uh_2, \vh_1, \vh_1', \vh_2 \}$$
and $f(\xh_1) = \bar \xh_1, f(\xh_2) = f(\xh_2') = \bar \xh_2,$ etc.

We define maps $\bar A, \bar B: \bar \ze \to \{1, 2, \dots, m\}$ by taking 
an intersection point in $\alpha_i \cap \beta_j$ to $(i, j).$ The set 
$$\bar \gz = (\alpha_1 \times \dots \times \alpha_m) \cap (\beta_1 \times 
\dots \times \beta_m) \subset \conf^m(\cc)$$ 
consists of unordered $m$-tuples $(x_1, \dots, x_m)$ of elements of $\bar 
\ze$ with $A(x_i) \neq A(x_j)$ and $B(x_i) \neq B(x_j)$ for all $i \neq 
j.$

Similarly, we can consider the compositions $A= \bar A \circ f, B = \bar B 
\circ f: \ze \to \{1, 2, \dots, m \}.$ The elements of the resulting set 
$\gz$ are exactly the Bigelow generators. The map (\ref{ff}) induces  
natural map from $\gz$ to $\bar \gz,$ which we still denote by $f.$

In the trefoil example we have 18 Bigelow generators:
$$\gz = \{\ \xh_1\yh_1,\ \xh_1\yh_2,\ \xh_1\yh_2',\ \xh_2\yh_1,\ 
\xh_2\yh_2,\ \xh_2\yh_2',\  
\xh_2'\yh_1,\ \xh_2'\yh_2,\ \xh_2'\yh_2'\ \} $$ 
$$ \cup \ \{ \ \uh_1\vh_1,\ \uh_1\vh_1',\ \uh_1\vh_2,\  \uh_1'\vh_1,\ 
\uh_1'\vh_1',\  
\uh_1'\vh_2,\  \uh_2\vh_1,\ \uh_2\vh_1',\ \uh_2\vh_2\ \}. $$

\subsection {Gradings.}
\label {sec:grad}

An absolute grading on the Bigelow generators is a map $F: \gz \to \zz.$ 
An affine grading is an equivalence class of absolute gradings under 
the equivalence relation $F_1 \sim F_2$ if $F_1 = F_2 + k$ for some $k \in 
\zz.$ The following two definitions make sense for both absolute and 
affine gradings.

\begin {definition}
\label {deff}
A grading $F: \gz \to \zz$ is called {\bf additive} if it comes as a
summation of gradings on $\ze,$ i.e. there is a grading $F^*: \ze \to \zz$
and $F$ is defined by $F(e_1,\dots, e_m) = F^*(e_1) + \dots + F^*(e_m).$
\end {definition}

\begin {definition} 
\label {defs} A grading $F: \gz \to \zz$ is called {\bf stable} if
it can be expressed as a composition $\bar F \circ f$ for a grading $\bar
F: \bar \gz \to \zz.$ \end {definition}

Our interest lies in the absolute Jones grading from 
Definition~\ref{defy}, which we denote by $J: \gz \to \zz.$ In section 3
of \cite{B} Bigelow indicated how to compute $J$ from the flattened braid
diagram.  If we look at equation (\ref{jones}), we see that affinely $J$
is minus twice the difference of two gradings $Q$ and $T$ which correspond
to the pair $(a, b) \in \zz \oplus \zz$ describing the covering
transformation. To fix $Q, T$ and $J$ as absolute gradings we use a
distinguished element $\nu \in \gz.$ Specifically, we let $\nu_1, \dots,
\nu_m \in \ze$ be the preimages of the $m$ even puncture points $\mu_2,
\dots, \mu_{2m} \in \bar \ze$ under $f,$ and set $\nu = (\nu_1, \dots,
\nu_m).$ In our trefoil example, the distinguished element is $\nu = 
(\uh_2 \vh_2).$

The grading $Q$ is defined to be additive. Let us explain what is the
corresponding $Q^*:  \ze \to \zz.$ We start by setting $Q^*(\nu_k)= 0$ for
all $k=1, \dots, m.$ Then, for every $e \in \alpha_i \cap E_j \subset
\ze,$ we set $Q^*(e)$ to be the total winding number around the $2m$
puncture points of the following loop in $\cc$: we start at $\nu_j \in
E_j,$ we follow one branch of the figure eight $E_j$ up to $e,$ then we go
along $\alpha_i$ to its midpoint, we follow the handle $h_i$ down to
$\eta_i,$ then move to $\eta_j$ along the lower half of the boundary
$\partial D.$ Next we go up the handle $h_j$ to the midpoint of $\alpha_j$
and back to $\nu_j$ along $\alpha_j.$ This completes the loop, and since
every figure eight has total winding number $+1-1=0$ with the punctures, it
doesn't really matter which route on $E_j$ we followed from $\nu_j$ to
$e.$

\begin {figure}
\begin {center}
\input {loop.pstex_t}
\caption {Loop used to compute the $Q^*$ grading for $x_1.$}
\label {fig:loop}
\end {center}
\end {figure}
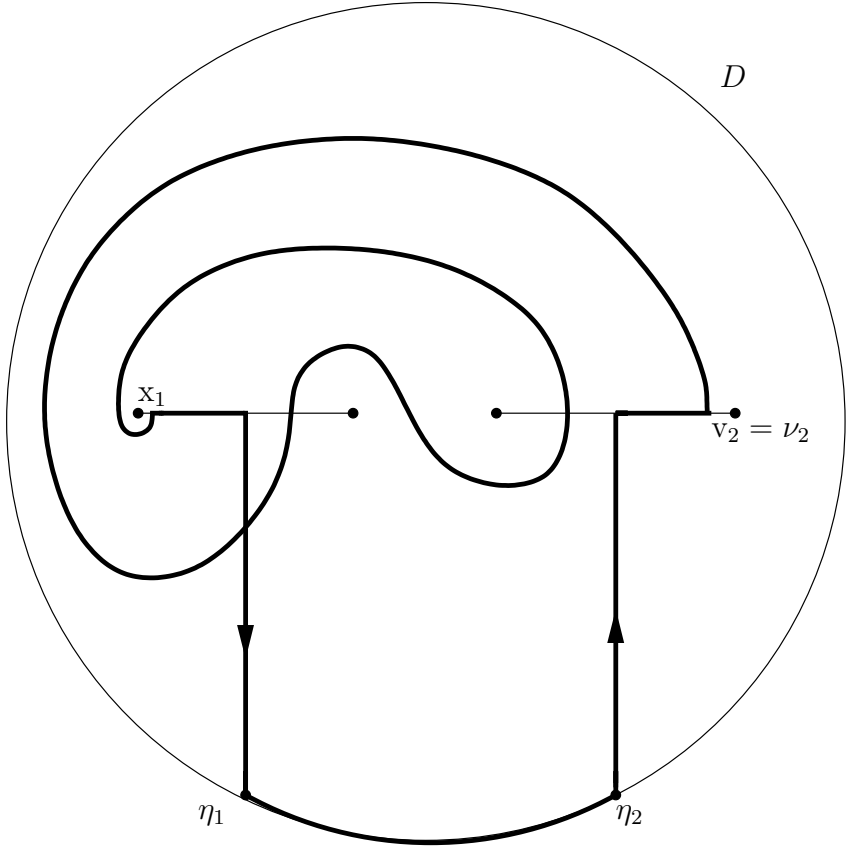

Figure~\ref{fig:loop} shows the corresponding loop for $\nu_2 = \vh_2$ 
and $e=\xh_1$ in the trefoil example. The loop has total winding number 
five
around the puncture points, so $Q^*(\xh_1) = 5.$ If we do similar
computations for all the intersection points we find that $Q^*(\yh_2) = 
-1;
\ Q^*(\vh_2) = Q^*(\yh_2') = Q^*(\uh_2) = 0;$ 
 $Q^*(\yh_1) = 1;  \ Q^*(\uh_1) = Q^*(\vh_1) = 2;$ 
 $Q^*(\vh_1') = Q^*(\xh_2) = Q^*(\uh_1') =3;$ 
$Q^*(\xh_2') = 4; \  Q^*(\xh_1) = 5.$

This gives the $Q$ grading on Bigelow generators:

\medskip
\begin {center}
\begin {tabular} {|llllll|}
\hline
$Q=6$ &:&  $\uh_1'\vh_1'$ & $\xh_1\yh_1$ & &  \\
$Q=5$ &:&  $\uh_1\vh_1'$ & $\uh_1'\vh_1$ & $\xh_2'\yh_1$ & $\xh_1\yh_2'$  
\\
$Q=4$ &:&  $\uh_1\vh_1$ &  $\xh_2\yh_1$ & $\xh_2'\yh_2'$ &  $\xh_1 \yh_2$   
\\
$Q=3$ &:&  $\uh_2\vh_1'$ & $\uh_1'\vh_2$ & $\xh_2\yh_2'$ & $\xh_2'\yh_2$    
\\
$Q=2$ &:&  $\uh_2\vh_1$ & $\uh_1\vh_2$ & $\xh_2\yh_2$ &  \\
$Q=1$ &:&  $-$  & & &  \\
$Q=0$ &:& $\uh_2\vh_2$ & & & \\    
\hline
\end {tabular}
\end {center}
\medskip

The second grading $T$ on $\gz$ is not additive but is stable. The
corresponding $\bar T: \bar \gz \to \zz$ is defined as follows. Let $x =
(x_1, \dots, x_m)$ and $y= (y_1, \dots, y_m)$ be two elements of $\bar
\gz.$ Each of them is formed from $m$ intersections of the $\alpha$ and
$\beta$ curves. To compute the difference $\bar T(x) - \bar T(y)$ we
consider the following loop in in $\sym^m(D) \subset \sym^m(\cc) \cong
\cc^m.$ We start at $x,$ go along the alpha curves to $y,$ then go back to
$x$ along the beta curves. Then $\bar T(x) - \bar T(y)$ is the linking
number of this loop with the diagonal $\Delta$ in $\sym^m(\cc),$ where
$\sym^m(\cc)$ and the diagonal are taken with their complex orientations.
Once we know $\bar T,$ we get $T$ by composing with the natural map $f:  
\gz \to \bar \gz.$ To fix $T$ as an absolute grading, we set $T =
0$ on the distinguished element $\nu.$

In practice, the linking number with the diagonal records the twisting of
the points around each other.  For example, a half twist such as the one
between $\bar \xh_1 \bar \yh_1$ and $\bar \uh_1 \bar \vh_1$ in
Figure~\ref{fig:flat2} gives a difference of $1$ in their $\bar T$
gradings. In fact there are enough half twists of this form in
Figure~\ref{fig:flat2} to relate any two elements of $\bar \gz.$ Thus it
is easy to write down the $T$ grading of the Bigelow generators by
calculating $\bar T$ and then composing with $f.$ We obtain:

\medskip
\begin {center}
\begin {tabular} {|llllllllll|}
\hline
$T=3$ &:&  $\xh_1\yh_1$ & &  & & & & & \\
$T=2$ &:&  $\uh_1\vh_1'$ & $\uh_1'\vh_1$ & $\uh_1\vh_1'$ & $\uh_1'\vh_1'$ 
& & & & \\
$T=1$ &:&  $\xh_2\yh_1$ & $\xh_2'\yh_1$  &  $\xh_1 \yh_2$ & $\xh_1\yh_2'$ 
& $\xh_2\yh_2$ & $\xh_2'\yh_2$ & $\xh_2\yh_2'$ & $\xh_2'\yh_2'$  \\
$T=0$ &:&  $\uh_2\vh_1$ & $\uh_2\vh_1'$ & $\uh_1\vh_2$ & $\uh_1'\vh_2$ & 
$\uh_2\vh_2$ & & & \\
\hline
\end {tabular}
\end {center}
\medskip
 
Looking back at (\ref{jones}), we get that the $T$ and $Q$ gradings 
determine $J$ by the formula $ J = 2(T-Q) + m + w. $ Doing this 
computation for the trefoil, and taking into account the factor 
$m + w = 2-3 = -1,$ we get:

\medskip
\begin {center}
\begin {tabular} {|llllllllll|}
\hline
$J=-1$ &:&  $\uh_2\vh_2$ & & & & & & & \\
$J=-3$ &:&  $\xh_2\yh_2$ & & & & & & & \\
$J=-5$ &:&  $\uh_1\vh_1$ & $\uh_1\vh_2$ & $\uh_2\vh_1$ & $\xh_2'\yh_2$ &  
$\xh_2 \yh_2'$ & & &  \\
$J=-7$ &:&  $\uh_2\vh_1'$ & $\uh_1'\vh_2$ & $\uh_1\vh_1'$ & $\uh_1'\vh_1$ 
& $\xh_1\yh_1$ & $\xh_1\yh_2$ & $\xh_2'\yh_2'$ & $\xh_2\yh_1$    \\
$J=-9$ &:& $\uh_1'\vh_1'$ & $\xh_1\yh_2'$ & $\xh_2'\yh_1$ & & & & & \\
\hline
\end {tabular}
\end {center}
\medskip

\subsection {The sign.} \label {sec:signy} Now that we understand $J,$ let us turn our
attention to the sign $\sigma : \gz \to \{\pm 1\}$ appearing in the formula (\ref{ujon}).
Looking at (\ref{jones}), we see that $\sigma$ has a contribution $(-1)^{b + m}$ from the
$T$ grading and the number of strands, and a contribution from comparing the orientations
in the intersection product. The sign of an intersection point $\gamma \in \gz$ is also
composed of two parts $\sigma'$ and $\sigma''.$ The sign $\sigma'= (-1)^b$ represents the parity of
the permutation of $\{1, \dots, m\}$ associated to the Bigelow generator $\gamma.$ The
second factor $\sigma''$ is obtained by giving orientations to the figure eights, checking
whether the figure eight hits the alpha curve from above or from below at each point in
$\ze$, and then multiplying the local intersection signs of all the $m$ points that form
$\gamma.$ In \cite{B} the alpha curves are oriented from left to right, and the figure
eights are oriented so that they hit the alpha curves from above at the endpoints $\nu_j.$
This means that $\sigma''(\nu) = (-1)^m.$

Putting these together we get the formula:
\begin {equation}
\label {signy}
 \sigma (\gamma) = (-1)^{b+m} \sigma'(\gamma) \sigma''(\gamma) = (-1)^m 
\sigma''(\gamma).
\end {equation}
Note that this overall sign $\sigma$ is always $+$ on the distinguished
element $\nu.$

\medskip
\section {Floer cohomology and a bijective correspondence}
\label {sec:flob}

In this section we prove Theorem~\ref{bijection} and define a new grading
on the Bigelow generators. Some care needs to be taken about the 
transversality of intersections.

\subsection {Clean intersections and Floer cochains.} We work with the two 
Lagrangians
$$ \kk = \Sigma_{\alpha_1} \times \Sigma_{\alpha_2} \times \dots
\times \Sigma_{\alpha_m} \ , \
 \kk' = \Sigma_{\beta_1} \times \Sigma_{\beta_2} \times \dots
\times \Sigma_{\beta_m} \subset Y=\yy_{m, \tau}$$
from the statement of Theorem~\ref{lagr}. We would like to describe a set 
of generators for the Floer cochain $CF^*(\kk, \kk').$ 

However, the intersection $\nn = \kk \cap \kk'$ is not transverse, and we
will need to isotope one of the Lagrangians in order to achieve
transversality. Let us understand the structure of $\nn.$ The $2$-spheres 
$\Sigma_{\alpha_i}$ are Lagrangians in $S_{\tau}=(u^2 + v^2 + P_{\tau}(z) 
= 0)$ and they map to the corresponding alpha curve under the projection 
to the $z \in \cc$ coordinate. Similarly $\Sigma_{\beta_j}$ map to the 
beta curves. Therefore the intersection $\Sigma_{\alpha_i} \cap 
\Sigma_{\beta_j}$ is a disjoint union of circles and points: each 
point $x \in \alpha_i \cap \beta_j \subset \bar \ze$ contributes a point 
$e_x$ to  $\Sigma_{\alpha_i} \cap \Sigma_{\beta_j}$ if $x = \mu_k$ for 
some $k$ and a circle $S^1_x$ otherwise. It follows that the intersection 
$\nn = \kk \cap \kk'$ consists of a disjoint union of tori $T^k = 
S^1 \times \dots \times S^1$ of various dimensions, one for each point 
in $\bar \ze.$ Usually not all tori are trivial, so the intersection is 
not transverse. 

In order to control the isotopy necessary for transversality, we observe
that $\nn = \kk \cap \kk'$ is a {\it clean intersection} in the sense of
Pozniak \cite{P}, i.e. $T\nn = (T\kk|_{\nn}) \cap (T\kk'|_{\nn}).$ A model
for isotoping a Lagrangian in a neighborhood $V$ of a clean intersection
is due to Weinstein \cite{We} and was used by Khovanov and Seidel in the
proof of Proposition~5.15 in \cite{KS}. Basically, if one starts with a
Morse-Smale function $g$ on $\nn,$ then $\kk'$ can be isotoped into a
Lagrangian $\kk''$ which is identical to $\kk'$ outside $V$ and intersects
$\kk \cap V$ exactly at the critical points of $g.$ Moreover, given two
critical points $x, y$ of $g,$ the difference in their cohomological
grading in $CF^*(\kk, \kk'')$ is the same as the difference in their Morse
indices on $\nn.$

We apply this to our case. Take the standard height function on every 
circle $S^1_x$ with two critical points: a minimum $e_x$ and a maximum 
$e_x'.$ Deform the Lagrangian spheres $\Sigma_{\beta_j}$ into 
$\Sigma''_{\beta_j}$ accordingly and set
$$ \kk'' =  \Sigma_{\beta_1}'' \times \Sigma_{\beta_2}'' \times \dots
\times \Sigma_{\beta_m}'' \subset Y=\yy_{m, \tau}.$$

This makes the intersection $\nn'' = \kk \cap \kk''$ transverse. The
intersection points in $\Sigma_{\alpha_i} \cap \Sigma''_{\beta_j}$ are now
all of the form $e_x, e_x',$ exactly in 1-to-1 correspondence with the
elements of $\ze.$ In turn this induces a 1-to-1 correspondence between
the points in $\nn''$ and the Bigelow generators. Thus we can think of
$CF^*(\kk, \kk') = CF^*(\kk, \kk'')$ as being generated by the elements of
$\gz.$ This proves Theorem~\ref{bijection}.

\subsection {The projective grading.} \label {sec:proj} Recall from section~\ref{sec:floer}
that the Floer cochain complex $CF^*(\kk, \kk')$ is absolutely $\zz$-graded. 

\begin {definition} We denote by $\tilde P$ the cohomological grading on Bigelow
generators, induced by the identification in Theorem~\ref{bijection}. We also set $P =
\tilde P - m-w.$ \end {definition} 

The renormalization by $m+w$ is made to account for the
shift in degree in Theorem~\ref{theoss}. Of course, $P$ and $\tilde P$ are identical as
affine gradings. We would like to understand $P$ in terms of a flattened braid diagram.

To define the absolute grading we started with complex volume forms
for the slice $\ss_m$ and the base $\sym^{2m}_0(\cc),$ which gave us a
family of forms on $\yy_{m, \tau}$ for $\tau \in \conf^{2m}_0(\cc).$ We
want to specify the volume form on the open sets $U_{m, \tau} \cap \yy_{m,
\tau}$ from Corollary~\ref{cor}, which we can identify with open subsets
(\ref{vmt}) in $\sym^m(S_{\tau} - \Delta)$ using the Hilbert-Chow morphism 
$\pi.$ 

We begin by considering the standard volume form $ \omega =
(dv \wedge dz)/2u$ on the hypersurface 
$S_{\tau} \in \cc^3$ given by the equation $u^2 + v^2 + P_{\tau}(z)=0.$ 
Then we use the following lemma, which is subsumed in the proof of 
Theorem 1.17 in \cite{N2}. The argument there is in fact due to Beauville 
\cite{Be}.

\begin {lemma}
\label {beau}
Let $\omega$ be a complex volume form on a smooth surface $S.$ Clearly 
$\omega^{m}$ is a complex volume form on the Cartesian product 
$S^{\times m}$ invariant under the 
action of the symmetric group, and therefore descends to a complex 
volume form $\omega_m$ on $\sym^m(S) - \Delta.$ 

Then there is an extension of the pullback $\pi^*\omega_m$ from
$\pi^{-1}(\sym^m(S) - \Delta)$ to a complex volume form on the whole
Hilbert scheme $\hi^m(S).$ \end {lemma}

Taking $S= S_{\tau}$ as $\tau$ varies over $\conf^{2m}_0(\cc)$ we get a
family of volume forms on the Hilbert schemes $\hi^m(S_{\tau}),$ and by
restriction one on the family of $\yy_{m, \tau} \subset \hi^m(S_{\tau}).$
We can use these volume forms to define the absolute grading. On the open
set (\ref{vmt}), in terms of the coordinates $u_j, v_j, z_j,$ the volume
form is
\begin {equation}
 \Theta = \prod_{j=1}^m \frac{dv_j \wedge dz_j}{2u_j}.
\end {equation}

At a point $x \in \kk = \Sigma_{\alpha_1} \times \dots \times
\Sigma_{\alpha_m}$ the coordinates $u_j, v_j, z_j$ satisfy $z_j = 
\alpha_j(t_j)$ for some $t_j \in [0,1]$ and $u_j, v_j \in 
\sqrt{-P_{\tau}(\alpha_j(t_j))} \rr.$ Hence the resulting square phase map 
(\ref{sqp}) on the Lagrangian $\kk$ is given by
\begin {equation}
\theta_{\kk}: \kk \to \cc/\rr_+, \ \theta_{\kk}(x) = \prod_{j=1}^m 
\frac{-P_{\tau}(\alpha_j(t_j)) \cdot 
\alpha'(t_j)^2}{-P_{\tau}(\alpha_j(t_j))} 
= \prod_{j=1}^m \alpha_j'(t_j)^2. 
\end {equation}

This is of course constant because the alpha curves are horizontal. On the 
other hand, the square phase map on $\kk'$ is 
\begin {equation}
\label {betas}
\theta_{\kk'}: \kk' \to \cc/\rr_+, \ \theta_{\kk'}(x) = \prod_{j=1}^m 
\beta_j'(t_j)^2,
\end {equation}
which is nontrivial. Equation~(\ref{betas}) and the standard additivity 
properties of the Morse index imply that the $\tilde P$ grading on the Bigelow 
generators is additive in the sense of Definition~\ref{deff} and can be 
computed from the flattened braid diagram. The 
corresponding grading $\tilde P^*: \ze \to \zz$ on the intersection points is 
described by the phase map $x=\beta(t) \to \beta'(t)^2$ on the beta curves 
in the plane. 

This is the projective grading considered by Khovanov and Seidel in 
\cite[section 3d]{KS}. If we take into account the difference in grading 
of $1$ between $e_x'$ and $e_x$ coming from the Morse-Smale function on 
the clean intersections, we get a very simple description of the affine 
grading $\tilde P^*=P^*$ in terms of Figure~\ref{fig:fig8}, where we replaced the 
beta curves by figure-eights. Specifically, we assume that all the figure 
eights intersect the alpha curves at $90^{\circ}$ angles. Each figure 
eight is the image of an immersion $\gamma_j: S^1 \to \cc.$ We get a map 
$$ \varepsilon_j : S^1 \to S^1, \  t \to \gamma_j'(t)^2/ | 
\gamma_j'(t)|^2.$$

This map has degree $0$ so it can be lifted to a real valued map $\tilde 
\varepsilon_j : S^1 \to \rr.$ We choose the lift so that $\tilde 
\varepsilon_j(t) = 0$ when $\gamma_j(t)$ is the puncture point $\mu_{2j}.$
Then $\varepsilon_j(t)$ is an integer whenever $\gamma_j(t) \in \gz,$ and 
that integer is its $\tilde P^*$ grading. 

In our trefoil example we have $\tilde P^*(\uh_2)=\tilde P^*(\vh_2) = \tilde P^*(\yh_2) =0;  
\ \tilde P^*(\yh_2') = 1; \ \tilde P^*(\xh_2) = \tilde P^*(\vh_1) = \tilde P^*(\uh_1) =
\tilde P^*(\yh_1) = 2;  \ \tilde P^*(\xh_2') = \tilde P^*(\vh_1') = \tilde P^*(\uh_1') = 3$
and $\tilde P^*(\xh_1)  =4.$

We can add these together and get the $\tilde P$ grading on $\zz.$ Taking into 
account the factor $m+w = 2-3 = -1$ we can list the $P$ grading on Bigelow 
generators:
                                                                                
\medskip
\begin {center}
\begin {tabular} {|llllll|}
\hline
$P=7$ &:&   $\uh_1'\vh_1'$ & $\xh_1\yh_1$ & &  \\
$P=6$ &:&  $\uh_1\vh_1'$ & $\uh_1'\vh_1$ & $\xh_2'\yh_1$ & $\xh_1\yh_2'$
\\
$P=5$ &:&  $\uh_1\vh_1$ &  $\xh_2\yh_1$ & $\xh_2'\yh_2'$ &  $\xh_1 \yh_2$
\\
$P=4$ &:&  $\uh_2\vh_1'$ & $\uh_1'\vh_2$ & $\xh_2\yh_2'$ & $\xh_2'\yh_2$
\\
$P=3$ &:&  $\uh_2\vh_1$ & $\uh_1\vh_2$ & $\xh_2\yh_2$ &  \\
$P=2$ &:&  $-$  & & &  \\
$P=1$ &:& $\uh_2\vh_2$ & & & \\
\hline
\end {tabular}
\end {center}

\begin {remark}
\label {sign}
In section \ref{sec:signy} we described the sign $\sigma = (-1)^m \sigma''$ 
appearing in front of a Bigelow generator in the formula (\ref{ujon}) for 
the Jones polynomial. It is now straightforward to see that: $$\sigma = (-1)^{\tilde P} 
= (-1)^{P+m+w}.$$ \end {remark}

\subsection {The $P-Q$ grading.} The $P$ and $Q$ gradings are additive but
not stable. If we look at the corresponding tables for our trefoil example
we may get the impression that $P$ and $Q$ are identical as affine
gradings.  This is just a coincidence. If we replaced the beta curve with
a figure eight on the right hand side of Figure~\ref{fig:unknot}, for
example, we would get several Bigelow generators for the unknot that have
the same $Q$ grading but different projective gradings.

Nevertheless, we can still say something about the difference $P-Q$ for 
any diagram:

\begin {lemma}
\label {pqs}
The grading $P - Q: \gz \to \zz$ is stable in the sense of 
Definition~\ref{defs}.
\end {lemma}

\noindent {\it Proof.} Since $P$ and $Q$ are additive it suffices to show 
that $(P^*-Q^*)(e_x) = (P^*-Q^*)(e_x')$ for every $x \in \ze$ which is not 
a puncture. But this follows from the fact that $P^*(e'_x) - P^*(e_x) = 
Q^*(e'_x) - Q^*(e_x) = 1$ for any such $x.$ $\hfill \square$

\subsection {Comparison with Khovanov cohomology.}

In the introduction we mentioned the conjecture made by Seidel and Smith 
about the equivalence between their theory and the cohomology theory of 
Khovanov \cite{Kh}:
\begin {equation}
\label {conjy}
 Kh_{symp}^k(\ka) \cong \bigoplus_{i-j=k} Kh^{i,j}(\ka) \hskip14pt (?)
\end {equation}

The grading $j$ in Khovanov cohomology is the ``Jones grading'' which 
describes which coefficient in the unnormalized Jones polynomial comes 
from the Euler characteristic:
$$ \sum_{i, j \in \zz} (-1)^i q^j \dim(Kh^{i,j} \otimes \qq) = J_L(q).$$ 

The grading $i$ describes the cohomological degree, while $k$ on the 
Seidel-Smith side of (\ref{conjy}) is the projective grading $P,$ which 
is supposed to correspond to $i-j.$

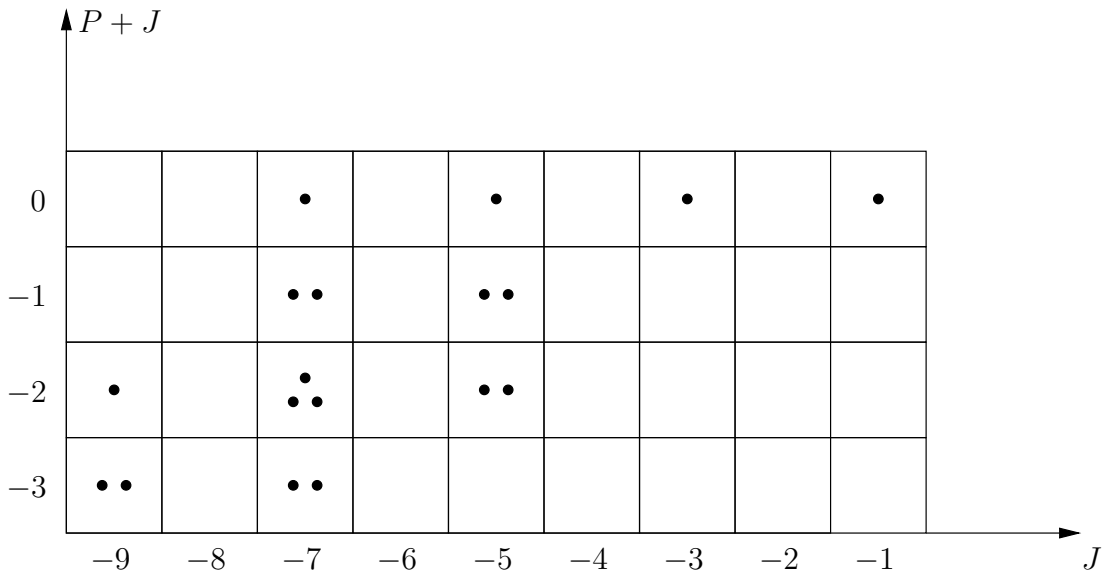
\begin {figure}
\begin {center}
\input {big.pstex_t}
\caption {Bigelow generators for the trefoil.}
\label {fig:big}
\end {center}
\end {figure}

Let us plot the Bigelow generators in our example for the left-handed 
trefoil, with the $J$ and $P+J$ grading on the axes. The result is 
Figure~\ref{fig:big}, where each dot represents a generator.  

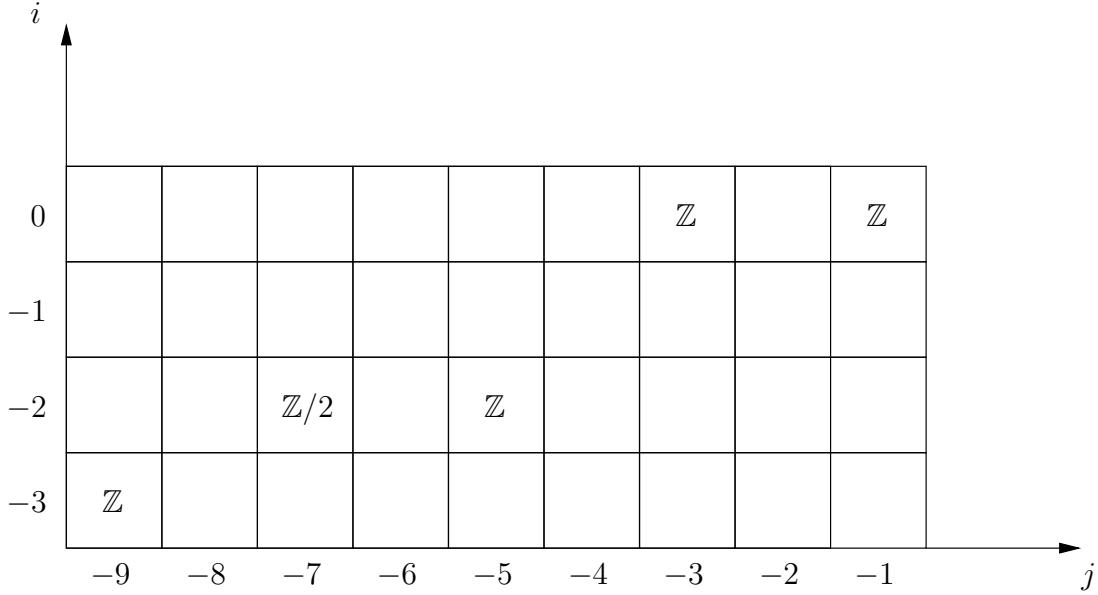
\begin {figure}
\begin {center}
\input {khov.pstex_t}
\caption {Khovanov cohomology of the trefoil.}
\label {fig:khov}
\end {center}
\end {figure}

Let us also plot the Khovanov cohomology of the trefoil, calculated in 
\cite[section 7]{Kh}. We get Figure~\ref{fig:khov}. 

Note that Conjecture (\ref{conjy}) was verified for the trefoil by Seidel 
and Smith in \cite{SS}. They computed $Kh^*_{symp}$ of the trefoil 
to be $(\zz, 0, \zz^2, 0, \zz/2, \zz)$ in degrees from $1$ to $6$ 
respectively, and $0$ in all other degrees.

Let us compare Figures \ref{fig:big} and \ref{fig:khov}. We know that the
Seidel-Smith cochain complex is generated by the dots in
Figure~\ref{fig:big}, and that the differentials have to increase the $P$
grading by $1.$ Although we do not deal with differentials in this paper,
let us observe that there is enough room for them to produce the abelian
groups in Figure~\ref{fig:khov}, even assuming that they preserve the $J$
grading. This lends support to the conjecture that the bigrading $(P, J)$
on Bigelow generators should descend to a bigrading on $Kh^*_{symp}$
similar to that on $Kh^*,$ with the Jones grading $J$ playing the role of
$j.$
 
\section {Heegaard Floer homology of the double branched cover}

Seidel and Smith (\cite{SS}, \cite{S}) have considered an interesting 
involution acting on the manifold $\yy_{m, \tau}$ and described how this 
could be used to relate their theory to the Heegaard Floer homology of 
$\dd(L),$ the double cover of $S^3$ branched over $L.$ In this section 
we explore this relation in terms of the Bigelow generators.

\subsection {An involution.} On the Milnor fiber $S_{\tau}$ given by the 
equation $u^2 + v^2  + P_{\tau}(z) = 0$ there is an involution 
\begin {equation}
\label {invo}
 (u, v, z) \to (u, -v, z).
\end {equation}

Its fixed point set $\sig_{\tau}$ is the affine complex curve in $\cc^2$ 
with equation $u^2 + P_{\tau}(z)=0.$ The map $\sig_{\tau} \to \cc, (u,z) 
\to z$ presents $\sig_{\tau}$ as a double branched cover 
of the plane, with branch points at the $2m$ components of $\tau.$

The involution (\ref{invo}) induces one on the Hilbert 
scheme 
$\hi^m(S_{\tau}).$ We denote it by $\sigma.$ The fixed point set of 
$\sigma$ parametrizes the closed 0-dimensional 
subschemes in $S_{\tau}$ which are invariant under the involution. 
In particular Fix$(\sigma)$ contains $\hi^m(\sig_{\tau}),$ which by Fact~\ref{fact1} 
is the same as the symmetric product $\sym^m(\sig_{\tau}).$

Theorem~\ref{hilb} presents $\yy_{m, \tau}$ as an open subset of 
$\hi^m(S_{\tau})$ via the morphism (\ref{jaduv}). The involution $\sigma$ 
maps $\yy_{m, \tau}$ to itself. In terms of the polynomials $A(t), 
U(t), V(t)$ 
appearing in (\ref{jaduv}) we have
$$ \sigma: (A(t), U(t), V(t)) \to (A(t), U(t), - V(t)).$$

Going back to (\ref{bcuv}) and (\ref{eq}) we can infer the effect of $\sigma$ on the
polynomials $A(t),$ $B(t),$ $C(t), D(t)$ which record the coefficients of a matrix in the
slice $\ss_m:$ $$ (A(t), B(t), C(t), D(t)) \ \to \ (A(t), C(t), B(t), D(t)).$$ This is
exactly the involution considered by Seidel and Smith \cite{S}.

The following result appears in \cite{S}:
\begin {proposition}
\label {fset}
The fixed point set of $\sigma|_{\yy_{m, \tau}}$ is the complement in 
$\sym^m(\sig_{\tau})$ of the anti-diagonal
$$ \nabla = \{(u_k, z_k), \ k=1, \dots, m:\ u_k^2 + P_{\tau}(z_k) =0, \  
(u_i, z_i) = (-u_j, z_j) \text{ for some } i\neq j \}.$$
\end {proposition}

\noindent {\it Proof.} Looking at equation~(\ref{eq}), the fixed point set
is given by $B(t) = C(t),$ so its points can be described as a triad of
polynomials $A(t), B(t), D(t)$ satisfying \begin {equation} \label {triad}
A(t)D(t) - B(t)^2 = P_{\tau}(t), \end {equation} with $B(t)$ of degree
$m-1,$ $A(t)$ and $D(t)$ monic of degree $m$ and such that the
coefficients of $t^{m-1}$ in $A(t)$ and $D(t)$ sum up to zero.

Recall that the embedding $j: \yy_{m, \tau} \hookrightarrow \hi^m(S_{\tau})$ from (\ref{jaduv}) is
given by taking the roots $z_k$ of $A(t)$ and then setting $u_k = B(z_k)$ for all $k=1, \dots, m.$
It follows that the image of Fix$(\sigma|_{\yy_{m, \tau}})$ under this embedding lies in
$\sym^m(\sig_{\tau}).$ To find the image, pick a point $x \in \sym^m(\sig_{\tau})$ given by
coordinates $(u_k, z_k), k=1,\dots,m.$ This lies in the image of $j$ if and only if we can find a
polynomial $B(t)$ of degree $\leq m-1$ such that $B(z_k) = u_k$ for all $k.$ The necessary and
sufficient condition for the existence of $B$ is that $u_i$ must equal $u_j$ whenever $z_i = z_j.$
This is the same as saying that $x$ is not on the anti-diagonal. $\hfill \fin$

\medskip

Let us turn our attention to the Lagrangians $\kk$ and $\kk'$ from 
Theorem~\ref{lagr}. Note that the Lagrangian 2-spheres 
$\Sigma_{\alpha_k},$ 
are preserved by (\ref{invo}) and the resulting fixed point sets are 
the simple closed curves
\begin {equation}
\label {alfas}
\hat \alpha_k = \{(u, z) \in \cc^2 : z = \alpha_k(t) \text{ for 
some } t \in [0,1]; \ u= \pm \sqrt{-P_{\tau}(z)} \} \subset \sig_{\tau}.
\end {equation}

The same holds true for the beta curves and gives a set of other $m$ 
simple closed curves $\hat \beta_k, \ k=1, \dots, m$ on $\sig_{\tau}.$ 

Consequently, we have:
\begin {proposition}
\label {tori}
The Lagrangians $\kk, \kk' \subset \yy_{m, \tau} \subset 
\hi^m(S_{\tau})$ are mapped into themselves by the involution $\sigma.$ 
The fixed point sets of $\sigma$ restricted to $\kk$ and $\kk'$ are the 
tori
$$ \tt_{\hat \alpha} = \hat \alpha_1 \times \hat \alpha_2 \times \cdots
\times \hat \alpha_m; \ \tt_{\hat \beta} = \hat \beta_1 \times \hat 
\beta_2 \times \cdots \times \hat \beta_m \subset \sym^m(\sig_{\tau}) - 
\nabla,$$ respectively.
\end {proposition}

\subsection {Heegaard Floer homology.} 
\label {sec:osz}
Heegaard Floer homology is a 
powerful tool in low-dimensional topology introduced by Ozsv\'ath and 
Szab\'o in \cite{OS1}, \cite{OS2}. We will be interested in only one 
aspect of their theory, namely the invariant $\widehat{HF}$ of 
3-manifolds.
We will work with cohomology so that we are consistent with our previous 
conventions.

Let $M$ be a closed, oriented 3-manifold. Then $M$ can be described by a
{\it Heegaard diagram}, i.e. a triple $(\Sigma, \hat{\aalpha},
\hat{\bbeta})$ with $\Sigma=\Sigma_g$ a closed oriented surface 
of genus
$g$ and $\hat {\aalpha}= (\hat \alpha_1, \dots, \hat \alpha_g),$
$\hat {\bbeta}= (\hat \beta_1, \dots, \hat \beta_g)$ two collections
of $g$ simple closed curves on $\Sigma,$ such that the $g$ curves in
each collection are linearly independent in $H_1(\Sigma; \zz)$ and
disjoint from the other curves in the same collection.

Given a Heegaard diagram, we can reconstruct $M$ with the help of a
Heegaard splitting, i.e. a decomposition $M=H \cup_{\Sigma} H'$ into two
handlebodies with oriented boundary $\partial H = \Sigma = -\partial H'.$
The handlebody $H$ is obtained from $\Sigma$ by first attaching $g$
two-handles along the $\hat \alpha$ curves, and then attaching one
three-handle (this last step can be done in an essentially unique way).  
Similarly, $H'$ is constructed from $\Sigma$ by first attaching $g$
two-handles along the $\hat \beta$ curves, and then attaching one
three-handle.

Consider the tori $$ \tt_{\hat {\aalpha}} = \hat \alpha_1 \times
\hat \alpha_2 \times \cdots \times \hat \alpha_g; \ \tt_{\hat
{\bbeta}} = \hat \beta_1 \times \hat \beta_2 \times \cdots \times \hat
\beta_g \subset \sym^g(\Sigma)$$ similar to the ones appearing in
Proposition~\ref{tori}. If we pick $w \in \Sigma$ a basepoint disjoint
from the $\hat \alpha$ and $\hat \beta$ curves, we can also think of $
\tt_{\hat {\aalpha}}$ and $\tt_{\hat{\bbeta}}$ as living in 
$\sym^g(\Sigma - w).$ We call $(\Sigma, \hat{ \aalpha},
\hat{\bbeta}, w)$ a {\it pointed Heegaard diagram.}

Of course, there are many Heegaard diagrams that give rise to the same
three-manifold. Ozsv\'ath and Szab\'o have defined an abelian group 
$\widehat{HF}(M)$ by applying a version of Lagrangian Floer cohomology to 
the submanifolds $\tt_{\hat {\aalpha}}$ and $\tt_{\hat{\bbeta}}$
in $\sym^g(\Sigma - w).$ Then they proved that $\widehat{HF}(M)$ is a 
well-defined invariant of the three-manifold $M,$ in the sense that it 
does not depend on the Heegaard diagram chosen to represent $M.$

Let us outline the aspects in the construction of $\widehat{HF}$ which are
of interest to us. The first observation is that $\tt_{\hat
{\aalpha}}$ and $\tt_{\hat{\bbeta}}$ are actually not quite
Lagrangians. If we have a K\"ahler form $\eta$ on $\Sigma,$ then
$\eta^{\times m}$ is a K\"ahler form on $\Sigma^m$ invariant under the
action of the symmetric group.  Hence it descends to a K\"ahler form
$\omega_0$ on $\Sigma^m - \Delta,$ and the tori $\tt_{\hat{
\aalpha}}$ and $\tt_{\hat{\bbeta}}$ are Lagrangian for $\omega_0.$

However, in general $\omega_0$ cannot be extended over the diagonal
$\Delta,$ so a modification of the usual construction is needed. This is
done by using a class of almost complex structures $J_s$ on
$\sym^g(\Sigma),$ which are chosen to tame $\omega_0$ in a neighborhood of
the tori. It follows that $\tha$ and
$\thb$ are totally real submanifolds of
$(\sym^g(\Sigma - w), J_s).$ Under a certain admissibility condition for
the pointed Heegaard diagram (explained below), one can still
count pseudo-holomorphic disks in $(\sym^g(\Sigma - w), J_s)$ and get a
Floer cohomology group 
$$HF(\tha,\thb) = \widehat{HF}(M).$$ 

The generators of the cochain complex $CF^*(\tha, \thb)$ are still given by the
intersection points between $\tha$ and $\thb,$ provided that their intersection is made
transverse after a small isotopy. The Floer cohomology groups $\widehat{HF}(M) =
\widehat{HF}^*(M)$ admit a $\zz/2$ grading, given by the sign of these intersection points.

The admissibility condition for pointed Heegaard diagrams needed in the
construction of $\widehat{HF}(M)$ is called {\it weak admissibility} in
\cite[Definition 4.10]{OS1}. The definition there is given in terms of a
decomposition of $\widehat{HF}(M)$ according to $\text{spin}^c$ structures
on $M.$ It is then proved that for every $\text{spin}^c$ structure $\spc$
on $M$ there is a pointed Heegaard diagram which is weakly admissible with
respect to $\spc$ and that diagram can be used to define a group
$\widehat{HF}(M, \spc).$ Then we take the direct sum of $\widehat{HF}(M,
\spc)$ over all $\spc$ to obtain $\widehat{HF}(M).$

For our purposes, it is only important that a pointed Heegaard 
diagram is weakly admissible for all $\spc$ if it satisfies the 
requirements in the following definition:

\begin {definition} 
\label {admy}
Let $(\Sigma, \hat {\aalpha},\hat{\bbeta}, w)$ be a
pointed Heegaard diagram. We denote by $D_1, \dots, D_s$ the closures of
the components of $\Sigma - \hat \alpha_1 - \dots - \hat\alpha_g - \hat
\beta_1 - \dots - \hat \beta_g,$ with the convention that $D_s$ is the
component containing $w.$ We say that a two-chain $\pp = \sum_{i=1}^{s-1} 
n_i D_i, n_i \in \zz$ is a {\bf periodic domain} if its boundary 
is a sum of $\hat \alpha$ and $\hat \beta$ curves.

The pointed Heegaard diagram $(\Sigma, \hat {\aalpha},\hat{\bbeta}, w)$ 
is called {\bf admissible} if every nontrivial periodic domain 
admits both negative and positive coefficients among the $n_i$'s. 
\end {definition}

\subsection {Double branched covers.} \label {sec:dbc} 
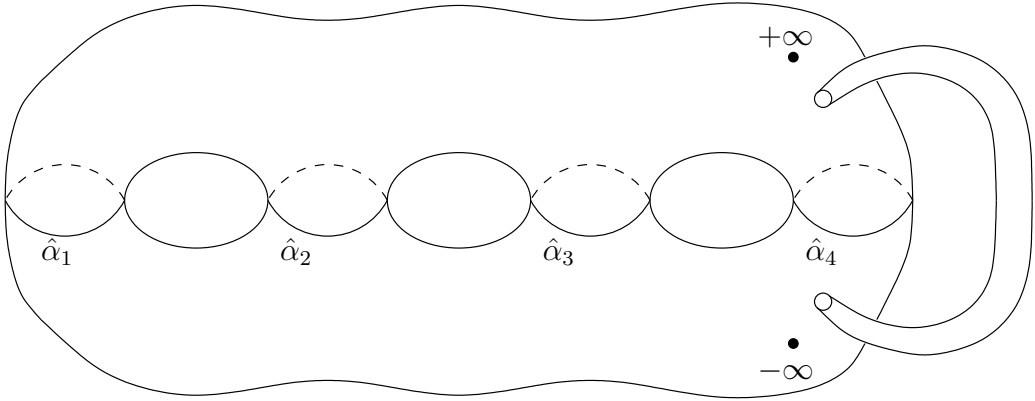
\begin {figure}
\begin {center}
\input {hbody.pstex_t}
\end {center}
\caption {The surface $\sigc_{\tau}$ for $m=4.$}
\label {fig:hbody}
\end {figure}

The collection of curves $\hat
\aalpha = (\hat \alpha_1, \dots, \hat \alpha_m)$ from (\ref{alfas})  and
its analogue $\hat \bbeta$ can be used as a set of attaching curves in a
Heegaard diagram. We need to be careful though. Note that $\sig_{\tau}$ is
the complement of two points at infinity $\pm \infty$ in the closed
surface $\Sigma_{m-1},$ the double cover of $S^2$ branched over $2m$
points. A collection of $m$ curves needs to be put on a surface of genus
$m,$ so the cure is to add an extra handle. We do this by removing two 
small disks from neighborhoods of $\pm \infty$ (but such that the disks 
do not contain $\pm\infty$) and joining their boundaries by a handle, as 
shown in Figure~\ref{fig:hbody}. We call the resulting surface 
$\sigc_{\tau}= \Sigma_m.$

\begin {proposition} 
\label {prohf}
\label {dbc} $(\sigc_{\tau}, \hat \aalpha, \hat
\bbeta, +\infty)$ is an admissible Heegaard diagram for the manifold
$\dd(L)$ $\#$ $(S^1 \times S^2),$ where $\dd(L)$ indicates 
the double cover of $S^3$ branched over the link $L$ and $\#$ denotes 
connected sum. \end{proposition}

\begin {figure}
\begin {center}
\input {disks.pstex_t}
\end {center}
\caption {Lifting the $\alpha$ curves into the ball.}
\label {fig:disks}
\end {figure}
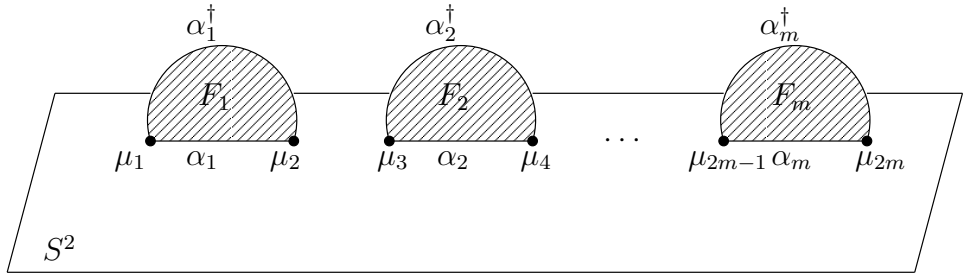

\noindent{\it Proof.} The $\alpha$ curves form a crossingless matching of
$2m$ points $\mu_1, \dots, \mu_{2m}$ in the complex plane $\cc.$ We add a
point at infinity to $\cc$ to obtain $S^2$ and think of that as the
boundary of a ball $B^3.$ We ``lift'' the $\alpha$ curves to form $m$ 
segments $\alpha_1^{\dag}, \dots, 
\alpha_m^{\dag}$ in $B^3$ joining the points $\mu_k$ in pairs just like 
the $\alpha$'s do, but such that their interiors do not intersect the boundary 
$\partial B^3 = S^2.$ Then there are $m$ small disks $F_1, \dots, F_m$ in 
$B^3$ with boundaries $\partial F_k = \alpha_k \cup (-\alpha_k^{\dag}).$ 
This is shown in Figure~\ref{fig:disks}.

We have a similar picture for the beta curves, with corresponding lifts
$\beta_1^{\dag}, \dots, \beta_m^{\dag}.$ We can form $S^3$ by joining the
two balls $B^3$ along their common boundaries. One half contains the
$\alpha^{\dag}$ curves, and the other half the $\beta^{\dag}$ curves.
The union of the $\alpha^{\dag}$ and $\beta^{\dag}$ is exactly the link $L 
\subset S^3.$

If we take the double cover of $B^3$ branched over the curves 
$\alpha_k^{\dag}$
we obtain a handlebody $H_0$ of genus $m-1$ with $\partial H_0 =
\sig_{\tau}.$ The preimage of $\alpha_k$ is $\hat \alpha_k,$ for $k=1, 
\dots, m.$ The homology classes $[\hat \alpha_k] \in H_1(\partial H_0;  
\zz)$ add up to zero. Also, each curve $\hat \alpha_k$ bounds a disk in
$H_0$ which is the preimage of $F_k$ under the double covering map. Adding
one extra handle to $H_0$ we obtain a genus $m$ handlebody $H.$ The
picture is exactly the one in Figure~\ref{fig:hbody}, and the
$\hat\alpha_k$ curves can serve as attaching circles for $H.$ 

The same construction works for the beta curves on the other half of
$B^3.$ The result is a handlebody $H'$ with attaching circles $\hat 
\beta_k.$ Gluing $H$ and $H'$ together along the boundary we obtain 
$\dd(L)$ with one handle attached, which is $\dd(L) \# (S^1 \times S^2),$  
as desired.

We are left to check admissibility. Let $\pp = \sum_{i=1}^s n_i D_i$ 
be a nontrivial periodic domain as in Definition~\ref{admy}, with $n_s = 
0$ as required. The involution (\ref{invo}) induces an involution $i \to 
\bar i$ on the set $\{1, 2, \dots, s \},$ according to how a domain $D_i$ 
on $\Sigma_m = \sigc_{\tau}$ is taken to another domain $D_{\bar i}.$ 
Because of our choice of the basepoint near infinity we have $\bar s = s.$ 
We claim that 
\begin {equation}
\label {nii}
n_i + n_{\bar i} = 0
\end {equation}
for all $i=1, \dots, s.$ 

Let $D_i$ and $D_j$ be two adjacent components, i.e. so that they have an
edge $E$ in common. The edge $E$ could be either part of a $\hat \alpha$
curve or of a $\hat \beta$ curve.  If we pick orientations on
$\sigc_{\tau}$ and on the $\hat \alpha$ and $\hat \beta$ curves
such that the oriented boundary of $D_i$ has an $E$ term, then $\partial
D_j$ has a term $-E.$ The involution (\ref{invo}) takes $E$ to another
edge $\bar E,$ which lies on the same $\hat \alpha$ (or $\hat \beta$)  
curve, but comes with the opposite orientation. Note that $\partial
D_{\bar i}$ has a term $\bar E$ and $\partial D_{\bar j}$ a term $-\bar
E.$ Since the boundary of $\pp = \sum n_i D_i$ must be a sum of alpha and
beta curves, $E$ and $\bar E$ must appear with opposite signs as part of 
$\partial \pp.$ Therefore, 
$$n_i - n_j =  - n_{\bar i} + n_{\bar j}.$$

It follows that if $i$ satisfies (\ref{nii}), then so does $j.$ Since 
the surface $\sigc_{\tau}$ is connected and we already know that 
$n_s = n_{\bar s} = 0,$ we get recursively that $(\ref{nii})$ is true for 
all $i.$ Since not all $n_i$'s are zero by assumption, there must be at 
least one positive and one negative integer among them. This means that 
$ (\sigc_{\tau}, \hat \aalpha, \hat \bbeta, +\infty)$ is admissible.
$\hfill \fin $
\medskip

A quick corollary of Proposition~\ref{dbc} is a description of a set of 
generators for the cochain complex $CF^*(\tha, \thb) = 
\widehat{CF}^*(\dd(L) \# (S^1 \times S^2)).$ They are the intersection 
points in $$ \hat \gz = \tha \cap \thb.$$ 

In section~\ref{sec:bigg} we denoted by $\bar \ze$ the union of all 
intersection points between the $\alpha$ and the $\beta$ curves in the 
plane. Now let $\hat \ze$ be the union of all
intersection points between the $\hat \alpha$ and the $\hat \beta$ curves 
in the surface $\sig_{\tau},$ or equivalently in $\sigc_{\tau}.$ 
There is a natural map coming from the double cover
\begin {equation}
\label {fine}
 \hat f: \hat \ze \to \bar \ze.
\end {equation}

An element $x \in \bar \ze$ has one preimage $\hat e_x \in \hat \ze$ under
$f$ in case it is one of the puncture points $\mu_k,$ and two preimages
$\hat e_x, \hat e_x'$ otherwise. This situation is very similar to that of
the map $f: \ze \to \bar \ze$ from (\ref{ff}). However, there is no
canonical way of identifying the set $\{e_x, e_x'\}$ with $\{\hat e_x,
\hat e_x'\}.$ Hence $\ze$ and $\hat \ze$ have the same number of elements, 
but they cannot be identified in a natural way. 

The set $\hat \gz$ can be recovered from $\hat \ze$ similarly to how $\gz$ 
and $\bar \gz$ were constructed from $\ze$ and $\bar \ze$ in section 
\ref{sec:bigg}. It follows that if we fix identifications of $\{e_x, 
e_x'\}$ with $\{\hat e_x, \hat e_x'\}$ for all the $x \in \ze$ that are 
not punctures, we get an identification of $\hat \gz$ with the set of Bigelow 
generators $\gz.$ 

\subsection {Grading.} \label {sec:gradf} The elements of $\hat \gz$ form generators for
$\widehat{CF}(\dd(L)\# (S^1 \times S^2)),$ and we would like to understand
their cohomological grading. As explained in section~\ref{sec:osz}, this
grading is well-defined only modulo $2.$ However, it can be improved to a $\zz$ grading in the 
following way. 

First, note that because our identification of $\hat \gz$ with the set $\gz$ of Bigelow 
generators was not canonical, we cannot talk about all the gradings 
defined in sections \ref{sec:bigj} and \ref{sec:flob} on $\gz$ as gradings 
on $\hat \gz.$ Nevertheless, if a grading $F: \gz \to \zz$ is stable in 
the sense of Definition~\ref{defs}, then it comes from a grading $\bar F$ 
on $\bar \gz,$ and by composing with (\ref{fine}) we can think of $F$ as a 
well-defined grading on $\hat \gz.$ This is the case of the gradings $T$ 
and $P-Q$ (or $\tilde P - Q$) as seen in section~\ref{sec:bigg} and in Lemma~\ref{pqs}.

Now, instead of considering the tori $\tha, \thb$ as embedded in $\sym^m(\sigc_{\tau}),$ we view
them as totally real submanifolds of $W = \sym^m(\sig_{\tau}) - \nabla$ as in 
Proposition~\ref{tori}.
Recall from the proof of Proposition~\ref{fset} that $W$ is an affine algebraic variety given by the
equations (\ref{triad}) in the coefficients of $A, B$ and $D.$ Let us equip $W = $
Fix$(\sigma|_{\yy_{m, \tau}})$ with the restriction of the K\"ahler form $\tilde \Omega$ on $\yy_{m,
\tau}.$ We apply the formalism in section~\ref{sec:floer} (the totally real case) to $Y = W, \ttt =
\tha$ and $\ttt' = \thb$ .

\begin {proposition}
\label {bach}
There exists a complex volume form $\Theta$ on $W$ so that we can endow $\tha$ and $\thb$ with 
gradings in the sense of Definition~\ref{deaf}. The resulting absolute Maslov grading 
on the elements of $\hat \gz = \tha \cap \thb$ is $\tilde P - Q + T.$
\end {proposition}

For the proof we need the following:
\begin {lemma}
\label {vander}
Let $w_1, \dots, w_n$ be $n$ formal variables, and set $Q(t) = \prod 
(t-w_i) = t^n - s_1 t^{n-1} + \dots + (-1)^n s_n,$ so that $s_1, s_2, \dots, s_n$ are the 
symmetric polynomials in $w_1, \dots, w_n.$
Then:
$$ ds_1 \wedge \dots \wedge ds_n = \prod_{i<j} (w_i - w_j) \cdot dw_1 \wedge \dots 
\wedge dw_n.$$ 
\end {lemma}

\noindent{\it Proof.} Differentiating the relation $Q(w_i) = 0$ we get $\sum_{j=1}^n 
(-1)^j w_i^{n-j} ds_j = -Q'(w_i)dw_i$ for all $i=1, \dots, n.$ The formula for the Vandermonde 
determinant gives
\begin {equation}
\label {color} 
(-1)^{n(n+1)/2}\prod_{i<j} (w_i - w_j)\cdot ds_1 \wedge \dots \wedge ds_n = (-1)^n \prod_i Q'(w_i) 
\cdot dw_1 
\wedge \dots \wedge dw_n.
\end{equation}
Observe that $Q'(w_i) = \prod_{j \neq i} (w_j - w_i).$ The result now follows after dividing both 
sides of (\ref{color}) by $(-1)^{n(n+1)/2}\prod_{i<j} (w_i - w_j). \hfill \fin$

\bigskip

\noindent{\it Proof of Proposition~\ref{bach}.} We start with the construction of the complex volume 
form $\Theta$ on $W.$ Recall that $\sig_{\tau}$ is given by the equation $u^2 + P_{\tau}(z) =0$ and 
this gives a set of coordinates $u_j, z_j \ (j=1, \dots, m)$ on its symmetric product. Set
\begin {equation}
\label {thety}
 \Theta = \prod_{1\leq i < j \leq m} (z_i - z_j) \cdot \prod_{j=1}^m 
\frac{dz_j}{u_j}. \end {equation}

Note that $dz_j/ u_j$ is a well-defined form on $\sig_{\tau},$ and hence
the expression (\ref{thety}) is well-defined as a form on the Cartesian
product $(\sig_{\tau})^m.$ It is invariant under the action of the
symmetric group on $m$ elements, because if we switch $z_i$ and $z_j$ this
produces a minus sign in both of the products appearing in (\ref{thety}).
Thus it descends to a complex $m$-form on $\sym^m(\sig_{\tau})$ which
clearly has no zeros or poles outside $\Delta \cup \nabla,$ and therefore
gives a good volume form there. 

We claim that $\Theta$ extends to a volume form on all of $W = \sym^m(\sig_{\tau}) - \nabla.$ We
need to check this for points $x \in \Delta - \nabla.$ A point of this type is an $m$-tuple of pairs
$(u_j, z_j)$ with $u_j^2 = -P_{\tau}(z_j) \neq 0$ for all $j = 1, \dots, m,$ and $(u_i, z_i) = (u_j,
z_j)$ for some $i \neq j.$ The point $x$ determines a partition of $\{1, \dots, m\}$ into blocks of
the form $\{j_1, \dots, j_n\}$ such that $z_{j_1} = \dots = z_{j_n}$ and the $z_j$'s are different
for $j$ in different blocks. Putting together the symmetric functions $s_1, \dots, s_n$ in $z_{j_1}, 
\dots, z_{j_n}$ for each block we form a set of local coordinates on $W$ around $x.$ Using 
Lemma~\ref{vander} and the fact that $u_j \neq 0$ we get that $\Theta$ is a nonzero multiple of the 
product of all $ds_1 \wedge \dots \wedge ds_n,$ taken over all blocks and in a neighborhood of $x.$ 
Thus $\Theta$ extends to a well-defined complex volume form near any $x \in W.$

Now we can do a computation similar to that in section~\ref{sec:proj}. A 
point $x \in \thb$ has coordinates $(u_j, z_j),$ with $z_j = \beta_j(t_j)$ 
for some $t_j \in [0,1]$ and $u_j = \pm \sqrt{-P_{\tau}(\beta_j(t_j))}.$
The resulting square phase map (\ref{sqp}) on $\thb$ is
$$ \theta_{\thb}: \thb \to \cc^*/\rr_+, \ \theta_{\thb}(x) =  \prod_{1\leq i 
< j 
\leq m} 
\bigl( \beta_i(t_i) - \beta_j(t_j)\bigr)^2 \cdot \prod_{j=1}^m
\frac{\beta_j'(t_j)^2}{-P_{\tau}(\beta_j(t_j))}.$$

Since $H^1(\thb)$ is nontrivial, we did not know {\it a priori} that the 
square phase map lifted to a real-valued function. Now we know that this 
is actually the case, because it factors through the projection to the
contractible space $\beta_1 \times \dots \times \beta_m.$

There is a similar square phase map for $\tha,$ and that turns out to be
null-homotopic because the $\alpha$ curves are horizontal. Thus, we can endow both $\tha$ and $\thb$ 
with gradings, i.e. with lifts of the square phase maps to $\tilde \theta_{\tha}: \tha \to \rr, 
\tilde \theta_{\thb}: \thb \to \rr.$ We choose them just like in the Seidel-Smith picture 
(section~\ref{sec:floer}), i.e. we obtain the grading on $\thb$ from that on $\tha$ by following 
continuously the family of crossingless matchings in the plane determined by the braid $b \times 
1^m.$

As explained in section~\ref{sec:floer}, this induces an absolute Maslov grading $\tilde R(x) \in
\zz$ on the points $x \in \tha \cap \thb.$ To make it more explicit, recall that in the standard
picture from section~\ref{sec:braid} all alpha curves are subsets of the real line. We can assume
that they always intersect the beta curves at $90$ degree angles. Also, at the endpoints $\mu_j$ the
beta curves need to be parametrized so that their derivatives vanish. Then $\theta_{\thb}(x)= -1 \in
S^1$ for every $x \in \tha \cap \thb.$ Hence $\tilde \theta_{\thb}(x)$ is an odd integer.  If the
lift $\tilde \theta_{\thb}$ is chosen so that the special point $(\mu_2, \dots, \mu_{2m})$ is mapped
to $1,$ then the Maslov index $\tilde R(x)$ is equal to $k \in \zz$ at the points where $\tilde
\theta_{\thb}(x)=2k+1.$

Just like in section~\ref{sec:proj}, we get that the resulting Maslov grading can be computed from a 
flattened braid diagram. We use the standard additivity properties of the index. The factor $\prod
\bigl( \beta_i(t_i) - \beta_j(t_j)\bigr)^2$ counts the twisting of the points around each other, 
which is the $T$ grading. 

The second factor produces a grading which is clearly both additive and stable. We decompose it 
into contributions coming from 
\begin {equation}
\label {fbb}
f_j(t_j)=\frac{\beta_j'(t_j)^2}{-P_{\tau}(\beta_j(t_j))}
\end {equation}
for each $j.$ The value of the respective contribution $C^*(p)$ at a point $p=\beta_j(t_j) \in 
\alpha_i \cap \beta_j$ can be computed as follows. At the endpoint $\mu_{2j}$ of $\beta_j$ we have
$C^*(\mu_{2j}) =0.$ We then follow the curve $\beta_j$ in reverse until we hit $p,$ and look at the 
lift $\tilde f_j$ of $f_j$ to a real-valued function. Then $C^*(p)$ describes on which sheet we end 
up, i.e. it equals half the difference $f_j(t_j) - f_j(1).$ It is a bit hard 
to see the result because at the endpoints of $\beta$ both the numerator and the denominator in
(\ref{fbb}) are zero. However, the Maslov index is invariant under deformation, so we can replace
$\beta_j$ with one half $H_j$ of the corresponding figure-eight $E_j.$ (It does not matter which of
the two halves we choose.) $H_j$ is an arc ending at $\nu_j$ near $\mu_{2j}$ and starting at the 
intersection point near the other endpoint of $\beta_j.$ Now $C^*$ decomposes into a contribution 
from the numerator $H_j'(t_j)^2,$ which is the projective grading $\tilde P^*,$ and one from the 
denominator
$(-1)^m \prod P_{\tau}(H_j(t_j)),$ which records the twisting around the puncture points and hence
gives the $Q^*$ grading. Summing up all these contributions we get $\tilde R = \tilde P- Q +
T,$ as desired. $\hfill \fin$

\medskip

\noindent{\it Proof of Theorem~\ref{newth}.}
The identification claimed in the first part of the theorem was explained at the end of 
section~\ref{sec:dbc}. The second part of the theorem is basically Proposition~\ref{bach}. 

\medskip

\noindent{\it Proof of Corollary~\ref{coro}.} Recall that in the introduction we denoted 
$$R = \tilde R - (m+w)/2  = P + (J/2).$$

The expression for the Jones polynomial in terms of $P$ and $R$ follows readily from (\ref{ujones}),
Bigelow's formula (\ref{ujon}), and Remark~\ref{sign}. $\hfill \fin$

\subsection {Reduced theories.} Seidel and Smith suggested the 
involution $\tau$ as a key to a geometric interpretation of the spectral 
sequence in \cite{OS3}, which relates Khovanov homology to 
$\widehat{HF}(\dd(L) \# (S^1 \times S^2)).$ 

The work of Ozsv\'ath and Szab\'o in \cite{OS3} involves the reduced
Khovanov homology $\widetilde{Kh}(L).$ This was defined in \cite{Kh3}. It
is a homology theory combinatorially defined starting with a plane diagram
for the link $L,$ just like $Kh(L).$ It is an invariant of $L$ only when
considered with $\zz/2$ coefficients. With $\zz$ coefficients
$\widetilde{Kh}(L)$ depends on a distinguished component of $L;$ in
particular it still gives a well-defined invariant of knots, for example.
The Euler characteristic of $\widetilde{Kh}$ is the usual Jones polynomial 
$V_L(t).$

The spectral sequence in \cite{OS3} is defined with $\zz/2$ 
coefficients only, has as $E^2$ term the reduced Khovanov 
homology of the mirror of $L,$ and converges to $E^{\infty} = 
\widehat{HF}(\dd(L); \zz/2).$ (Taking the mirror only has the effect 
of changing the signs in the bidegree of $Kh$ and $\widetilde{Kh}$ in a 
well-understood manner \cite{Kh}, so we will not worry about it.) If we 
add an unlinked unknot $O$ to $L,$ we have $\widetilde{Kh}(L \cup O) = 
Kh(L)$ and, according to \cite[Proposition 6.4]{OS2}: 
$$\widehat{HF}^*(\dd(L \cup O)) = \widehat{HF}^*(\dd(L) \# (S^1 \times 
S^2)) = \widehat{HF}^*(\dd(L)) \otimes H^*(S^1).$$

Thus we have a spectral sequence from $Kh$ of the mirror of $L$ to 
$\widehat{HF}(\dd(L) \# (S^1 \times S^2)),$ which is what we mentioned 
before.

It is worthwhile seeing how many of the constructions in this paper admit
a ``reduced'' version as well. First of all, if in Proposition~\ref{prohf}
we eliminate a pair of corresponding $\hat \alpha$ and $\hat \beta$
curves, say $\hat \alpha_m$ and $\hat \beta_m,$ and use 
$\Sigma_{m-1} = \sig_{\tau} \cup \{\pm \infty \}$
without the handle instead of $\sigc_{\tau},$ then the same proof applies
and we get an admissible Heegaard diagram for $\dd(L).$ Thus we can
describe a set of generators for $\widehat{HF}(\dd(L)).$

Jacob Rasmussen pointed out to the author that there is a reduced variant
of Bigelow's picture. We need a mild condition on the flattened braid
diagram, namely that there is a path from $\mu_{2m}$ to infinity that
does not intersect any of the $\alpha$ and $\beta$ curves in any other
point. Under this condition one can eliminate $\alpha_m$ and $\beta_m$
from the picture, get a set of reduced Bigelow generators, compute their
$Q, T$ and $J$ gradings as before, and obtain the Jones polynomial $V_L$
instead of the unnormalized one $J_L.$ (For the $Q$ grading, one still 
needs to consider the total winding number around all the $\mu_k$'s, 
including $\mu_{2m-1}$ and $\mu_{2m}$.)

Finally, on the Seidel-Smith side a reduced version is yet to be
developed. We expect that the construction from \cite{SS} works also for
the slice $\ss_{m-1}$ instead of $\ss_m,$ and that one can define a
similar Lagrangian Floer cohomology theory on $\yy_{m-1, \tau}$ for $\tau
\in \conf^{2m}_0(\cc).$ The result is likely to be a well-defined
invariant for knots with $\zz$ coefficients, and for links with $\zz/2$
coefficients, as is the case for reduced Khovanov homology. The reduced 
Bigelow generators should then give a set of generators for the reduced 
Seidel-Smith cochain complex.

\end{document}

%% file: circ.pstex_t
\begin{picture}(0,0)%
\includegraphics{circ.pstex}%
\end{picture}%
\setlength{\unitlength}{3947sp}%
\begingroup\makeatletter\ifx\SetFigFont\undefined%
\gdef\SetFigFont#1#2#3#4#5{%
  \reset@font\fontsize{#1}{#2pt}%
  \fontfamily{#3}\fontseries{#4}\fontshape{#5}%
  \selectfont}%
\fi\endgroup%
\begin{picture}(2775,3168)(826,-3820)
\put(3601,-2686){\makebox(0,0)[lb]{\smash{\SetFigFont{12}{14.4}{\rmdefault}{\mddefault}{\updefault}{\color[rgb]{0,0,0}$1^m$}%
}}}
\put(826,-2761){\makebox(0,0)[lb]{\smash{\SetFigFont{12}{14.4}{\rmdefault}{\mddefault}{\updefault}{\color[rgb]{0,0,0}b}%
}}}
\end{picture}

%% file: match.pstex_t
\begin{picture}(0,0)%
\includegraphics{match.pstex}%
\end{picture}%
\setlength{\unitlength}{3947sp}%
\begingroup\makeatletter\ifx\SetFigFont\undefined%
\gdef\SetFigFont#1#2#3#4#5{%
  \reset@font\fontsize{#1}{#2pt}%
  \fontfamily{#3}\fontseries{#4}\fontshape{#5}%
  \selectfont}%
\fi\endgroup%
\begin{picture}(4763,439)(376,-869)
\put(376,-811){\makebox(0,0)[lb]{\smash{\SetFigFont{12}{14.4}{\rmdefault}{\mddefault}{\updefault}{\color[rgb]{0,0,0}$\mu_1$}%
}}}
\put(2926,-811){\makebox(0,0)[lb]{\smash{\SetFigFont{12}{14.4}{\rmdefault}{\mddefault}{\updefault}{\color[rgb]{0,0,0}$\mu_4$}%
}}}
\put(3976,-811){\makebox(0,0)[lb]{\smash{\SetFigFont{12}{14.4}{\rmdefault}{\mddefault}{\updefault}{\color[rgb]{0,0,0}$\mu_{2m-1}$}%
}}}
\put(5026,-811){\makebox(0,0)[lb]{\smash{\SetFigFont{12}{14.4}{\rmdefault}{\mddefault}{\updefault}{\color[rgb]{0,0,0}$\mu_{2m}$}%
}}}
\put(4501,-586){\makebox(0,0)[lb]{\smash{\SetFigFont{12}{14.4}{\rmdefault}{\mddefault}{\updefault}{\color[rgb]{0,0,0}$\alpha_{m}$}%
}}}
\put(2401,-586){\makebox(0,0)[lb]{\smash{\SetFigFont{12}{14.4}{\rmdefault}{\mddefault}{\updefault}{\color[rgb]{0,0,0}$\alpha_2$}%
}}}
\put(1351,-811){\makebox(0,0)[lb]{\smash{\SetFigFont{12}{14.4}{\rmdefault}{\mddefault}{\updefault}{\color[rgb]{0,0,0}$\mu_2$}%
}}}
\put(2026,-811){\makebox(0,0)[lb]{\smash{\SetFigFont{12}{14.4}{\rmdefault}{\mddefault}{\updefault}{\color[rgb]{0,0,0}$\mu_3$}%
}}}
\put(901,-586){\makebox(0,0)[lb]{\smash{\SetFigFont{12}{14.4}{\rmdefault}{\mddefault}{\updefault}{\color[rgb]{0,0,0}$\alpha_1$}%
}}}
\end{picture}

%% file: tref.pstex_t
\begin{picture}(0,0)%
\includegraphics{tref.pstex}%
\end{picture}%
\setlength{\unitlength}{3947sp}%
\begingroup\makeatletter\ifx\SetFigFont\undefined%
\gdef\SetFigFont#1#2#3#4#5{%
  \reset@font\fontsize{#1}{#2pt}%
  \fontfamily{#3}\fontseries{#4}\fontshape{#5}%
  \selectfont}%
\fi\endgroup%
\begin{picture}(4307,3249)(419,-3595)
\put(2626,-2536){\makebox(0,0)[lb]{\smash{\SetFigFont{12}{14.4}{\familydefault}{\mddefault}{\updefault}{\color[rgb]{0,0,0}$ \mu_2$}%
}}}
\put(4726,-2536){\makebox(0,0)[lb]{\smash{\SetFigFont{12}{14.4}{\familydefault}{\mddefault}{\updefault}{\color[rgb]{0,0,0}$\mu_4$}%
}}}
\put(1126,-1186){\makebox(0,0)[lb]{\smash{\SetFigFont{12}{14.4}{\familydefault}{\mddefault}{\updefault}{\color[rgb]{0,0,0}$\beta_1$}%
}}}
\put(1126,-661){\makebox(0,0)[lb]{\smash{\SetFigFont{12}{14.4}{\familydefault}{\mddefault}{\updefault}{\color[rgb]{0,0,0}$\beta_2$}%
}}}
\put(1201,-2536){\makebox(0,0)[lb]{\smash{\SetFigFont{12}{14.4}{\rmdefault}{\mddefault}{\updefault}{\color[rgb]{0,0,0}$\mu_1$}%
}}}
\put(3376,-2536){\makebox(0,0)[lb]{\smash{\SetFigFont{12}{14.4}{\familydefault}{\mddefault}{\updefault}{\color[rgb]{0,0,0}$ \mu_3$}%
}}}
\put(4051,-2386){\makebox(0,0)[lb]{\smash{\SetFigFont{12}{14.4}{\rmdefault}{\mddefault}{\updefault}{\color[rgb]{0,0,0}$\alpha_2$}%
}}}
\put(1951,-2386){\makebox(0,0)[lb]{\smash{\SetFigFont{12}{14.4}{\familydefault}{\mddefault}{\updefault}{\color[rgb]{0,0,0}$\alpha_1$}%
}}}
\end{picture}

%% file: unknot.pstex_t
\begin{picture}(0,0)%
\includegraphics{unknot.pstex}%
\end{picture}%
\setlength{\unitlength}{3947sp}%
\begingroup\makeatletter\ifx\SetFigFont\undefined%
\gdef\SetFigFont#1#2#3#4#5{%
  \reset@font\fontsize{#1}{#2pt}%
  \fontfamily{#3}\fontseries{#4}\fontshape{#5}%
  \selectfont}%
\fi\endgroup%
\begin{picture}(4088,689)(376,-891)
\put(3301,-811){\makebox(0,0)[lb]{\smash{\SetFigFont{12}{14.4}{\rmdefault}{\mddefault}{\updefault}{\color[rgb]{0,0,0}$\mu_1$}%
}}}
\put(4276,-811){\makebox(0,0)[lb]{\smash{\SetFigFont{12}{14.4}{\rmdefault}{\mddefault}{\updefault}{\color[rgb]{0,0,0}$\mu_2$}%
}}}
\put(376,-811){\makebox(0,0)[lb]{\smash{\SetFigFont{12}{14.4}{\rmdefault}{\mddefault}{\updefault}{\color[rgb]{0,0,0}$\mu_1$}%
}}}
\put(1351,-811){\makebox(0,0)[lb]{\smash{\SetFigFont{12}{14.4}{\rmdefault}{\mddefault}{\updefault}{\color[rgb]{0,0,0}$\mu_2$}%
}}}
\end{picture}

%% file: skein.pstex_t
\begin{picture}(0,0)%
\includegraphics{skein.pstex}%
\end{picture}%
\setlength{\unitlength}{3947sp}%
\begingroup\makeatletter\ifx\SetFigFont\undefined%
\gdef\SetFigFont#1#2#3#4#5{%
  \reset@font\fontsize{#1}{#2pt}%
  \fontfamily{#3}\fontseries{#4}\fontshape{#5}%
  \selectfont}%
\fi\endgroup%
\begin{picture}(3920,1195)(889,-1844)
\put(1276,-1786){\makebox(0,0)[lb]{\smash{\SetFigFont{12}{14.4}{\rmdefault}{\mddefault}{\updefault}{\color[rgb]{0,0,0}$L_+$}%
}}}
\put(2776,-1786){\makebox(0,0)[lb]{\smash{\SetFigFont{12}{14.4}{\rmdefault}{\mddefault}{\updefault}{\color[rgb]{0,0,0}$L_-$}%
}}}
\put(4276,-1786){\makebox(0,0)[lb]{\smash{\SetFigFont{12}{14.4}{\rmdefault}{\mddefault}{\updefault}{\color[rgb]{0,0,0}$L_0$}%
}}}
\end{picture}

%% file: fig8.pstex_t
\begin{picture}(0,0)%
\includegraphics{fig8.pstex}%
\end{picture}%
\setlength{\unitlength}{3947sp}%
\begingroup\makeatletter\ifx\SetFigFont\undefined%
\gdef\SetFigFont#1#2#3#4#5{%
  \reset@font\fontsize{#1}{#2pt}%
  \fontfamily{#3}\fontseries{#4}\fontshape{#5}%
  \selectfont}%
\fi\endgroup%
\begin{picture}(7146,4329)(687,-5206)
\put(5626,-1711){\makebox(0,0)[lb]{\smash{\SetFigFont{10}{12.0}{\rmdefault}{\mddefault}{\updefault}{\color[rgb]{0,0,0}$E_1$}%
}}}
\put(6076,-1186){\makebox(0,0)[lb]{\smash{\SetFigFont{10}{12.0}{\rmdefault}{\mddefault}{\updefault}{\color[rgb]{0,0,0}$E_2$}%
}}}
\put(2551,-3886){\makebox(0,0)[lb]{\smash{\SetFigFont{9}{10.8}{\rmdefault}{\mddefault}{\updefault}{\color[rgb]{0,0,0}$\uh_1'$}%
}}}
\put(2926,-3661){\makebox(0,0)[lb]{\smash{\SetFigFont{9}{10.8}{\rmdefault}{\mddefault}{\updefault}{\color[rgb]{0,0,0}$\uh_1$}%
}}}
\put(3676,-3811){\makebox(0,0)[lb]{\smash{\SetFigFont{9}{10.8}{\rmdefault}{\mddefault}{\updefault}{\color[rgb]{0,0,0}$\uh_2$}%
}}}
\put(5926,-3661){\makebox(0,0)[lb]{\smash{\SetFigFont{9}{10.8}{\rmdefault}{\mddefault}{\updefault}{\color[rgb]{0,0,0}$\yh_1$}%
}}}
\put(6226,-3811){\makebox(0,0)[lb]{\smash{\SetFigFont{9}{10.8}{\rmdefault}{\mddefault}{\updefault}{\color[rgb]{0,0,0}$\vh_1$}%
}}}
\put(6751,-3811){\makebox(0,0)[lb]{\smash{\SetFigFont{9}{10.8}{\rmdefault}{\mddefault}{\updefault}{\color[rgb]{0,0,0}$\yh_2$}%
}}}
\put(6526,-3661){\makebox(0,0)[lb]{\smash{\SetFigFont{9}{10.8}{\rmdefault}{\mddefault}{\updefault}{\color[rgb]{0,0,0}$\vh_1'$}%
}}}
\put(7126,-3661){\makebox(0,0)[lb]{\smash{\SetFigFont{9}{10.8}{\rmdefault}{\mddefault}{\updefault}{\color[rgb]{0,0,0}$\yh_2'$}%
}}}
\put(7426,-3811){\makebox(0,0)[lb]{\smash{\SetFigFont{9}{10.8}{\rmdefault}{\mddefault}{\updefault}{\color[rgb]{0,0,0}$\vh_2$}%
}}}
\put(2251,-3661){\makebox(0,0)[lb]{\smash{\SetFigFont{9}{10.8}{\rmdefault}{\mddefault}{\updefault}{\color[rgb]{0,0,0}$\xh_1$}%
}}}
\put(3526,-3661){\makebox(0,0)[lb]{\smash{\SetFigFont{9}{10.8}{\rmdefault}{\mddefault}{\updefault}{\color[rgb]{0,0,0}$\xh_2'$}%
}}}
\put(3151,-3811){\makebox(0,0)[lb]{\smash{\SetFigFont{9}{10.8}{\rmdefault}{\mddefault}{\updefault}{\color[rgb]{0,0,0}$\xh_2$}%
}}}
\end{picture}

%% file: handles.pstex_t
\begin{picture}(0,0)%
\includegraphics{handles.pstex}%
\end{picture}%
\setlength{\unitlength}{3947sp}%
\begingroup\makeatletter\ifx\SetFigFont\undefined%
\gdef\SetFigFont#1#2#3#4#5{%
  \reset@font\fontsize{#1}{#2pt}%
  \fontfamily{#3}\fontseries{#4}\fontshape{#5}%
  \selectfont}%
\fi\endgroup%
\begin{picture}(5142,5329)(3130,-5444)
\put(3226,-2836){\makebox(0,0)[lb]{\smash{\SetFigFont{12}{14.4}{\rmdefault}{\mddefault}{\updefault}{\color[rgb]{0,0,0}$\mu_1$}%
}}}
\put(5776,-2836){\makebox(0,0)[lb]{\smash{\SetFigFont{12}{14.4}{\rmdefault}{\mddefault}{\updefault}{\color[rgb]{0,0,0}$\mu_4$}%
}}}
\put(6826,-2836){\makebox(0,0)[lb]{\smash{\SetFigFont{12}{14.4}{\rmdefault}{\mddefault}{\updefault}{\color[rgb]{0,0,0}$\mu_{2m-1}$}%
}}}
\put(7876,-2836){\makebox(0,0)[lb]{\smash{\SetFigFont{12}{14.4}{\rmdefault}{\mddefault}{\updefault}{\color[rgb]{0,0,0}$\mu_{2m}$}%
}}}
\put(7351,-2611){\makebox(0,0)[lb]{\smash{\SetFigFont{12}{14.4}{\rmdefault}{\mddefault}{\updefault}{\color[rgb]{0,0,0}$\alpha_{m}$}%
}}}
\put(5251,-2611){\makebox(0,0)[lb]{\smash{\SetFigFont{12}{14.4}{\rmdefault}{\mddefault}{\updefault}{\color[rgb]{0,0,0}$\alpha_2$}%
}}}
\put(4201,-2836){\makebox(0,0)[lb]{\smash{\SetFigFont{12}{14.4}{\rmdefault}{\mddefault}{\updefault}{\color[rgb]{0,0,0}$\mu_2$}%
}}}
\put(4876,-2836){\makebox(0,0)[lb]{\smash{\SetFigFont{12}{14.4}{\rmdefault}{\mddefault}{\updefault}{\color[rgb]{0,0,0}$\mu_3$}%
}}}
\put(3751,-2611){\makebox(0,0)[lb]{\smash{\SetFigFont{12}{14.4}{\rmdefault}{\mddefault}{\updefault}{\color[rgb]{0,0,0}$\alpha_1$}%
}}}
\put(7126,-3736){\makebox(0,0)[lb]{\smash{\SetFigFont{12}{14.4}{\rmdefault}{\mddefault}{\updefault}{\color[rgb]{0,0,0}$h_{m}$}%
}}}
\put(3601,-4636){\makebox(0,0)[lb]{\smash{\SetFigFont{12}{14.4}{\rmdefault}{\mddefault}{\updefault}{\color[rgb]{0,0,0}$\eta_1$}%
}}}
\put(5176,-5386){\makebox(0,0)[lb]{\smash{\SetFigFont{12}{14.4}{\rmdefault}{\mddefault}{\updefault}{\color[rgb]{0,0,0}$\eta_2$}%
}}}
\put(7726,-1036){\makebox(0,0)[lb]{\smash{\SetFigFont{12}{14.4}{\rmdefault}{\mddefault}{\updefault}{\color[rgb]{0,0,0}$D$}%
}}}
\put(7501,-4711){\makebox(0,0)[lb]{\smash{\SetFigFont{12}{14.4}{\rmdefault}{\mddefault}{\updefault}{\color[rgb]{0,0,0}$\eta_{m}$}%
}}}
\put(3976,-3736){\makebox(0,0)[lb]{\smash{\SetFigFont{12}{14.4}{\rmdefault}{\mddefault}{\updefault}{\color[rgb]{0,0,0}$h_1$}%
}}}
\put(5476,-3736){\makebox(0,0)[lb]{\smash{\SetFigFont{12}{14.4}{\rmdefault}{\mddefault}{\updefault}{\color[rgb]{0,0,0}$h_2$}%
}}}
\end{picture}

%% file: flat2.pstex_t
\begin{picture}(0,0)%
\includegraphics{flat2.pstex}%
\end{picture}%
\setlength{\unitlength}{3947sp}%
\begingroup\makeatletter\ifx\SetFigFont\undefined%
\gdef\SetFigFont#1#2#3#4#5{%
  \reset@font\fontsize{#1}{#2pt}%
  \fontfamily{#3}\fontseries{#4}\fontshape{#5}%
  \selectfont}%
\fi\endgroup%
\begin{picture}(5419,4084)(432,-4448)
\put(2026,-3211){\makebox(0,0)[lb]{\smash{\SetFigFont{10}{12.0}{\rmdefault}{\mddefault}{\updefault}{\color[rgb]{0,0,0}$\bar \uh_1$}%
}}}
\put(3226,-3136){\makebox(0,0)[lb]{\smash{\SetFigFont{10}{12.0}{\familydefault}{\mddefault}{\updefault}{\color[rgb]{0,0,0}$ \bar \uh_2$}%
}}}
\put(4152,-3117){\makebox(0,0)[lb]{\smash{\SetFigFont{10}{12.0}{\familydefault}{\mddefault}{\updefault}{\color[rgb]{0,0,0}$\bar \yh_1$}%
}}}
\put(4726,-3211){\makebox(0,0)[lb]{\smash{\SetFigFont{10}{12.0}{\familydefault}{\mddefault}{\updefault}{\color[rgb]{0,0,0}$\bar \vh_1$}%
}}}
\put(5401,-3211){\makebox(0,0)[lb]{\smash{\SetFigFont{10}{12.0}{\familydefault}{\mddefault}{\updefault}{\color[rgb]{0,0,0}$\bar \yh_2$}%
}}}
\put(5851,-3061){\makebox(0,0)[lb]{\smash{\SetFigFont{10}{12.0}{\familydefault}{\mddefault}{\updefault}{\color[rgb]{0,0,0}$\bar \vh_2$}%
}}}
\put(1576,-3136){\makebox(0,0)[lb]{\smash{\SetFigFont{10}{12.0}{\rmdefault}{\mddefault}{\updefault}{\color[rgb]{0,0,0}$\bar \xh_1$}%
}}}
\put(2551,-3211){\makebox(0,0)[lb]{\smash{\SetFigFont{10}{12.0}{\rmdefault}{\mddefault}{\updefault}{\color[rgb]{0,0,0}$\bar \xh_2$}%
}}}
\end{picture}

%% file: loop.pstex_t
\begin{picture}(0,0)%
\includegraphics{loop.pstex}%
\end{picture}%
\setlength{\unitlength}{3947sp}%
\begingroup\makeatletter\ifx\SetFigFont\undefined%
\gdef\SetFigFont#1#2#3#4#5{%
  \reset@font\fontsize{#1}{#2pt}%
  \fontfamily{#3}\fontseries{#4}\fontshape{#5}%
  \selectfont}%
\fi\endgroup%
\begin{picture}(5282,5313)(1117,-5035)
\put(2326,-4861){\makebox(0,0)[lb]{\smash{\SetFigFont{12}{14.4}{\rmdefault}{\mddefault}{\updefault}{\color[rgb]{0,0,0}$\eta_1$}%
}}}
\put(4951,-4861){\makebox(0,0)[lb]{\smash{\SetFigFont{12}{14.4}{\rmdefault}{\mddefault}{\updefault}{\color[rgb]{0,0,0}$\eta_2$}%
}}}
\put(5603,-264){\makebox(0,0)[lb]{\smash{\SetFigFont{12}{14.4}{\rmdefault}{\mddefault}{\updefault}{\color[rgb]{0,0,0}$D$}%
}}}
\put(1951,-2236){\makebox(0,0)[lb]{\smash{\SetFigFont{12}{14.4}{\rmdefault}{\mddefault}{\updefault}{\color[rgb]{0,0,0}$\xh_1$}%
}}}
\put(5551,-2461){\makebox(0,0)[lb]{\smash{\SetFigFont{12}{14.4}{\rmdefault}{\mddefault}{\updefault}{\color[rgb]{0,0,0}$\vh_2 = \nu_2 $}%
}}}
\end{picture}

%% file: big.pstex_t
\begin{picture}(0,0)%
\includegraphics{big.pstex}%
\end{picture}%
\setlength{\unitlength}{3947sp}%
\begingroup\makeatletter\ifx\SetFigFont\undefined%
\gdef\SetFigFont#1#2#3#4#5{%
  \reset@font\fontsize{#1}{#2pt}%
  \fontfamily{#3}\fontseries{#4}\fontshape{#5}%
  \selectfont}%
\fi\endgroup%
\begin{picture}(6762,3586)(526,-3635)
\put(526,-3136){\makebox(0,0)[lb]{\smash{\SetFigFont{12}{14.4}{\rmdefault}{\mddefault}{\updefault}{\color[rgb]{0,0,0}$-3$}%
}}}
\put(526,-2536){\makebox(0,0)[lb]{\smash{\SetFigFont{12}{14.4}{\rmdefault}{\mddefault}{\updefault}{\color[rgb]{0,0,0}$-2$}%
}}}
\put(526,-1936){\makebox(0,0)[lb]{\smash{\SetFigFont{12}{14.4}{\rmdefault}{\mddefault}{\updefault}{\color[rgb]{0,0,0}$-1$}%
}}}
\put(676,-1336){\makebox(0,0)[lb]{\smash{\SetFigFont{12}{14.4}{\rmdefault}{\mddefault}{\updefault}{\color[rgb]{0,0,0}$0$}%
}}}
\put(1051,-3586){\makebox(0,0)[lb]{\smash{\SetFigFont{12}{14.4}{\rmdefault}{\mddefault}{\updefault}{\color[rgb]{0,0,0}$-9$}%
}}}
\put(1651,-3586){\makebox(0,0)[lb]{\smash{\SetFigFont{12}{14.4}{\rmdefault}{\mddefault}{\updefault}{\color[rgb]{0,0,0}$-8$}%
}}}
\put(2251,-3586){\makebox(0,0)[lb]{\smash{\SetFigFont{12}{14.4}{\rmdefault}{\mddefault}{\updefault}{\color[rgb]{0,0,0}$-7$}%
}}}
\put(2851,-3586){\makebox(0,0)[lb]{\smash{\SetFigFont{12}{14.4}{\rmdefault}{\mddefault}{\updefault}{\color[rgb]{0,0,0}$-6$}%
}}}
\put(3451,-3586){\makebox(0,0)[lb]{\smash{\SetFigFont{12}{14.4}{\rmdefault}{\mddefault}{\updefault}{\color[rgb]{0,0,0}$-5$}%
}}}
\put(4051,-3586){\makebox(0,0)[lb]{\smash{\SetFigFont{12}{14.4}{\rmdefault}{\mddefault}{\updefault}{\color[rgb]{0,0,0}$-4$}%
}}}
\put(4651,-3586){\makebox(0,0)[lb]{\smash{\SetFigFont{12}{14.4}{\rmdefault}{\mddefault}{\updefault}{\color[rgb]{0,0,0}$-3$}%
}}}
\put(5251,-3586){\makebox(0,0)[lb]{\smash{\SetFigFont{12}{14.4}{\rmdefault}{\mddefault}{\updefault}{\color[rgb]{0,0,0}$-2$}%
}}}
\put(5851,-3586){\makebox(0,0)[lb]{\smash{\SetFigFont{12}{14.4}{\rmdefault}{\mddefault}{\updefault}{\color[rgb]{0,0,0}$-1$}%
}}}
\put(7276,-3586){\makebox(0,0)[lb]{\smash{\SetFigFont{12}{14.4}{\rmdefault}{\mddefault}{\updefault}{\color[rgb]{0,0,0}$J$}%
}}}
\put(976,-211){\makebox(0,0)[lb]{\smash{\SetFigFont{12}{14.4}{\rmdefault}{\mddefault}{\updefault}{\color[rgb]{0,0,0}$P+J$}%
}}}
\end{picture}

%% file: khov.pstex_t
\begin{picture}(0,0)%
\includegraphics{khov.pstex}%
\end{picture}%
\setlength{\unitlength}{3947sp}%
\begingroup\makeatletter\ifx\SetFigFont\undefined%
\gdef\SetFigFont#1#2#3#4#5{%
  \reset@font\fontsize{#1}{#2pt}%
  \fontfamily{#3}\fontseries{#4}\fontshape{#5}%
  \selectfont}%
\fi\endgroup%
\begin{picture}(6762,3715)(526,-3644)
\put(526,-3136){\makebox(0,0)[lb]{\smash{\SetFigFont{12}{14.4}{\rmdefault}{\mddefault}{\updefault}{\color[rgb]{0,0,0}$-3$}%
}}}
\put(526,-2536){\makebox(0,0)[lb]{\smash{\SetFigFont{12}{14.4}{\rmdefault}{\mddefault}{\updefault}{\color[rgb]{0,0,0}$-2$}%
}}}
\put(526,-1936){\makebox(0,0)[lb]{\smash{\SetFigFont{12}{14.4}{\rmdefault}{\mddefault}{\updefault}{\color[rgb]{0,0,0}$-1$}%
}}}
\put(676,-1336){\makebox(0,0)[lb]{\smash{\SetFigFont{12}{14.4}{\rmdefault}{\mddefault}{\updefault}{\color[rgb]{0,0,0}$0$}%
}}}
\put(1051,-3586){\makebox(0,0)[lb]{\smash{\SetFigFont{12}{14.4}{\rmdefault}{\mddefault}{\updefault}{\color[rgb]{0,0,0}$-9$}%
}}}
\put(1651,-3586){\makebox(0,0)[lb]{\smash{\SetFigFont{12}{14.4}{\rmdefault}{\mddefault}{\updefault}{\color[rgb]{0,0,0}$-8$}%
}}}
\put(2251,-3586){\makebox(0,0)[lb]{\smash{\SetFigFont{12}{14.4}{\rmdefault}{\mddefault}{\updefault}{\color[rgb]{0,0,0}$-7$}%
}}}
\put(2851,-3586){\makebox(0,0)[lb]{\smash{\SetFigFont{12}{14.4}{\rmdefault}{\mddefault}{\updefault}{\color[rgb]{0,0,0}$-6$}%
}}}
\put(3451,-3586){\makebox(0,0)[lb]{\smash{\SetFigFont{12}{14.4}{\rmdefault}{\mddefault}{\updefault}{\color[rgb]{0,0,0}$-5$}%
}}}
\put(4051,-3586){\makebox(0,0)[lb]{\smash{\SetFigFont{12}{14.4}{\rmdefault}{\mddefault}{\updefault}{\color[rgb]{0,0,0}$-4$}%
}}}
\put(4651,-3586){\makebox(0,0)[lb]{\smash{\SetFigFont{12}{14.4}{\rmdefault}{\mddefault}{\updefault}{\color[rgb]{0,0,0}$-3$}%
}}}
\put(5251,-3586){\makebox(0,0)[lb]{\smash{\SetFigFont{12}{14.4}{\rmdefault}{\mddefault}{\updefault}{\color[rgb]{0,0,0}$-2$}%
}}}
\put(5851,-3586){\makebox(0,0)[lb]{\smash{\SetFigFont{12}{14.4}{\rmdefault}{\mddefault}{\updefault}{\color[rgb]{0,0,0}$-1$}%
}}}
\put(676,-61){\makebox(0,0)[lb]{\smash{\SetFigFont{12}{14.4}{\rmdefault}{\mddefault}{\updefault}{\color[rgb]{0,0,0}$i$}%
}}}
\put(7276,-3586){\makebox(0,0)[lb]{\smash{\SetFigFont{12}{14.4}{\rmdefault}{\mddefault}{\updefault}{\color[rgb]{0,0,0}$j$}%
}}}
\put(5926,-1336){\makebox(0,0)[lb]{\smash{\SetFigFont{12}{14.4}{\rmdefault}{\mddefault}{\updefault}{\color[rgb]{0,0,0}$\zz$}%
}}}
\put(2251,-2536){\makebox(0,0)[lb]{\smash{\SetFigFont{12}{14.4}{\rmdefault}{\mddefault}{\updefault}{\color[rgb]{0,0,0}$\zz/2$}%
}}}
\put(1126,-3136){\makebox(0,0)[lb]{\smash{\SetFigFont{12}{14.4}{\rmdefault}{\mddefault}{\updefault}{\color[rgb]{0,0,0}$\zz$}%
}}}
\put(4726,-1336){\makebox(0,0)[lb]{\smash{\SetFigFont{12}{14.4}{\rmdefault}{\mddefault}{\updefault}{\color[rgb]{0,0,0}$\zz$}%
}}}
\put(3526,-2536){\makebox(0,0)[lb]{\smash{\SetFigFont{12}{14.4}{\rmdefault}{\mddefault}{\updefault}{\color[rgb]{0,0,0}$\zz$}%
}}}
\end{picture}

%% file: hbody.pstex_t
\begin{picture}(0,0)%
\includegraphics{hbody.pstex}%
\end{picture}%
\setlength{\unitlength}{3947sp}%
\begingroup\makeatletter\ifx\SetFigFont\undefined%
\gdef\SetFigFont#1#2#3#4#5{%
  \reset@font\fontsize{#1}{#2pt}%
  \fontfamily{#3}\fontseries{#4}\fontshape{#5}%
  \selectfont}%
\fi\endgroup%
\begin{picture}(6474,2506)(889,-1914)
\put(5626,314){\makebox(0,0)[lb]{\smash{\SetFigFont{12}{14.4}{\rmdefault}{\mddefault}{\updefault}{\color[rgb]{0,0,0}$+\infty$}%
}}}
\put(1126,-1036){\makebox(0,0)[lb]{\smash{\SetFigFont{12}{14.4}{\rmdefault}{\mddefault}{\updefault}{\color[rgb]{0,0,0}$\hat\alpha_1$}%
}}}
\put(2626,-1036){\makebox(0,0)[lb]{\smash{\SetFigFont{12}{14.4}{\rmdefault}{\mddefault}{\updefault}{\color[rgb]{0,0,0}$\hat \alpha_2$}%
}}}
\put(4276,-1036){\makebox(0,0)[lb]{\smash{\SetFigFont{12}{14.4}{\rmdefault}{\mddefault}{\updefault}{\color[rgb]{0,0,0}$\hat\alpha_3$}%
}}}
\put(5926,-1036){\makebox(0,0)[lb]{\smash{\SetFigFont{12}{14.4}{\rmdefault}{\mddefault}{\updefault}{\color[rgb]{0,0,0}$\hat\alpha_4$}%
}}}
\put(5626,-1786){\makebox(0,0)[lb]{\smash{\SetFigFont{12}{14.4}{\rmdefault}{\mddefault}{\updefault}{\color[rgb]{0,0,0}$-\infty$}%
}}}
\end{picture}

%% file: disks.pstex_t
\begin{picture}(0,0)%
\includegraphics{disks.pstex}%
\end{picture}%
\setlength{\unitlength}{3947sp}%
\begingroup\makeatletter\ifx\SetFigFont\undefined%
\gdef\SetFigFont#1#2#3#4#5{%
  \reset@font\fontsize{#1}{#2pt}%
  \fontfamily{#3}\fontseries{#4}\fontshape{#5}%
  \selectfont}%
\fi\endgroup%
\begin{picture}(6024,1668)(1189,-1573)
\put(1876,-886){\makebox(0,0)[lb]{\smash{\SetFigFont{12}{14.4}{\rmdefault}{\mddefault}{\updefault}{\color[rgb]{0,0,0}$\mu_1$}%
}}}
\put(4426,-886){\makebox(0,0)[lb]{\smash{\SetFigFont{12}{14.4}{\rmdefault}{\mddefault}{\updefault}{\color[rgb]{0,0,0}$\mu_4$}%
}}}
\put(5476,-886){\makebox(0,0)[lb]{\smash{\SetFigFont{12}{14.4}{\rmdefault}{\mddefault}{\updefault}{\color[rgb]{0,0,0}$\mu_{2m-1}$}%
}}}
\put(6526,-886){\makebox(0,0)[lb]{\smash{\SetFigFont{12}{14.4}{\rmdefault}{\mddefault}{\updefault}{\color[rgb]{0,0,0}$\mu_{2m}$}%
}}}
\put(2851,-886){\makebox(0,0)[lb]{\smash{\SetFigFont{12}{14.4}{\rmdefault}{\mddefault}{\updefault}{\color[rgb]{0,0,0}$\mu_2$}%
}}}
\put(3526,-886){\makebox(0,0)[lb]{\smash{\SetFigFont{12}{14.4}{\rmdefault}{\mddefault}{\updefault}{\color[rgb]{0,0,0}$\mu_3$}%
}}}
\put(1426,-1486){\makebox(0,0)[lb]{\smash{\SetFigFont{12}{14.4}{\rmdefault}{\mddefault}{\updefault}{\color[rgb]{0,0,0}$S^2$}%
}}}
\put(3901,-511){\makebox(0,0)[lb]{\smash{\SetFigFont{12}{14.4}{\rmdefault}{\mddefault}{\updefault}{\color[rgb]{0,0,0}$F_2$}%
}}}
\put(6001,-511){\makebox(0,0)[lb]{\smash{\SetFigFont{12}{14.4}{\rmdefault}{\mddefault}{\updefault}{\color[rgb]{0,0,0}$F_m$}%
}}}
\put(2326,-61){\makebox(0,0)[lb]{\smash{\SetFigFont{12}{14.4}{\rmdefault}{\mddefault}{\updefault}{\color[rgb]{0,0,0}$\alpha_1^{\dag}$}%
}}}
\put(3826,-61){\makebox(0,0)[lb]{\smash{\SetFigFont{12}{14.4}{\rmdefault}{\mddefault}{\updefault}{\color[rgb]{0,0,0}$\alpha_2^{\dag}$}%
}}}
\put(5926,-61){\makebox(0,0)[lb]{\smash{\SetFigFont{12}{14.4}{\rmdefault}{\mddefault}{\updefault}{\color[rgb]{0,0,0}$\alpha_{m}^{\dag}$}%
}}}
\put(2326,-886){\makebox(0,0)[lb]{\smash{\SetFigFont{12}{14.4}{\rmdefault}{\mddefault}{\updefault}{\color[rgb]{0,0,0}$\alpha_1$}%
}}}
\put(3901,-886){\makebox(0,0)[lb]{\smash{\SetFigFont{12}{14.4}{\rmdefault}{\mddefault}{\updefault}{\color[rgb]{0,0,0}$\alpha_2$}%
}}}
\put(6001,-886){\makebox(0,0)[lb]{\smash{\SetFigFont{12}{14.4}{\rmdefault}{\mddefault}{\updefault}{\color[rgb]{0,0,0}$\alpha_{m}$}%
}}}
\put(2401,-511){\makebox(0,0)[lb]{\smash{\SetFigFont{12}{14.4}{\rmdefault}{\mddefault}{\updefault}{\color[rgb]{0,0,0}$F_1$}%
}}}
\end{picture}

%% file: hs.bbl
\begin {thebibliography}{99999}

\bibitem {B}
S. Bigelow,  {\it A homological definition of the Jones polynomial,} {in 
Invariants of knots and 3-manifolds (Kyoto, 2001), 29-41,  
Geom. Topol. Monogr. {\bf 4}, Geom. Topol. Publ., Coventry, 2002.}

\bibitem {Be}
A. Beauville, {\it Vari\'et\'es k\"ahleriennes dont la premi\`ere classe 
de Chern est nulle,} J. Diff. Geom. {\bf 18} {(1983), 755-782.}

\bibitem {F}
A. Floer, {\it Morse theory for Lagrangian intersections, } {J. Diff. 
Geom.} {\bf 28} {(1988), 513-547.}

\bibitem {Fo}
J. Fogarty, {\it Algebraic families on an algebraic surface,} {Amer. J. 
Math.} {\bf 90} {(1968), 511-521.}

\bibitem {FOOO}
K. Fukaya, Y.-G. Oh, H. Ohta, and K. Ono, {\it Lagrangian intersection Floer theory - anomaly and 
obstruction,} preprint (2000), available at \url{www.math.kyoto-u.ac.jp/~fukaya/fukaya.html}

\bibitem {G}
A. Grothendieck, {\it Techniques de construction et the\'eoremes
d'existence en g\'eom\'etrie algebrique, IV: Les sch\'emas de Hilbert,}
{S\'em. Bourbaki} {\bf 221} {(1960/61).}

\bibitem {J}
V. Jones,  {\it A polynomial invariant for knots via von Neumann 
algebras,} {Bull. Amer. Math. Soc.} {\bf 12} {(1985), 103-111.}

\bibitem {KL}
S. Katz and C.-C. M. Liu, {\it Enumerative geometry of stable maps with Lagrangian boundary 
conditions and multiple covers of the disc,} {Adv. Theor. Math. Phys.} {\bf 5}  {(2001)  no.1, 
1-49.}

\bibitem {Kh}
M. Khovanov, {\it A categorification of the Jones polynomial,} 
{Duke Math. J.} {\bf 101} {(2000) no. 3, 359-426.}

\bibitem {Kh2} 
M. Khovanov, {\it A functor-valued invariant of tangles,}  Algebr. Geom. 
Topol. {\bf 2} {(2002), 665-741.}

\bibitem {Kh3}
M. Khovanov, {\it  Patterns in knot cohomology I,} Experiment. Math. {\bf 
12} {(2003), no. 3, 365-374.}

\bibitem {KS}
M. Khovanov and P. Seidel, {\it Quivers, Floer cohomology, and braid group 
actions,} {J. Amer. Math. Soc.} {\bf 15} {(2002), 203-271.}

\bibitem {Ko}
M. Kontsevich, {\it Homological algebra of mirror symmetry,} in 
Proceedings of the International Congress of Mathematicians (Z\"urich, 
1994), p. 120-139, Birkh\"auser, 1995.

\bibitem {Kr}
P. B. Kronheimer, {\it The construction of ALE spaces as hyper-K\"ahler 
quotients,} { J. Diff. Geom.} {\bf 29} {(1989), 665-683.}

\bibitem {Kr2}
P. B. Kronheimer, {\it Instantons and the geometry of the nilpotent 
variety,} { J. Diff. Geom.} {\bf 32} {(1990), 473-490.}

\bibitem {KN}
P. B. Kronheimer and H. Nakajima, {\it Yang-Mills instantons on ALE 
gravitational instantons,} {Math. Ann.} {\bf 288} {(1990), 263-307.}

\bibitem {L}
R. J. Lawrence, {\it Homological representations of the Hecke algebra,} 
{Comm. Math. Phys.} {\bf 135} {(1990), no. 1, 141-191.}

\bibitem {N1}
H. Nakajima, {\it Instantons on ALE spaces, quiver varieties, and 
Kac-Moody algebras,} {Duke Math. J.} {\bf 76} {(1994) no. 2, 365-416.}

\bibitem {N2}
H. Nakajima, {\it Lectures on Hilbert schemes of points on surfaces,} 
{Univ. Lect. Series} {\bf 18,} {Amer. Math. Soc, 1999.}

\bibitem {O1}
Y.-G. Oh, {\it Riemann-Hilbert problem and application to the perturbation theory of analytic 
discs,} {Kyungpook Math. J.} {\bf 35} {(1995) no. 1, 39-75.}

\bibitem {O2}
Y.-G. Oh, {\it On the structure of pseudo-holomorphic discs with totally real boundary conditions, }
{J. Geom. Anal.} {\bf 7} {(1997) no. 2, 305-327.}

\bibitem {OS1}
P. Ozsv\'ath and Z. Szab\'o, {\it Holomorphic disks and topological 
invariants for closed three-manifolds,} {Annals of Math. } {\bf 159} {(2004), no. 3, 
1027-1158.}

\bibitem {OS2}
P. Ozsv\'ath and Z. Szab\'o, {\it Holomorphic disks and three-manifold 
invariants: properties and applications,} {Annals of Math.} {\bf 159} {(2004), no. 3,  
1159-1245.}

\bibitem {OS3}
P. Ozsv\'ath and Z. Szab\'o, {\it On the Heegaard Floer homology of 
branched double-covers,} {Adv. Math.} {\bf 194} {(2005), 1-33.}

\bibitem {P}
M. Pozniak, {\it Floer homology, Novikov rings and clean intersections,} 
in {Northern California Symplectic Geometry Seminar, p. 119-181, Amer. 
Math. Soc., 1999.}

\bibitem {RS}
J. Robbin and D. Salamon, {\it The Maslov index for paths,} {Topology } {\bf 32} {(1993) no. 4, 
827-844.}

\bibitem {QW}
Z. Qin and W. Wang, {\it Hilbert schemes of points on the minimal 
resolution and soliton equations,} {preprint, math.QA/0404540.}   

\bibitem {R}
J. Rasmussen, {\it Khovanov homology and the slice genus,} {preprint, 
math.GT/0402131.}

\bibitem {Se}
P. Seidel, {\it Graded Lagrangian submanifolds,} {Bull. Soc. Math. 
France} {\bf 128} {(2000), 103-146.}

\bibitem {SS}
P. Seidel and I. Smith, {\it A link invariant from the symplectic 
geometry of nilpotent slices,} {preprint, math.SG/0405089.}

\bibitem {S}
I. Smith, {\it Symplectic geometry of the adjoint quotient I (joint with 
Paul Seidel),} lectures at MSRI in Spring 2004, 
\url {www.msri.org/publications/ln/msri/2004/symplecticgeom/smith/1/index.html}

\bibitem {W}
W. Wang, {\it Hilbert schemes, wreath products, and the McKay 
correspondence,} preprint, math.AG/9912104. 

\bibitem {We}
A. Weinstein, {\it Lagrangian submanifolds and hamiltonian systems,} 
{Annals of Math.} {\bf 98} {(1973), 377-410.} 

\end{thebibliography}
